\documentclass[12pt]{amsart}
\date{14 August 2012}
\usepackage{amscd,amssymb}
\usepackage{framed}
\usepackage{mathdots}
\usepackage{hyperref}
\input xy
\xyoption{all}

\let\oldmarginpar\marginpar
\renewcommand\marginpar[1]{\oldmarginpar{\tiny\bf\begin{flushleft} #1
\end{flushleft}}}

\numberwithin{table}{section}

\newcommand{\la}{\langle}
\newcommand{\ra}{\rangle}

\newcommand{\xra}{\xrightarrow}
\theoremstyle{plain}  
\newtheorem{theorem}{Theorem}[section]

\newtheorem*{theorem*}{Theorem}
\newtheorem{corollary}[theorem]{Corollary}
\newtheorem{lemma}[theorem]{Lemma}
\newtheorem{proposition}[theorem]{Proposition}

\theoremstyle{definition}
\newtheorem{definition}[theorem]{Definition}
\newtheorem{notation}[theorem]{Notation}

\theoremstyle{remark}

\newtheorem{remark}[theorem]{Remark}
\newtheorem*{remark*}{Remark}

\newtheorem*{claim*}{Claim}

\newcommand{\curly}{\mathcal}

\newcommand{\bE}{{\bf{E}}}
\newcommand{\dbar}{\bar{\partial}}
\newcommand{\cR}{\mathcal{R}}
\newcommand{\cM}{\mathcal{M}}
\renewcommand{\AA}{\mathbb{A}}

\newcommand{\CC}{\mathbb{C}}

\newcommand{\EE}{\mathbb{E}}

\newcommand{\HH}{\mathbb{H}}

\newcommand{\RR}{\mathbb{R}}
\newcommand{\VV}{\mathbb{V}}
\newcommand{\WW}{\mathbb{W}}
\newcommand{\ZZ}{\mathbb{Z}}
\newcommand{\lie}{\mathfrak}

\newcommand{\alie}{\mathfrak{a}}

\newcommand{\clie}{\mathfrak{c}}
\newcommand{\glie}{\mathfrak{g}}
\newcommand{\hlie}{\mathfrak{h}}

\newcommand{\llie}{\mathfrak{l}}
\newcommand{\mlie}{\mathfrak{m}}
\newcommand{\olie}{\mathfrak{o}}
\newcommand{\plie}{\mathfrak{p}}

\newcommand{\tlie}{\mathfrak{t}}
\newcommand{\ulie}{\mathfrak{u}}
\newcommand{\zlie}{\mathfrak{z}}
\newcommand{\sllie}{\mathfrak{sl}}
\newcommand{\splie}{\mathfrak{sp}}

\newcommand{\gllie}{\mathfrak{gl}}

\newcommand{\lieg}{\mathfrak{g}}
\newcommand{\lieh}{\mathfrak{h}}
\newcommand{\liehc}{\mathfrak{h}^{\CC}}

\newcommand{\liem}{\mathfrak{m}}

\newcommand{\liep}{\mathfrak{p}}

\newcommand{\lieu}{\mathfrak{u}}
\newcommand{\liez}{\mathfrak{z}}

\newcommand{\Ad}{\operatorname{Ad}}
\newcommand{\aut}{\operatorname{aut}}
\newcommand{\Aut}{\operatorname{Aut}}

\newcommand{\CH}{\operatorname{CH}}

\newcommand{\End}{\operatorname{End}}

\newcommand{\GL}{\operatorname{GL}}

\newcommand{\Gr}{\operatorname{Gr}}

\newcommand{\Hom}{\operatorname{Hom}}
\newcommand{\Id}{\operatorname{Id}}

\newcommand{\Ker}{\operatorname{Ker}}

\newcommand{\Lie}{\operatorname{Lie}}

\renewcommand{\O}{\operatorname{O}}

\newcommand{\rk}{\operatorname{rk}}

\newcommand{\SL}{\operatorname{SL}}
\newcommand{\SO}{\operatorname{SO}}

\newcommand{\Sp}{\operatorname{Sp}}

\newcommand{\SU}{\operatorname{SU}}

\newcommand{\Tr}{\operatorname{Tr}}

\newcommand{\U}{\operatorname{U}}

\renewcommand{\exp}{\operatorname{exp}}

\renewcommand{\H}{\operatorname{H}}
\newcommand{\ov}{\overline}

\newcommand{\Mg}{\mathcal{M}^{\operatorname{gauge}}}

\newcommand{\AAA}{{\curly A}}
\newcommand{\CCC}{{\curly C}}

\newcommand{\EEE}{{\curly E}}
\newcommand{\FFF}{{\curly F}}
\newcommand{\GGG}{{\curly G}}
\newcommand{\HHH}{{\curly H}}
\newcommand{\III}{\mathcal{I}}
\newcommand{\JJJ}{\mathcal{J}}

\newcommand{\OOO}{{\curly O}}
\newcommand{\PPP}{{\curly P}}

\newcommand{\SSS}{{\curly S}}

\newcommand{\VVV}{{\curly V}}
\newcommand{\WWW}{\mathcal{W}}
\newcommand{\qu}{/\kern-.7ex/}
\newcommand{\exh}{\to\kern-1.8ex\to}

\newcommand{\VP}{{\curly V}\kern-0.9ex\PPP}

\newcommand{\imag}{{\mathbf i}}

\newcommand{\suchthat}{\;|\;}

\newcommand{\lto}{\longrightarrow}
\newcommand{\into}{\hookrightarrow}
\newcommand{\ad}{\mathrm{ad}}

\newcommand{\Kto}{\ar@{*+=[o][F]{\scriptscriptstyle{K}}->}[r]}

\newcommand{\HC}{H^{\CC}}
\newcommand{\hclie}{\hlie^{\CC}}
\newcommand{\mclie}{\mlie^{\CC}}

\newcommand{\ddt}{\dfrac{d}{dt}}

\title[The Hitchin--Kobayashi correspondence]
{The Hitchin--Kobayashi correspondence, Higgs pairs and surface
group representations}

\author{O. Garc\'{\i}a-Prada, P.~B. Gothen, I. Mundet i Riera}



\thanks{
  Members of the Research Group VBAC (Vector Bundles on Algebraic
  Curves).
  Research partially supported by Ministerio de Educaci\'{o}n y
  Ciencia, CSIC, Conselho de Reitores das Universidades Portuguesas
  and FCT (Portugal) through Spain--Portugal bilateral research
  projects.
  First and Third authors partially supported by Ministerio de
  Educaci{\'o}n y Ciencia (Spain) through Projects
  MTM2004-07090-C03-01 and MTM2007-67623.
  Second author partially supported by the FCT (Portugal) with EU
  (FEDER/COMPETE) and national funds through the projects
  PTDC/MAT/099275/2008 and PTDC/MAT/098770/2008, and through Centro de
  Matem\'atica da Universidade do Porto (PEst-C/MAT/UI0144/2011).
}

\begin{document}

\begin{abstract}
  We develop a complete Hitchin--Kobayashi correspondence for twisted
  pairs on a compact Riemann surface $X$. Here a twisted pair consists
  of a principal bundle and a section of an associated vector bundle
  twisted by a fixed line bundle. The main novelty lies in a careful
  study of the the notion of polystability, required for having a
  bijective correspondence between solutions to the Hermite--Einstein
  equations, on one hand, and polystable pairs, on the other. Our
  results allow us to establish rigorously the identification between
  the moduli space of polystable $G$-Higgs bundles on $X$ and the
  character variety for representations of the fundamental group of
  $X$ in $G$ for any connected real reductive Lie group $G$. We also
  study in detail several interesting examples of the correspondence
  for particular groups and show how to significantly simplify the
  general stability condition in these cases.
\end{abstract}

\maketitle



\section{Introduction}

In this paper we study the Hitchin--Kobayashi correspondece for
$L$-twisted pairs on a compact Riemann surface $X$. The main
motivation for our study comes from non-abelian Hodge theory on $X$
for a real semisimple Lie group $G$. Our resuls allow us to establish
a one-to-one correspondence between the moduli space of $G$-Higgs
bundles over $X$ and the moduli space of reductive representations of
the fundamental group of $X$ in $G$.

The non-abelian Hodge theory correspondence has two fundamental
ingredients: one ingredient is the Theorem of Corlette
\cite{corlette:1988} and Donaldson \cite{donaldson:1987} on the
existence of harmonic metrics in flat bundles, and the other grows out
of the Hitchin--Kobayashi correspondence between polystable Higgs
bundles and solutions to Hitchin's gauge theoretic equations,
established by Hitchin \cite{hitchin:1987a} and Simpson
\cite{simpson:1988,simpson:1992,simpson:1994,simpson:1995}. While the
Corlette--Donaldson Theorem applies directly in our context, for the
Hitchin--Kobayashi we need to work in the general setting of stable
pairs treated in \cite{banfield,bradlow-garcia-prada-mundet:2003}.
One of the  main contributions of the present paper is to establish the
extension of this general correspondence to strictly polystable
pairs. This is required for having a complete correspondence with
solutions to the gauge theoretic equations and is essential for the
application of the theory to moduli of representations of surface
groups. The other main contribution lies in a careful study of the
general stability condition in several important special cases.
This leads to a simplification of the stability condition which makes
it practical for applications of the theory.

We describe now briefly the content of the different sections of the
paper.

In order to establish the full Hitchin--Kobayashi correspondence,
in Section~\ref{appendix} we review the general theory of
$L$-twisted pairs and the Hitchin--Kobayashi correspondence over a
compact Riemann surface $X$. By an $L$-twisted pair over $X$ we
mean a pair $(E,\varphi)$ consisting of a holomorphic
$H^\CC$-principal bundle, where $H^\CC$ is a complex reductive Lie
group and $\varphi$ is a holomorphic section of $E(B)\otimes L$,
where $E(B)$ is the vector bundle associated to a complex
representation $H^\CC\to \GL(B)$ and $L$ is a holomorphic line
bundle over $X$. We study in full the notion of polystability and
prove the correspondence between polystable pairs and solutions to
the corresponding Hermite--Einstein equations for a reduction of
the structure group of $E$ to $H$ --- the maximal compact subgroup
of $H^\CC$. This extends the correspondence for stable pairs of
\cite{banfield,bradlow-garcia-prada-mundet:2003} to the strictly
polystable case and solves the problem of completely
characterizing the pairs which support solutions to the equations.
The Hermite--Einstein equations combine the curvature term of the
classical Hermite--Einstein equation for polystable vector bundles
and a quadratic term on the Higgs field, which can be interpreted
as a moment map (see Theorem \ref{hk-twisted-pairs}). When the
general Hermite--Einstein equation is considered for $G$-Higgs
bundles, we call it the Hitchin equation.

In Section~\ref{reps} we study non-abelian Hodge theory over a compact
Riemann surface $X$ for a general connected semisimple Lie group $G$.
Let $G$ be a reductive real Lie group with maximal compact subgroup
$H\subset G$, let $K$ be the canonical line bundle over $X$ and let
$\glie=\hlie\oplus\mlie$ be the Cartan decomposition of $\glie$. Then
a $G$-Higgs bundle is a pair $(E,\varphi)$, consisting of a
holomorphic $\HC$-principal bundle $E$ over $X$ and a holomorphic
section $\varphi$ of $E(\mclie)\otimes K$. Here $E(\mclie)$ is the
$\mclie$-bundle associated to $E$ via the isotropy representation
$H^\CC\to \GL(\liem^\CC)$.  These objects are a particular case of the
general twisted pairs introduced in Section~\ref{appendix}.  We study
the deformations and the moduli spaces of $G$-Higgs bundles.  An
important result is the correspondence between the moduli space of
polystable $G$-Higgs bundles and the moduli space of solutions to the
Hitchin equations. This is well-known when $G$ is complex
\cite{hitchin:1987a, simpson:1988,simpson:1992} or compact
\cite{narasimhan-seshadri:1965,ramanathan:1975}, and a proof for any
$G$ follows from \cite{bradlow-garcia-prada-mundet:2003} in the case
of \emph{stable} $G$-Higgs bundles.  In this paper, we prove the
general case of polystable $G$-Higgs bundles.  The result (given by
Theorem \ref{hk}) is a consequence of the more general
Hitchin--Kobayashi correspondence given in Theorem
\ref{hk-twisted-pairs}.

We then consider the moduli space of reductive representations of the
fundamental group of a compact Riemann surface $X$ in a Lie group $G$.
By a representation we mean a homomorphism from $\pi_1(X)$ to $G$. A
representation is said to be reductive means that if its composition
with the adjoint representation of $G$ is completely reducible; when
$G$ is algebraic this is equivalent to the image of the representation
of $\pi_1(X)$ in $G$ having reductive Zariski closure.  Combining
Theorem~\ref{hk} with Corlette's existence theorem for harmonic metrics
\cite{corlette:1988}, we establish in Theorem \ref{na-Hodge} the
correspondence between this moduli space and the moduli space of
polystable $G$-Higgs bundles when $G$ is connected and semisimple.

In Section~\ref{twisted-higgs} we study how the stability condition
stated in general in Section~\ref{appendix} simplifies for $G$-Higgs
bundles for various groups.  This includes $G=\Sp(2n,\RR)$ ---
the group of linear transformations of $\RR^{2n}$ which preserve the
standard symplectic form --- and also other groups that naturally
contain $\Sp(2n,\RR)$, like $\Sp(2n,\CC)$, and $\SL(2n,\CC)$, as well
as $\GL(n,\RR)$.

The notion of an $L$-twisted $G$-Higgs pair is a slight generalization
of that of a $G$-Higgs bundle, where one allows a general line bundle
$L$ to play the role of the canonical bundle in the definition. Some
(though not all) of the results of Sections \ref{reps} and
\ref{twisted-higgs} apply in the setting of $L$-twisted $G$-Higgs
pairs at no extra cost and in these cases we choose to work in this
more general setting.

\subsubsection*{Acknowledgements}

The authors thank Olivier Biquard, Steven Bradlow, Carlos Florentino, 
Bill Goldman,
Nigel Hitchin, S.~Ramanan, Roberto Rubio and Alexander Schmitt for
useful conversations and shared insights. We also thank Olivier
Guichard for carefully reading an earlier version of this paper and
pointing out a number of typos and mistakes.

\section{Stability of twisted pairs and
Hitchin--Kobayashi correspondence}
 \label{appendix}

In this section we introduce a general notion of polystability for
pairs of the form $(E,\phi)$, where $E$ is a holomorphic principal
bundle and $\phi$ is a section of an associated vector bundle, and
we prove a Hitchin--Kobayashi correspondence for polystable pairs.
There have appeared in the literature several papers
(e.g.\ \cite{banfield, bradlow-garcia-prada-mundet:2003, mundet:2000,
simpson:1988}) with extensions of the original Hitchin--Kobayashi
correspondence due to Uhlenbeck and Yau \cite{uhlenbeck-yau:1988},
obtaining different levels of generality. Lest the reader think
that we have any pretension of founding a new literary genre on
slight variations of the Hitchin--Kobayashi correspondence, we now
briefly describe what are the new aspects which we consider,
compared to the previous existing papers.

The main novelty of the present paper regarding the
Hitchin--Kobayashi correspondence is the introduction and study of
a general notion of polystability which is equivalent, without any
additional hypothesis, to the existence of solution to the
Hermite--Einstein equations corresponding to the type of pair
considered. Polystability was of course well understood in the
case of vector bundles and some of their generalizations as
vortices, triples or Higgs bundles. However, the extensions of the
Hitchin--Kobayashi correspondence to general pairs which have
appeared so far deal only with stable objects (i.e., those for
which the degree inequalities are always strict) satisfying a
certain simplicity condition, and in this sense they are
unnecessarily restricted, as the intuition obtained from the case
of vector bundles suggests.

Roughly speaking, a pair $(E,\phi)$ is polystable if it is
semistable and the structure group of $E$ can be reduced to a
smaller subgroup so as to give rise to a stable pair (this
corresponds, in the vector bundle case, to the process of looking
at a polystable vector bundle as a direct sum of stable vector
bundles of the same slope). Our actual definition of polystability
(see Section \ref{ss:app-pol,sem,stability}) is not expressed
in this way, but rather in terms of reductions of the structure
group from parabolic subgroups to their Levi subgroups. The
existence of a reduction of the structure group leading to a
stable object is proved to be a consequence of polystability in
Section \ref{ss:Jordan-Holder}. We also prove the uniqueness of
such a reduction (which we call, following the usual terminology,
the Jordan--H\"older reduction).

Strictly polystable vector bundles can be distinguished from
stable vector bundles by the fact that their automorphism group
contains elements which are not homotheties. In Section
\ref{sec:inf-aut} we prove that something similar happens for
general pairs. The Hitchin--Kobayashi correspondence for
polystable pairs is proved in Section
\ref{ss:Hitchin-Kobayashi-correspondence}. Our strategy is to
reduce the proof to the case of stable pairs, for which we refer
to the result in \cite{bradlow-garcia-prada-mundet:2003}. Finally,
we prove in Section \ref{ss:simple-infinitesimally-simple} that
the automorphism group of a polystable pair is reductive. This is
a consequence of two facts: first, that the group of gauge
transformations which preserve a pair $(E,\phi)$ and the reduction
of $E$ solving the Hermite--Einstein equation is compact and,
second, that the full group of automorphisms of $(E,\phi)$ is the
complexification of the previous group (this is a general fact,
which follows formally from the moment map interpretation of the
equations).

We have included in this section some material on parabolic
subgroups which is perhaps classical but for which we did not find
any reference adapted to our point of view. These results are most
of the times only sketched, but we have tried to be careful in
setting the notation, so that all the notions which we are using
are clearly defined.

\subsection{Standard parabolic subgroups}

\label{ss:standard-parabolic-subgroups}

Let $H$ be a compact and connected Lie group and let $H^{\CC}$ be
its complexification. Parabolic subgroups of $H^{\CC}$ can be
defined in several different but equivalent ways. Here we list
some of them: (1) the subgroups $P\subset H^{\CC}$ such that the
homogeneous space $H^{\CC}/P$ is a projective variety, (2) any
subgroup containing a maximal closed and connected solvable
subgroup of $H^{\CC}$ (i.e., a Borel subgroup), (3) the
stabilizers of points at infinity of the visual compactification
of the symmetric space $H\backslash H^{\CC}$. Here we use a more
constructive definition: we first define standard parabolic
subgroups with respect to a root space decomposition, and then we
define a parabolic subgroup to be any subgroup which is conjugate
to a standard parabolic subgroup. The reader meeting this notion
for the first time is advised to think as an example on the
parabolic subgroups of $\GL(n,\CC)$, which are simply the
stabilizers of any partial flag $0\subset V_1\subset\dots
V_r\subset\CC^n$.

Here is some notation which will be used:


\begin{align*}
  H &- \text{a compact and connected Lie group }\\
  \HC &- \text{the complexification of $H$ } \\
  T\subset H &- \text{a maximal torus }\\
  \hlie &- \text{the Lie algebra of $H$ }\\
  \hclie &- \text{the Lie algebra of $\HC$ }\\
  \tlie\subset\hlie &- \text{the Lie algebra of $T$}\\
  \alie\subset\hclie &- \text{the complexification of $\tlie$,
    $\alie=\tlie\otimes_{\RR}\CC$}\\
  \hlie_s=[\hlie,\hlie] &- \text{the semisimple part of $\hlie$
  }\\
  \hclie_s=[\hclie,\hclie] &- \text{the semisimple part of $\hclie$
  }\\
  \zlie\subset\alie &- \text{the center of $\hclie$ }\\
  \clie\subset\hclie_s &- \text{the Cartan subalgebra of $\hclie_s$
    defined as $\clie=\alie\cap \hclie_s$ }\\
  \la\cdot,\cdot\ra &- \text{a non-degenerate invariant $\CC$-bilinear
    pairing on
    $\hclie$}\\
  R\subset\clie^*=\Hom_{\CC}(\clie,\CC) &-\text{the roots of
    $\hclie_{s}$}\\
  \hlie_{\delta}\subset\hclie &-\text{the root space corresponding
    to $\delta\in R$}\\
  \Delta\subset R&-\text{a choice of simple roots.}
\end{align*}


Using the previous notation we can write the root space
decomposition of $\hclie$ as:
$$\hclie=\zlie\oplus\clie\oplus\bigoplus_{\delta\in
R}\hlie_{\delta}.$$ For any $A\subset \Delta$ define $R_A$ to be
the set of roots of the form
$\delta=\sum_{\beta\in\Delta}m_{\beta}\beta\in R$ with
$m_{\beta}\geq 0$ for all $\beta\in A$ (so if $A=\emptyset$ then
$R_A=R$). Then
$$\plie_A=\zlie\oplus\clie\oplus\bigoplus_{\delta\in
R_A}\hlie_{\delta}$$ is a Lie subalgebra of $\hclie$. Denote by
$P_A\subset \HC$ the connected subgroup whose Lie algebra is
$\plie_A$.

\begin{definition}
A {\bf standard parabolic subgroup} of $\HC$ is any subgroup of
the form $P_A$, for any choice of subset $A\subset R$. A {\bf
parabolic subgroup} of $\HC$ is any subgroup which is conjugate to
a standard parabolic subgroup.
\end{definition}


Define similarly $R_A^0$ as the set of roots
$\delta=\sum_{\beta\in\Delta}m_{\beta}\beta$ with $m_{\beta}=0$
for all $\beta\in A$. The vector space
\begin{equation}
\label{eq:def-llie-A}
\llie_A=\zlie\oplus\clie\oplus\bigoplus_{\delta\in
R_A^0}\hlie_{\delta}
\end{equation}
is a Lie subalgebra of $\plie_A$. Let $L_A$ be the connected
subgroup with Lie algebra $\llie_A$. Then $L_A$ is a {\bf Levi
subgroup} of $P_A$, i.e., a maximal reductive subgroup of $P_A$.
Finally,
\begin{equation}
\label{eq:def-ulie-A} \ulie_A=\bigoplus_{\delta\in R_A\setminus
R_A^0}\hlie_{\delta}
\end{equation}
is also a Lie subalgebra of $\plie_A$, and the connected Lie group
$U_A\subset P_A$ with Lie algebra $\ulie_A$ is the unipotent
radical of $P_A$. $U_A$ is a normal subgroup of $P_A$ and the
quotient $P_A/U_A$ is naturally isomorphic to $L_A$ so we have
\begin{equation}
\label{eq:Levi-decomposition} P_A=L_AU_A.
\end{equation}


\subsection{Antidominant characters of $\plie_A$}

Recall that a character of a complex Lie algebra $\glie$ is a
complex linear map $\glie\to\CC$ which factors through the
quotient map $\glie\to\glie/[\glie,\glie]$. Here we classify the
characters of parabolic subalgebras $\plie_A\subset\hclie$. We
will see that all these characters come from elements of the dual
of the center of the Levi subgroup $\llie_A\subset\plie_A$. Then
we define antidominant characters.

Let $Z$ be the center of $\HC$, and let
$$\Gamma=\Ker(\exp:\zlie\to Z).$$ Then
$\zlie_{\RR}=\Gamma\otimes_{\ZZ}\RR\subset \zlie$ is the Lie
algebra of the maximal compact subgroup of $Z$.

Let
$\zlie_{\RR}^*=\Hom_{\RR}(\zlie_{\RR},\imag\RR)$ and let
$\Lambda=\{\lambda\in\zlie_{\RR}^*\mid\lambda(\Gamma)\subset
2\pi\imag\ZZ\}$.
Let $\{\lambda_{\delta}\}_{\delta\in\Delta}\subset\clie^*$ be the
set of fundamental weights of $\hclie_s$, i.e., the duals with
respect to the Killing form of the coroots
$\{2\delta/\la\delta,\delta\ra\}_{\delta\in\Delta}$. We extend any
$\lambda\in\Lambda$ to a morphism of complex Lie algebras
$$\lambda:\zlie\oplus\clie\to\CC$$ by setting
$\lambda|_{\clie}=0$, and similarly for any $\delta\in A$ we
extend $\lambda_\delta:\clie\to\CC$ to
$$\lambda_{\delta}:\zlie\oplus\clie_A\to\CC$$ by setting
$\lambda_{\delta}|_{\zlie}=0.$

\begin{lemma}
Define $\zlie_A=\bigcap_{\beta\in\Delta\setminus
A}\Ker\lambda_{\beta}$ if $A\neq\Delta$ and let $\zlie_A=\clie$ if
$A=\Delta$.
\begin{enumerate}
\item $\zlie_A$ is equal to the center of $\llie_A$, \item we have
$(\plie_A/[\plie_A,\plie_A])^*\simeq\zlie_A^*.$
\end{enumerate}
\end{lemma}
\begin{proof}
Both (1) and (2) follow from the fact that for any
$\delta,\delta'\in R$ we have
$[h_{\delta},h_{\delta'}]=h_{\delta+\delta'}$ if
$\delta+\delta'\neq 0$ and
$[h_{\delta},h_{-\delta}]=(\Ker\lambda_{\delta})^{\perp}$ (see
Theorem 2 in Chapter VI of \cite{serre:1987}).
\end{proof}

Let $\clie_A=\zlie_A\cap \llie_A$, so that
$\zlie_A=\zlie\oplus\clie_A.$ By the previous lemma, the
characters of $\plie_A$ are in bijection with the elements in
$\zlie^*\oplus\clie_A^*$.

\begin{definition}
  An {\bf antidominant character} of $\plie_A$ is any element of
  $\zlie^*\oplus\clie_A^*$ of the form $\chi=z+\sum_{\delta\in
    A}n_{\delta}\lambda_{\delta}$, where $z\in \zlie_{\RR}^*$ and each
  $n_{\delta}$ is a nonpositive real number. If for each $\delta\in A$
  we have $n_{\delta}<0$ then we say that $\chi$ is {\bf strictly
    antidominant}.
\end{definition}

The restriction of the invariant form
$\la\cdot,\cdot\ra$ to $\zlie\oplus\clie_A$ is non-degenerate, so
it induces an isomorphism
$\zlie^*\oplus\clie_A^*\simeq\zlie\oplus\clie_A$. For any
antidominant character $\chi$ we define
$s_{\chi}\in\zlie\oplus\clie_A\subset\zlie\oplus\clie$ to be the
element corresponding to $\chi$ via the previous isomorphism. One
checks that $s_{\chi}$ belongs to $\imag\hlie$.

\subsection{Exponentiating characters of $\plie_A$ to characters
of $P_A$} \label{ss:caracters-parabolics}

A character of a complex Lie group $G$ is a morphism of Lie groups
$G\to\CC^*$. Any character of $G$ induces a character of $\glie$.
When a character of $\glie$ comes from a character of $G$ then we
say that it exponentiates. In general there are (many) characters
of $\glie$ which do not exponentiate, but here we prove that the
set characters of $\plie_A$ which exponentiate generate (as a
subset of a vector space) the space of all characters of
$\plie_A$. This will be used to give an algebraic definition of
the degree of parabolic reductions in Section
\ref{ss:degree-reduction-character}.

Let $Z_A$ be the identity component of the center of $L_A$, and
let $L_A^{ss}$ be the connected subgroup of $L_A$ whose Lie
algebra is $[\llie_A,\llie_A]$. Then $L_A^{ss}$ is semisimple.
Define
$$Z^{ss}(L_A):=Z_A\cap L_A^{ss}.$$
The group $Z^{ss}(L_A)$ is a subgroup of the center of $L_A^{ss}$.
The center of a semisimple group over $\CC$ is finite, because it
coincides with the center of any of its maximal compact subgroups.
Hence $Z^{ss}(L_A)$ is finite. The product map $Z_A\times
L_A^{ss}\to L_A$ induces an isomorphism $L_A\simeq
Z_A\times_{Z^{ss}(L_A)}L_A^{ss}$, and projection to the first
factor gives a map $L_A\to Z_A/Z^{ss}(L_A)$. Composing this
projection with the quotient map $P_A\to P_A/U_A\simeq L_A$ we
obtain a morphism of Lie groups
$$\pi_A:P_A\to Z_A/Z^{ss}(L_A).$$
In the following lemma we use the fact that $Z^{ss}(L_A)$ is
finite.

\begin{lemma}
\label{lemma:n-exponentiates} There exists some positive integer
$n$ (depending on the fundamental group of $L_A$) such that for
any $\lambda\in\Lambda$ and any $\delta\in A$ the morphisms of Lie
algebras $n\lambda:\zlie\oplus\clie_A\to\CC$ and
$n\lambda_{\delta}:\zlie\oplus\clie_A\to\CC$ exponentiate to
morphisms of Lie groups
$$\exp(n\lambda):Z_A/Z^{ss}(L_A)
\to\CC^{\times},\qquad\qquad\qquad
\exp(n\lambda_{\delta}):Z_A/Z^{ss}(L_A) \to\CC^{\times}.$$
\end{lemma}

Composing the morphisms given by the previous lemma with the
morphism $\pi_A$ we get for any $\lambda\in\Lambda$ and $\delta\in
A$ morphisms of Lie groups
$$\kappa_{n\lambda}:P_A\to\CC^{\times},\qquad\qquad\qquad
\kappa_{n\delta}:P_A\to\CC^{\times}.$$

\subsection{Recovering a parabolic subgroup from its antidominant
characters}


\begin{lemma}
\label{lemma:stand-parab} Let $s\in \imag\hlie$ and define the
sets
\begin{align*}
\plie_s &:=\{x\in\hclie\mid \Ad(e^{ts})(x)\text{ is bounded as
$t\to\infty$} \}\subset\hclie, \\
\llie_s &:=\{x\in\hclie\mid [x,s]=0\ \}\subset\hclie, \\
P_s & :=\{g\in \HC\mid e^{t s}ge^{-t s} \text{ is bounded as
$t\to\infty$} \}\subset\HC, \\
L_s & :=\{g\in\HC\mid \Ad(g)(s)=s\ \}\subset\HC.
\end{align*}
The following properties hold:
\begin{enumerate}
\item\label{item:stand-parab-preliminar} Both $\plie_s$ and
$\llie_s$ are Lie subalgebras of $\hclie$ and $P_s$ and $L_s$ are
subgroups of $\HC$. Furthermore $P_s$ and $L_s$ are connected.
\item \label{item:stand-parab-1} Let $\chi$ be an antidominant
character of $P_A$. There are inclusions $\plie_A\subset
\plie_{s_{\chi}}$, $\llie_A\subset \llie_{s_{\chi}}$, $P_A\subset
P_{s_{\chi}}$ and $L_A\subset L_{s_{\chi}}$, with equality if
$\chi$ is strictly antidominant. \item \label{item:stand-parab-2}
For any $s\in\imag\hlie$ there exists $h\in H$ and a standard
parabolic subgroup $P_A$ such that $P_s=hP_Ah^{-1}$ and
$L_s=hL_Ah^{-1}$. Furthermore, there is an antidominant character
$\chi$ of $P_A$ such that $s=h s_{\chi} h^{-1}$.
\end{enumerate}
\end{lemma}
\begin{proof}
That $\llie_s$, $\plie_s$ are subalgebras and $L_s,P_s$ are
subgroups is immediate from the definitions. Let $T_s$ be the
closure of $\{e^{\imag ts}\mid t\in\RR\}$. Then $L_s$ is the
centralizer of the torus $T_s$ in $\HC$, so by Theorem 13.2 in
\cite{borel:1956}
$L_s$ is connected. To prove that $P_s$ is also connected, note
that if $g$ belongs to $P_s$, so that $e^{ts}ge^{-ts}$ is bounded
as $t\to\infty$, then the limit of $\pi_s(g):=e^{ts}ge^{-ts}$ as
$t\to\infty$ exists and belongs to $L_s$. Note by the way that the
resulting map $\pi_s:P_s\to L_s$ is a morphism of Lie groups which
can be identified with the projection $P_s\to P_s/U_s\simeq L_s$,
where
$$U_s=\{g\in \HC\mid e^{t s}ge^{-t s} \text{ converges to $1$ as
$t\to\infty$} \}\subset P_s$$ is the unipotent radical of $U_s$.
So if $g\in P_s$ then the map $\gamma:[0,\infty)\to \HC$ defined
as $\gamma(t)=e^{ts}ge^{-ts}$ extends to give a path from $g$ to
$L_s$, and since $L_s$ is connected it follows that $P_s$ is also
connected. This proves (\ref{item:stand-parab-preliminar}). Let
now $\chi=z+\sum_{\beta\in\Delta} n_{\beta}\lambda_{\beta}$ be an
antidominant character of $P_A$. Let $\delta=\sum_{\beta\in
\Delta}m_{\beta}\beta$ be a root and let $u\in\hlie_{\delta}$. We
have $[s_{\chi},u]=\la s_{\chi},\delta\ra u=\la\chi,\delta\ra
u=(\sum_{\beta\in\Delta}m_{\beta}n_{\beta}\la\beta,\beta\ra/2)u$.
Hence $\Ad(e^{ts_{\chi}})(u)=(\sum_{\beta\in\Delta}\exp(t
n_{\beta} m_{\beta}\la\beta,\beta\ra/2))u$, so this remains
bounded as $t\to\infty$ if $m_{\beta}\geq 0$ for any $\beta$ such
that $n_{\beta}\leq 0$. This implies that $\plie_A\subset\plie_s$
and $\llie_A\subset\llie_s$ and that the inclusions are equalities
when $\chi$ is strictly dominant. The analogous statements for
$P_A,L_A,P_s,L_s$ follow from this, because the subgroups
$P_A,L_A,P_s,L_s$  are connected. Hence (\ref{item:stand-parab-1})
is proved. To prove (\ref{item:stand-parab-2}) take a maximal
torus $T_s$ containing $\{e^{\imag ts}\mid t\in\RR\}$ and choose
$h\in H$ such that $h^{-1}T_sh=T$ and $\Ad(h^{-1})(s)$ belongs to
the Weyl chamber in $\tlie$ corresponding to the choice of
$\Delta\subset R$. Then use (\ref{item:stand-parab-1}).
\end{proof}

\begin{lemma}
\label{lemma:tots-parabolics-P-s} Let $P\subset\HC$ be any
parabolic subgroup, conjugate to $P_A$. Let $\chi$ be an
antidominant character of $\plie_A$. There exists an element
$s_{P,\chi}\in\imag\hlie$, depending smoothly on $P$, which is
conjugate to $s_{\chi}$ and such that $P\subset P_{s_{P,\chi}}$,
with equality if and only if $\chi$ is strictly antidominant.
\end{lemma}
\begin{proof}
Assume that $P=gP_Ag^{-1}$ for some $g\in\HC$. From the well known
equality $\HC/P_A=H/(P_A\cap H)=H/(L_A\cap H)$ we deduce that
there exists some $h\in H$ such that $P=hP_Ah^{-1}$. Then we set
$s_{P,\chi}=hs_{\chi}h^{-1}$. This is well defined because $h$ is
unique up to multiplication on the right by elements of $L_A\cap
H$, and these elements commute with $s_{\chi}$.
\end{proof}

\subsection{Principal bundles and parabolic subgroups}

\label{ss:principal-bundles-parabolic}

If $E$ is a $\HC$-principal holomorphic bundle over $X$ and $M$ is
any set on which $\HC$ acts on the left, we denote by $E(M)$ the
twisted product $E\times_{\HC}M$, defined as the quotient of
$E\times M$ by the equivalence relation $(eh,m)\sim(e,hm)$ for any
$e\in E$, $h\in\HC$ and $m\in M$. The sections $\varphi$ of $E(M)$
are in natural bijection with the maps $\phi:E\to M$ satisfying
$\varphi(eh)=h^{-1}\varphi(e)$ for any $e\in E$ and $h\in\HC$ (we
call such maps antiequivariant). Furthermore, $\phi$ is
holomorphic if and only if $\varphi$ is holomorphic.

If $M$ is a vector space (resp.\ complex variety) and the action of
$\HC$ on $M$ is linear (resp. holomorphic) then $E(M)$ is a vector
bundle (resp.\ holomorphic fibration). In this situation, for any
complex line bundle $L\to X$ we can form a vector bundle
$E(M)\otimes L$ which can be identified with $E^L(M)$, where $E^L$
denotes the principal $\HC\times\CC^{\times}$ bundle
$E^L=\{(e,l)\in E\times_XL\mid l\neq 0\}$ and we form the
associated product by making $(h,\lambda)\in\HC\times\CC^{\times}$
act on $m\in M$ as $\lambda hm$. Consequently, the sections of
$E(M)\otimes L$ can be identified with antiequivariant maps
$E^L\to M$.

Let $B$ be a Hermitian vector space and let $\rho:H\to \U(B)$ be a
unitary representation. The morphism $\rho$ extends to a
holomorphic representation of $\HC$ in $\GL(B)$, which we denote
also by $\rho$. Suppose that $P_A\subset \HC$ is the parabolic
subgroup corresponding to a subset $A\subset\Delta$ and let $\chi$
be an antidominant character. Define
$$B_{\chi}^-=\{v\in B \mid \rho(e^{t s_{\chi}})v\text{ remains bounded
as $\RR\ni t\to\infty$} \}.$$ This is a complex vector subspace of
$B$ and by (\ref{item:stand-parab-1}) in Lemma
\ref{lemma:stand-parab} it is invariant under the action of $P_A$.
Define also
$$B_{\chi}^0=\{v\in B \mid \rho(e^{t s_{\chi}})v=v \text{ for any $t$ }
\}\subset B_{\chi}^-.$$ This is a complex subspace of $B^-_{\chi}$
and, using again (\ref{item:stand-parab-1}) in Lemma
\ref{lemma:stand-parab}, we deduce that $B^0_{\chi}$ is invariant
under the action of $L_A$.

Suppose that $\sigma$ is a holomorphic section of $E(\HC/P_A)$.
Since $E(\HC/P_A)\simeq E/P_A$ canonically and the quotient $E\to
E/P_A$ has the structure of a $P_A$-principal bundle, the pullback
$E_{\sigma}:=\sigma^*E$ is a $P_A$-principal bundle over $X$, and
we can identify canonically $E\simeq E_{\sigma}\times_{P_A}\HC$ as
principal $\HC$-bundles (hence, $\sigma$ gives a reduction of the
structure group of $E$ to $P_A$). Equivalently, we can look at
$E_{\sigma}$ as a holomorphic subvariety $E_{\sigma}\subset E$
invariant under the action of $P_A\subset\HC$ and inheriting a
structure of principal bundle. It follows that $E(B)\simeq
E_{\sigma}\times_{P_A}B$, so the vector bundle
$E_{\sigma}\times_{P_A}B_{\chi}^-$ can be identified with a
holomorphic subbundle $$E(B)_{\sigma,\chi}^-\subset E(B).$$

Now suppose that $\sigma_L$ is a holomorphic section of
$E_{\sigma}(P_A/L_A)$. This section induces, exactly as before, a
reduction of the structure group of $E_{\sigma}$ from $P_A$ to
$L_A$. So we obtain from $\sigma_L$ a principal $L_A$ bundle
$E_{\sigma_L}$ and an isomorphism $E_{\sigma}\simeq
E_{\sigma_L}\times_{L_A}P_A$. Hence $E(B)\simeq
E_{\sigma_L}\times_{L_A}B$, and we can thus identify the vector
bundle $E_{\sigma_L}\times_{L_A}B_{\chi}^0$ with a holomorphic
subbundle
$$E(B)_{\sigma_L,\chi}^0\subset E(B)_{\sigma,\chi}^-.$$

\subsection{Degree of a reduction and an
antidominant character} \label{ss:degree-reduction-character} Let
$\sigma$ denote a reduction of the structure group of $E$ to a
standard parabolic subgroup $P_A$ and let $\chi$ be an
antidominant character of $\plie_A$. Let us write
$\chi=z+\sum_{\delta\in A}n_{\delta}\lambda_{\delta}$, with
$z\in\zlie_{\RR}^*$, and $z=z_1\lambda_1+\dots+z_r\lambda_r$,
where $\lambda_1,\dots,\lambda_r\in\Lambda$ and the $z_j$ are real
numbers. Let $n$ be an integer as given by Lemma
\ref{lemma:n-exponentiates}. Using the characters
$\kappa_{n\lambda},\kappa_{n\delta}:P_A\to\CC^{\times}$ defined in
Section \ref{ss:caracters-parabolics} we can construct from the
principal $P_A$ bundle $E_{\sigma}$ line bundles
$E_{\sigma}\times_{\kappa_{n\lambda}}\CC$ and
$E_{\sigma}\times_{\kappa_{n\delta}}\CC$.

\begin{definition}
  \label{def:degree-E-sigma-chi}
  We define the {\bf degree} of the bundle $E$ with respect to the
  reduction $\sigma$ and the antidominant character $\chi$ to be the
  real number:
  \begin{equation}
    \label{eq:def-deg-sigma}
    \deg(E)(\sigma,\chi):=\frac{1}{n}\left(\sum_j z_j
      \deg(E_{\sigma}\times_{\kappa_{n\lambda_j}}\CC)+\sum_{\delta\in A}
      n_{\delta}\deg(E_{\sigma}\times_{\kappa_{n\delta}}\CC)\right).
  \end{equation}
  This expression is independent of the choice of the $\lambda_j$'s
  and the integer $n$.
\end{definition}
There is the following, more intrinsic, way to define the
degree with respect to a reduction to \emph{any} parabolic subgroup
$P$ and an antidominant character $\chi$  of $\liep$: let $n$ be such
that $\chi^n$ lifts to a character $\tilde\chi$ of $P$. Then
\begin{displaymath}
  \deg(E)(\sigma,\chi) = \frac{1}{n}\deg E_{\sigma}(\tilde\chi).
\end{displaymath}

We now give another definition of the degree in terms of the
curvature of connections, in the spirit of Chern--Weil theory.
This definition is shorter and more natural from the point of view
of proving the Hitchin--Kobayashi correspondence (but, as we said,
in this paper we do not give a complete proof of it: we just
reduce our general result to the one obtained in
\cite{bradlow-garcia-prada-mundet:2003} for simple stable pairs;
this is why the reader will not find any use of the following
formula in the present paper). On the other hand, the definition
in terms of Chern--Weil theory uses obviously transcendental
methods, so it is not satisfying from the point of view of
obtaining a polystability condition of purely algebraic nature.

Define $H_A=H\cap L_A$ and $\hlie_A=\hlie\cap\llie_A$. Then $H_A$
is a maximal compact subgroup of $L_A$, so the inclusions
$H_A\subset L_A$ is a homotopy equivalence. Since the inclusion
$L_A\subset P_A$ is also a homotopy equivalence, given a reduction
$\sigma$ of the structure group of $E$ from $\HC$ to $P_A$ one can
further restrict the structure group of $E$ to $H_A$ in a unique
way up to homotopy. Denote by $E'_{\sigma}$ the resulting $H_A$
principal bundle. Let $\pi_A:\plie_A\to\zlie\oplus\clie_A$ be the
differential of the projection $\pi_A$ defined in Section
\ref{ss:caracters-parabolics}. Let $\chi=z+\sum_{\delta\in
A}n_{\delta}\lambda_{\delta}$ be an antidominant character. Define
$\kappa_{\chi}=(z+\sum_{\delta}n_{\delta} \lambda_{\delta})\circ
\pi_A\in \plie_A^*$. Let $\hlie_A\subset\llie_A\subset\plie_A$ be
the Lie algebra of $H_A$. Then
$\kappa_{\chi}(\hlie_A)\subset\imag\RR$. Choose a connection $\AA$
on $E'_{\sigma}$ and denote by
$F_{\AA}\in\Omega^2(X,E'_{\sigma}\times_{\Ad}\hlie_A)$ its
curvature. Then $\kappa_{\chi}(F_{\AA})$ is a $2$-form on $X$ with
values in $\imag\RR$, and we have
$$\deg(E)(\sigma,\chi):=\frac{\imag}{2\pi}\int_X
\kappa_{\chi}(F_{\AA}).$$

\subsection{$L$-twisted pairs and stability}

\label{ss:app-pol,sem,stability}

Let $X$ be a closed Riemann surface and let $L$ be a holomorphic line
bundle over $X$. Let $\HC$ be a connected complex reductive Lie group
and let $\rho\colon \HC \to \GL(B)$ be a representation. We denote its
derivative by the same letter
\begin{math}
  \rho\colon \liehc \to \lie{gl}(B)
\end{math}
and define
\begin{displaymath}
  \lieh_{\rho} = \lieh_s + \ker(\rho_{|\liez_{\RR}})^{\perp},
\end{displaymath}
where the the orthogonal complement is taken with respect to the
restriction of $\la \cdot , \cdot \ra$ to $\liez_{\RR}$. Thus
\begin{displaymath}
  \lieh = \lieh_{\rho} + \ker(\rho_{|\liez_{\RR}}).
\end{displaymath}

\begin{definition}
  \label{def:L-twisted-pair}
  An \textbf{$L$-twisted pair} is a pair of the form $(E,\varphi)$, where $E$
  is a holomorphic $\HC$-principal bundle over $X$ and $\varphi$ is a
  holomorphic section of $E(B)\otimes L$. When it does not lead to
  confusion we say that $(E,\varphi)$ is a pair, instead of an
  $L$-twisted pair.
\end{definition}

\begin{definition}
  \label{def:L-twisted-pairs-stability}
  Let $(E,\varphi)$ be an $L$-twisted pair and let
  $\alpha\in\imag\zlie_{\RR}\subset\zlie$. We say that $(E,\varphi)$
  is:

  \begin{itemize}
  \item $\alpha$-{\bf semistable} if: for any parabolic subgroup
    $P\subset \HC$, any antidominant character $\chi$ of $\liep$, and
    any holomorphic section $\sigma\in\Gamma(E(\HC/P))$ such that
    $\varphi\in \H^0(E(B)_{\sigma,\chi}^-\otimes L)$, we have
$$\deg(E)(\sigma,\chi)-\la\alpha,\chi\ra\geq 0.$$

\item $\alpha$-{\bf stable} if it is $\alpha$-semistable and
  furthermore: for any parabolic subgroup $P\subset \HC$ which is of
  the form $P = P_s $ for some $s\in \imag \lieh$ such that $s\not\in
  \ker(\rho_{|\liez_{\RR}})$, any antidominant character $\chi$
  of $\liep$ such that $\chi\not\in\ker(\rho_{|\liez_{\RR}})^*$, and
  any holomorphic section $\sigma\in\Gamma(E(\HC/P))$ such that
  $\varphi\in \H^0(E(B)_{\sigma,\chi}^-\otimes L)$, we have
$$\deg(E)(\sigma,\chi)-\la\alpha,\chi\ra> 0.$$

\item $\alpha$-{\bf polystable} if it is $\alpha$-semistable and for
  any $P$, $\chi$ and $\sigma$ as in the definition of
  $\alpha$-stability, such that $\varphi\in
  \H^0(E(B)_{\sigma,\chi}^-\otimes L)$, and $\chi$ is
  strictly antidominant, and such that
$$\deg(E)(\sigma,\chi)-\la\alpha,\chi\ra=0,$$
there is a holomorphic reduction of the structure group
$\sigma_{L_s}\in\Gamma(E_{\sigma}(P/L_s))$, where $E_{\sigma}$ denotes
the principal $P$-bundle obtained from the reduction $\sigma$ of the
structure group and $L_s \subset P$ is the Levi. Furthermore, under
these hypothesis $\varphi$ is required to belong to
$\H^0(E(B)_{\sigma_{L_s},\chi}^0\otimes L)\subset
\H^0(E(B)_{\sigma,\chi}^-\otimes L)$.
\end{itemize}
\end{definition}


\begin{remark}
  If we had stated the previous conditions considering reductions to
  only standard parabolic subgroups of $\HC$ then we would have obtained
  the same definitions. Indeed, since any parabolic subgroup is (for
  us, by definition) conjugate to a standard parabolic subgroup, the
  reductions of the structure group of $E$ to arbitrary parabolic
  subgroups are essentially the same as the reductions to standard
  parabolic subgroups.
\end{remark}

\begin{remark}
  The readers who are familiar with the stability condition for
  principal bundles as studied by Ramanathan \cite{ramanathan:1975}
  might find it surprising that our stability condition refers to
  antidominant characters of the parabolic Lie subalgebra and not only
  to characters of the parabolic subgroups (there are much less of the
  latter than of the former). The reason is that in the course of
  proving the Hitchin--Kobayashi correspondence one is naturally led
  to consider arbitrary antidominant characters of Lie subalgebras. It
  might be the case that the previous conditions do not vary if we
  only consider characters of the parabolic subgroups, but this is not
  at all obvious. We hope to come back to this question in the future.

\end{remark}

\subsection{The stability condition in terms of filtrations}
\label{ss:app-filtrations}

In order to obtain a workable notion of
$\alpha$-(poly,semi)stability it is desirable to have a more
concrete way to describe, for any holomorphic $\HC$-principal
bundle $E$,
\begin{itemize}
\item the reductions of the structure group of $E$ to parabolic
subgroups $P\subset\HC$, and the (strictly or not) antidominant
characters of $P$, \item the subbundle
$E(B)_{\sigma,\chi}^-\subset E(B)$, \item the degree
$\deg(E)(\sigma,\chi)$ defined in (\ref{eq:def-deg-sigma}), \item
reductions to Levi factors of parabolic subgroups and the
corresponding vector bundle $E(B)_{\sigma_L,\chi}^0\subset
E(B)_{\sigma,\chi}^-$.
\end{itemize}

We now discuss how to obtain in some cases such concrete
descriptions, beginning with the notion of degree. In
\cite{bradlow-garcia-prada-mundet:2003} the degree
$\deg(E)(\sigma,\chi)$ is defined in terms of a so-called
auxiliary representation (see \S 2.1.2 in
\cite{bradlow-garcia-prada-mundet:2003}) and certain linear
combinations of degrees of subbundles. The following lemma implies
that definition (\ref{eq:def-deg-sigma}) contains the one given in
\cite{bradlow-garcia-prada-mundet:2003} as a particular case.
Suppose that $\rho_W:H\to \U(W)$ is a representation on a
Hermitian vector space, and denote the holomorphic extension
$\HC\to \GL(W)$ with the same symbol $\rho_W$. Let
$(\Ker\rho_W)^{\perp}\subset\hclie$ be the orthogonal with respect
to invariant pairing on $\hclie$ of the kernel of
$\rho_W:\hclie\to\gllie(W)$, and let
$\pi:\hclie\to(\Ker\rho_W)^{\perp}$ be the orthogonal projection.

\begin{lemma}
\label{lemma:parabolic-representation} Take some element
$s\in\imag\hlie$. Then $\rho_W(s)$ diagonalizes with real
eigenvalues $\lambda_1<\dots<\lambda_k$. Let $W_j=\Ker
(\lambda_j\Id_W-\rho_W(s))$ and define $W_{\leq
i}=\bigoplus_{j\leq i}W_j$.

\begin{enumerate}
\item\label{item:parabolic-representation-1} The subgroup
$P_{W,s}\subset\HC$ consisting of those $g$ such that
$\rho_W(g)(W_{\leq i})\subset W_{\leq i}$ for any $i$ is a
parabolic subgroup, which can be identified with $P_{\pi(s)}$. Let
$\chi\in(\zlie\oplus\clie)^*$ be a character such that
$s_{\chi}=s$. Then $\chi$ is strictly antidominant for $P_{W,s}$.

\item\label{item:parabolic-representation-2} Suppose that for any
$a,b\in(\Ker\rho_W)^{\perp}$ we have $\la
a,b\ra=\Tr\rho_W(a)\rho_W(b)$. Let $u\in(\Ker\rho_W)^{\perp}$ be
any element, and write $\rho_W(u)=\sum\rho_W(u)_{ij}$ the
decomposition in pieces $\rho_W(u)_{ij}\in\Hom(W_i,W_j)$. Then
\begin{equation}
\label{eq:pairing-dominant} \la
\chi,u\ra=\Tr(\rho_W(s)\rho_W(u))=\lambda_k\Tr\rho_W(u)+
\sum_{i=1}^{k-1}(\lambda_i-\lambda_{i+1})\Tr\rho_W(u)_{ii}.
\end{equation}

\item \label{item:comparing-degree} Suppose that $\rho_W$
satisfies the conditions of
(\ref{item:parabolic-representation-2}). Let $E$ be a holomorphic
$\HC$-principal bundle and let $\WWW=E(W)$ be the associated
holomorphic vector bundle. Let $\sigma$ be a reduction of the
structure group of $E$ to a parabolic subgroup $P$ and an let
$\chi$ be an antidominant character of $P$. The endomorphism
$\rho_W(s_{\chi})$ diagonalizes with constant eigenvalues, giving
rise to a decomposition $\WWW=\bigoplus_{j=1}^k\WWW_j$, where
$\rho_W(s_{\chi})$ restricted to $\WWW_j$ is multiplication by
$\lambda_j\in\RR$. Suppose that $\lambda_1<\dots<\lambda_k$. For
each $i$ the subbundle $\WWW_{\leq i}=\bigoplus_{j\leq
i}\WWW_j\subset\WWW$ is holomorphic. We have:
$$\deg(E)(\sigma,\chi)=\lambda_k\deg
\WWW+\sum_{i=1}^{k-1}(\lambda_i-\lambda_{i+1})\deg\WWW_{\leq i}.$$
\end{enumerate}
\end{lemma}

\begin{proof}
The first assertion and formula (\ref{eq:pairing-dominant})
follows from easy computations. (\ref{item:comparing-degree})
follows from (\ref{item:parabolic-representation-2}).
\end{proof}

\begin{remark}
\label{rem:parabolic-representation} Condition
(\ref{item:parabolic-representation-2}) of the lemma is satisfied
when $W=\hlie$, endowed with the invariant metric, and
$\rho_W:\hclie\to\End W$ is the adjoint representation, since the
invariant metric on $\hlie$ is supposed to extend the Killing
pairing in the semisimple part $\hlie_s$.
\end{remark}

To clarify the other ingredients in the definition of
(poly,semi)stability, we put ourselves in the situation where
$\HC$ is a classical group. Let $\rho:\HC\to\GL(N,\CC)$ be the
fundamental representation. Suppose that $E$ is an $\HC$-principal
bundle, and denote by $V$ the vector bundle associated to $E$ and
$\rho$. One can describe pairs $(\sigma,\chi)$ consisting of a
reduction $\sigma$ of the structure group of $E$ to a parabolic
subgroup $P\subset \HC$ and an antidominant character $\chi$ of
$P$ in terms of filtrations of vector bundles
\begin{equation}
\label{eq:filtration} \VVV=(0\subsetneq
V_1\subsetneq\dots\subsetneq V_{k-1}\subsetneq V_k=V),
\end{equation}
and increasing sequences of real numbers (usually called weights)
\begin{equation}
\label{eq:weights} \lambda_1<\dots<\lambda_k,
\end{equation}
which are arbitrary if $\HC=\GL(n,\CC)$, and which satisfy
otherwise:
\begin{itemize}
\item if $\HC=\O(n,\CC)$ then\footnote{Even though $\O(n,\CC)$ is a
    disconnected group, the stability condition makes perfect sense for
  this group.}, for any $i$,
$V_{k-i}=V_i^{\perp}=\{v\in V\mid \la v,V_i\ra=0\}$, where
$\la,\ra$ denotes the bilinear pairing given by the orthogonal
structure (we implicitly define $V_0=0$), and
$\lambda_{k-i+1}+\lambda_i=0$.

\item if $\HC=\Sp(2n,\CC)$ then,
for any $i$, $V_{k-i}=V_i^{\perp}=\{v\in V\mid \omega(v,V_i)=0\}$,
where $\omega$ is the symplectic form on $V$ (as before, $V_0=0$),
and furthermore $\lambda_{k-i+1}+\lambda_i=0$.
\end{itemize}
The resulting character $\chi$ is strictly antidominant if all the
inequalities in (\ref{eq:weights}) are strict.

Given positive integers $p,q$ define the vector bundle
$V^{p,q}=V^{\otimes p}\otimes (V^*)^{\otimes q}$. For any choice
of reduction and antidominant character $(\sigma,\chi)$ specified
by a filtration (\ref{eq:filtration}) and weights
(\ref{eq:weights}) we define $$(V^{p,q})_{\sigma,\chi}^-=
\sum_{\lambda_{i_1}+\dots+\lambda_{i_p}\leq
\lambda_{j_1}+\dots+\lambda_{j_q}}V_{i_1}\otimes\dots\otimes
V_{i_p}\otimes V_{j_1}^{\perp}\otimes\dots\otimes
V_{j_q}^{\perp}\subset V^{p,q},$$ where $V_j^{\perp}=\{v\in
V^*\mid \la v,V_j\ra=0\}$ and $\la,\ra$ is the natural pairing
between $V$ and $V^*$. Since $\HC$ is a classical group, there is
an inclusion of representations
$$B\subset (\rho^{\otimes p_1}\otimes (\rho^*)^{\otimes q_1})
\oplus\dots\oplus (\rho^{\otimes p_r}\otimes (\rho^*)^{\otimes
q_r}),$$ so that the vector bundle $E(B)$ is contained in
$V^{p_1,q_1}\oplus\dots\oplus V^{p_r,q_r}$. One then has
$$E(B)_{\sigma,\chi}^-=E(B)\cap
((V^{p_1,q_1})_{\sigma,\chi}^-\oplus\dots\oplus
(V^{p_r,q_r})_{\sigma,\chi}^-).$$ Suppose that the invariant
pairing $\la\;,\,\ra$ on the Lie algebra $\hclie$ is defined using the
fundamental representation as $\la x,y\ra=\Tr \rho(x)\rho(y)$.
This clearly satisfies the condition of
(\ref{item:parabolic-representation-2}) of Lemma
\ref{lemma:parabolic-representation}, so by
(\ref{item:comparing-degree}) in the same lemma we have

$$\deg(E)(\sigma,\chi)=\lambda_k\deg
V+\sum_{i=1}^{k-1}(\lambda_i-\lambda_{i+1})\deg V_i.$$

We now specify what it means to have a reduction to a Levi factor
of a parabolic subgroup, as appears in the definition of
polystability. Assume that $(\sigma,\chi)$ is a pair specified by
(\ref{eq:filtration}) and (\ref{eq:weights}), so that $\sigma$
defines a reduction of the structure group of $E$ to a parabolic
subgroup $P\subset\HC$, and that $\varphi\in \H^0(L\otimes
E(B)_{\sigma,\chi}^-)$ and $\deg(E)(\sigma,\chi)=0$. If the pair
$(E,\varphi)$ is $\alpha$-polystable all these assumptions imply
the existence of a further reduction $\sigma_L$ of the structure
group of $\HC$ from $P$ to a Levi factor $L\subset P$; this is
given explicitly by an isomorphism of vector bundles
$$V\simeq \Gr\VVV:=V_1\oplus V_2/V_1\oplus\dots\oplus V_k/V_{k-1}$$
with respect to which $V_j\subset V$ corresponds to
$V_1\oplus V_2/V_1\oplus\dots\oplus V_j/V_{j-1}$ for each $j$.
When $\HC=\GL(n,\CC)$ such isomorphism is arbitrary. When $\HC$ is
$\O(n,\CC)$ (resp. $\Sp(2n,\CC)$), it is also assumed that the
pairing  of an element of $V_j/V_{j-1}$ with an element of
$V_i/V_{i-1}$, using the scalar product (resp. symplectic form),
is always zero unless $j+i=k+1$. We finally describe the bundle
$E(B)_{\sigma_L,\chi}^0$ in this situation. Let
$$(\Gr \VVV^{p,q})_{\sigma_L,\chi}^0=
\sum_{\lambda_{i_1}+\dots+\lambda_{i_p}=
\lambda_{j_1}+\dots+\lambda_{j_q}}(V_{i_1}/V_{i_1-1})
\otimes\dots\otimes (V_{i_p}/V_{i_p-1})\otimes
(V_{j_1}^{\perp}/V_{j_1+1}^{\perp})\otimes\dots\otimes
(V_{j_q}^{\perp}/V_{j_q+1}^{\perp}).$$ Then
$$E(B)_{\sigma_L,\chi}^0=E(B)\cap((\Gr \VVV^{p_1,q_1})_{\sigma_L,\chi}^0
\oplus\dots\oplus (\Gr \VVV^{p_r,q_r})_{\sigma_L,\chi}^0).$$

\subsection{Infinitesimal automorphism space}
\label{sec:inf-aut}

For any pair $(E,\varphi)$ we define the infinitesimal
automorphism space of $(E,\varphi)$ as
$$\aut(E,\varphi)=\{s\in \H^0(E(\hclie))\mid \rho(s)(\varphi)=0\},$$
where we denote by $\rho:\hclie\to\End(B)$ the morphism of Lie
algebras induced by $\rho$. We similarly define the semisimple
infinitesimal automorphism space of $(E,\varphi)$ as
$$\aut^{ss}(E,\varphi)=\{s\in\aut(E,\varphi)\mid s(x)\text{ is
semisimple for any $x\in X$ }\}.$$

\begin{proposition}
\label{proposition:stable-simple} Suppose that $(E,\varphi)$ is a
$\alpha$-polystable pair. Then $(E,\varphi)$ is $\alpha$-stable if
and only if $\aut^{ss}(E,\varphi)\subset \H^0(E(\zlie))$.
Furthermore, if $(E,\varphi)$ is $\alpha$-stable then we also have
$\aut(E,\varphi)\subset\H^0(E(\zlie))$.
\end{proposition}
\begin{proof}
Suppose that $(E,\varphi)$ is $\alpha$-polystable and that
$\aut^{ss}(E,\varphi)=\H^0(E(\zlie))$. We prove that $(E,\varphi)$
is $\alpha$-stable by contradiction. If $(E,\varphi)$ were not
$\alpha$-stable, then there would exist a parabolic subgroup
$P_A\subsetneq\HC$, a holomorphic reduction
$\sigma\in\Gamma(E/P_A)$, a strictly antidominant character $\chi$
such that $\deg(E)(\sigma,\chi)-\la\alpha,\chi\ra=0$, and a
further holomorphic reduction $\sigma_L\in\Gamma(E_{\sigma}/L_A)$
to the Levi $L_A$ (here $E_{\sigma}$ is the principal $P_A$ bundle
given by $\sigma$, satisfying $E_{\sigma}\times_{P_A}\HC\simeq E$)
such that $\varphi\in\H^0(E(B)_{\sigma_L,\chi}^0\otimes L)$. Since
the adjoint action of $L_A$ on $\hclie$ fixes $s_{\chi}$, there is
an element $$s_{\sigma,\chi}\in\H^0(E_{\sigma_L}(\hclie))\simeq
\H^0(E(\hclie))$$ which coincides fiberwise with $s_{\chi}$. On
the other hand $s_{\chi}$ is semisimple because it belongs to
$\imag\hlie$. The condition that
$\varphi\in\H^0(E(B)_{\sigma_L,\chi}^0\otimes L)$ implies that
$\rho(s_{\sigma,\chi})(\varphi)=0$, so
$s_{\sigma,\chi}\in\aut^{ss}(E,\varphi)$. And the condition that
$P_A\neq\HC$ implies that $s_{\chi}\notin\zlie$. This contradicts
the assumption that $\aut^{ss}(E,\varphi)=\H^0(E(\zlie))$, so
$(E,\varphi)$ is $\alpha$-stable.

Now suppose that $(E,\varphi)$ is $\alpha$-stable. We want to
prove that $\aut(E,\varphi)=\H^0(E(\zlie))$. Let
$\xi\in\aut(E,\varphi)$. Since $\xi$ is a section of
$E\times_{\HC}\hclie$, it can be viewed as an antiequivariant
holomorphic map $\psi:E\to\hclie$. The bundle $E$ is algebraic (to
prove this, take a faithful representation $\HC\to\GL(n,\CC)$ and
use the fact that any holomorphic vector bundle over an algebraic
curve is algebraic), so by Chow's theorem $\psi$ is algebraic.
Hence $\psi$ induces an algebraic map $\varphi:X\to\hclie\qu \HC$,
where $\hclie\qu \HC$ denotes the affine quotient, which is an
affine variety. Since $X$ is proper, $\varphi$ is constant, hence
it is contained in a unique fiber $Y:=\pi^{-1}(y)\subset\hclie$,
where $\pi:\hclie\to\hclie\qu \HC$ is the quotient map.

By a standard results on affine quotients, there is a unique
closed $\HC$ orbit $\OOO\subset Y$, and by a theorem of Richardson
the elements in $\OOO$ are all semisimple.

Consider the map
$\sigma:Y\to\OOO$ which sends any $y\in Y$ to $y_s$, where
$y=y_s+y_n$ is the Jordan decomposition of $y$ (see for example
\cite{borel:1991}). We claim that this map is algebraic (note that
the Jordan decomposition, when defined on the whole Lie algebra
$\hclie$, is not even continuous). To prove the claim first
consider the case $\hclie=\gllie(n,\CC)$. Then
$Y\subset\gllie(n,\CC)$ is the set of $n\times n$ matrices with
characteristic polynomial equal to some fixed polynomial, say
$\prod(x-\lambda_i)^{m_i}$, with $\lambda_i\neq\lambda_j$ for
$i\neq j$. By the Chinese remainder theorem there exists a
polynomial $P\in\CC[t]$ such that $P\equiv\lambda_i \mod
(t-\lambda_i)^{m_i}$ and $P\equiv 0 \mod t$. Then the map
$\sigma:Y\to\OOO$ is given by $\sigma(A)=P(A)$, which is clearly
algebraic. The case of a general $\hclie$ can be reduced to the
previous one using the adjoint representation
$\ad:\hclie\to\End(\hclie)\simeq \gllie(\dim\hclie,\CC)$.

By construction $\sigma$ is equivariant, so it induces a
projection $p_E:\H^0(E(Y))\to \H^0(E(\OOO))$. We define
$\xi_s=p_E(\xi)$ and $\xi_n=\xi-\xi_s$. Note that the
decomposition $\xi=\xi_s+\xi_n$ is simply the fiberwise Jordan
decomposition of an element of the Lie algebra as the sum of a
semisimple element plus a nilpotent one. We claim that both
$\xi_s$ and $\xi_n$ belong to $\aut(E,\phi)$. To prove this we
have to check that $\rho(\xi_s)(\phi)=\rho(\xi_n)(\phi)=0$. But
$\rho(\xi)=\rho(\xi_s)+\rho(\xi_n)$ is fiberwise the Cartan
decomposition of $\rho(\xi)$, since Cartan decomposition commutes
with Lie algebra representations. In addition, if $f=f_s+f_n$ is
the Cartan decomposition of an endomorphism $f$ of a finite
dimensional vector space $V$ and $v\in V$ satisfies $fv=0$, then
$f_sv=f_nv=0$, as the reader can check putting $f$ in Jordan form.
This proves the claim.

We want to prove that $\xi_s\in \H^0(E(\zlie))$ and that
$\xi_n=0$. We will need for that the following lemma.

\begin{lemma}
\label{lemma:centralizers} Let $s\in\hclie$ be a semisimple
element. There exists some $h\in \HC$ such that:
\begin{enumerate}
\item if we write $u=\Ad(h^{-1})(s)=h^{-1}sh=u_r+\imag u_i$ with
$u_r,u_i\in\hlie$, then $[u_r,u_i]=0$; \item there exists an
element $a\in\hlie$ such that
$$\Ker \ad(s)=\Ad(h)(\Ker\ad(u_r)\cap \Ker\ad(u_i))
=\Ad(h)\Ker\ad(a).$$
\end{enumerate}
\end{lemma}

\begin{proof}
Using the decomposition $\hclie=\hlie\oplus\imag\hlie$ we define a
real valued scalar product on $\hclie$ as follows: given
$u_r+\imag u_i,v_r+\imag v_i\in\hclie$ we set $\la u_r+\imag
u_i,v_r+\imag v_i\ra_{\RR}:=-\la u_r,v_r\ra-\la u_i,v_i\ra$. The
bilinear pairing $\la\;,\,\ra$ restricted to $\hlie$ is negative
definite, so the pairing $\la\;,\,\ra_{\RR}$ is positive definite on
the whole $\hclie$ and hence the function
$\|\cdot\|^2:\hclie\to\RR$ defined by $\|s\|^2:=\la s,s\ra_{\RR}$
is proper. Let $\OOO_s$ be the adjoint orbit of $s$. Since $s$ is
semisimple, $\OOO_s$ is a closed subset of $\hclie$, and hence the
function $\|\cdot\|^2:\OOO_s\to\RR$ attains its minimum at some
point $u=u_r+\imag u_i\in\OOO_s$. That $u$ minimizes $\|\cdot\|^2$
on its adjoint orbit means that for any $v\in\hclie$ we have $\la
u,[v,u]\ra_{\RR}=0$, since we can identify $T_u\OOO_s=\{[v,u]\mid
v\in\hclie\}$. Now we develop for any $v=v_r+\imag v_i$, using the
invariance of $\la\;,\,\ra$ and Jacobi rule:
\begin{align*}
0 &= \la u_r+\imag u_i, [u_r+\imag u_i,v_r+\imag v_i]\ra_{\RR} \\
&= \la u_r+\imag u_i,
([u_r,v_r]-[u_i,v_i])+\imag([u_i,v_r]+[u_r,v_i])\ra_{\RR} \\
&=-\la u_r,[u_r,v_r]-[u_i,v_i]\ra -\la u_i,[u_i,v_r]+[u_r,v_i]\ra
\\ &=\la u_r,[u_i,v_i]\ra-\la u_i,[u_r,v_i]\ra \\
&= -2\la [u_i,u_r],v_i\ra.
\end{align*}
Since this holds for any choice of $v$, it follows that
$[u_i,u_r]=0$. So the endomorphisms $\ad(u_i)$ and $\ad(u_r)$
commute and hence diagonalize in the same basis with purely
imaginary eigenvalues (because they respect the pairing
$\la\cdot,\cdot\ra_{\RR}$). Hence $\Ker\ad(u)=\Ker\ad(u_r+\imag
u_i)=\Ker(\ad(u_r)+\imag\ad(u_i))=\Ker\ad(u_r)\cap\Ker\ad(u_i)$.
Since $u_r$ and $u_i$ commute, they generate a torus $T_u\subset
H$. Take $h$ such that $u=\Ad(h^{-1})(s)$ and choose $a\in\hlie$
such that the closure of $\{e^{ta}\mid t\in\RR\}$ is equal to
$T_u$. Then $\Ker\ad(a)=\Ker\ad(u_r)\cap\Ker\ad(u_i)$, so the
result follows.
\end{proof}

We now prove that $\xi_s$ is central. Let $u=u_r+\imag
u_i=h^{-1}y_s h$ be the element given by the previous lemma such
that $[u_r,u_i]=0$. Let $\psi_s:E\to\hclie$ be the antiequivariant
map corresponding to $\xi_s\in \H^0(E(\hclie))$, whose image
coincides with the adjoint orbit $\OOO_s$. Define $E_0=\{e\in
E\mid \psi_s(e)=u\}\subset E.$ Then $E_0$ defines a reduction of
the structure group of $E$ to the centralizer of $u$, which we
denote by $\HC_0=\{g\in\HC\mid \Ad(g)(u)=u\}$. Define the
subgroups $P^{\pm}=\{g\in\HC\mid e^{\pm \imag tu_i} g e^{\mp
t\imag u_i}\text{ is bounded as $t\to\infty$ }\}\subset\HC.$ By
(\ref{item:stand-parab-2}) in Lemma \ref{lemma:stand-parab},
$P^{\pm}$ are parabolic subgroups and $L_{u_i}=P^+\cap
P^-=\{g\in\HC\mid \Ad(g)(u_i)=u_i\}$ is a common Levi subgroup of
$P^+$ and $P^-$. By (\ref{item:stand-parab-preliminar}) in Lemma
\ref{lemma:stand-parab}, $\HC_0$ is a connected subgroup of $\HC$,
so by the same argument as in the end of the proof of Lemma
\ref{lemma:centralizers} we can identify $\HC_0$ with
$\{g\in\HC\mid \Ad(g)(u_i)=u_i, \ \Ad(g)(u_r)=u_r\}$. This implies
that $\HC_0\subset L_{u_i}$, hence $E_0$ induces a reduction
$\sigma^+$ (resp. $\sigma^-$) of the structure group of $E$ to
$P^+$ (resp. $P^-$). One the other hand, if $\chi$ corresponds to
$\imag u_i$ via the isomorphism $(\zlie\oplus\clie)^*\simeq
\zlie\oplus\clie$ (so that $s_{\chi}=\imag u_i$), then $\chi$ is
antidominant for $P^+$ and $-\chi$ is antidominant for $P^-$.

Let $\phi:E^L\to B$ be the antiequivariant map corresponding to
$\varphi$. Since $\rho(\xi_s)(\varphi)=0$ we have
$\rho(u)\phi(e)=0$ for any $e\in E_0$. Let $v\in B$ be any
element. Since $u_i$ and $u_r$ commute, the vectors $\rho(e^{\imag
t u_i})v$ are uniformly bounded as $t\to\infty$ if and only if the
vectors $\rho(e^{t u_r})\rho(e^{\imag t u_i})v=\rho(e^{tu})v$ are
bounded. It follows that $\varphi$ belongs both to
$\H^0(E(B)_{\sigma^+,\chi}^-\otimes L)$ and to
$\H^0(E(B)_{\sigma^-,-\chi}^-\otimes L)$. Applying the
$\alpha$-stability condition we deduce that
$$\deg E(\sigma^+,\chi)-\la\alpha,\chi\ra\geq
0,\qquad\text{and}\qquad \deg
E(\sigma^-,-\chi)-\la\alpha,-\chi\ra\geq 0.$$ These inequalities,
together with $\deg E(\sigma^+,\chi)-\la\alpha,\chi\ra=-(\deg
E(\sigma^-,-\chi)-\la\alpha,-\chi\ra)$, imply that $\deg
E(\sigma,\chi)-\la\alpha,\chi\ra=0$. Since we assume that
$(E,\varphi)$ is $\alpha$-stable, such a thing can only happen if
$\chi$, and hence any element in the image of $\psi_s$, is
central.

Finally, we prove that $\xi_n=0$ proceeding by contradiction.
Since the set of nilpotent elements $\hclie_n\subset\hclie$
contains finitely many adjoint orbits, which are locally closed in
the Zariski topology, and since $\xi_n$ is algebraic, there exists
a Zariski open subset $U\subset X$ and an adjoint orbit
$\OOO_n\subset\hclie_n$ such that $\xi_n(x)\in\OOO_n$ for any
$x\in U$. Assume that $\xi_n(x)\neq 0$ for $x\in U$ (otherwise
$\xi_n$ vanishes identically). Consider for any $x\in U$ the
weight filtration of the action of $\ad(\xi_n(x))$ on
$E(\hclie)_x$: $$\dots\subset W_x^{-k}\subset
W_x^{-k+1}\subset\dots\subset W_x^{k-1}\subset
W_x^k\subset\dots,$$ which is uniquely defined by the conditions:
$\ad(\xi_n(x))(W_x^j)\subset W_x^{j-2}$,
$\ad(\xi_n(x))^{j+1}(W_x^{j})=0$ and the induced map on graded
spaces $\Gr\ad(\xi_n(x))^j:\Gr W_x^{j}\to \Gr W_x^{-j}$ is an
isomorphism. As $x$ moves along $U$ the spaces $W_x^j$ give rise
to an algebraic filtration of vector bundles $\dots\subset
W_U^{-k}\subset W_U^{-k+1}\subset\dots\subset W_U^{k-1}\subset
W_U^k\subset\dots\subset E(\hclie)|_U.$ By the properness of the
Grassmannian of subspaces of $\hclie$ these vector bundles extend
to vector bundles defined on the whole $X$
\begin{equation}
\label{eq:weight-filtration} \dots\subset W^{-k}\subset
W^{-k+1}\subset\dots\subset W^{k-1}\subset W^k\subset\dots\subset
E(\hclie)
\end{equation}
and the induced map between graded bundles $\Gr\ad(\xi_n)^j:\Gr
W^j\to \Gr W^{-j}$ is an isomorphism away from finitely many
points. This implies that
\begin{equation}
\label{eq:grau-decreix} \deg \Gr W^j\leq \deg \Gr W^{-j}.
\end{equation}
By Jacobson--Morozov's theorem the weight filtration
(\ref{eq:weight-filtration}) induces a reduction $\sigma$ of the
structure group of $E$ to a parabolic subgroup $P\subset\HC$ (the
so-called Jacobson--Morozov's parabolic subgroup associated to the
nilpotent elements in the image of $\xi_n|_U$), and there exists
an antidominant character $\chi$ of $P$ such that $\ad(s_{\chi})$
preserves the weight filtration and induces on the graded piece
$\Gr W^j$ the map given by multiplication by $j$.

The subbundle $E(B)_{\sigma,\chi}^-\otimes L\subset E(B)\otimes L$
can be identified with the piece of degree $0$ in the weight
filtration on $E(B)\otimes L$ induced by the nilpotent
endomorphism $\rho(\xi_n)$. Since $\rho(\xi_n)(\phi)=0$, we have
$\phi\in\H^0(E(B)_{\sigma,\chi}^-\otimes L)$ (the kernel of a
nonzero nilpotent endomorphism is included in the piece of degree
zero of the weight filtration). Hence, by $\alpha$-stability,
$\deg(E)(\sigma,\chi)-\la\alpha,\chi\ra$ has to be positive. On
the other hand, the character $\chi$ can be chosen to be
perpendicular to $\zlie$, so by (\ref{item:comparing-degree}) in
Lemma \ref{lemma:parabolic-representation} we have
$$\deg E(\sigma,\chi)-\la\alpha,\chi\ra=
\sum_{j\in\ZZ} j\deg\Gr W^j.$$ By (\ref{eq:grau-decreix}) this is
$\leq 0$, thus contradicting the stability of $(E,\varphi)$.
\end{proof}

\subsection{Jordan--H\"older reduction} \label{ss:Jordan-Holder}

In this subsection we associate to each $\alpha$-polystable pair
$(E,\varphi)$ an $\alpha$-stable pair for a different group. This is
accomplished by picking an appropriate subgroup $H'\subset H$ (defined
as the centralizer of a torus in $H$) and by choosing a reduction of
the structure group of $E$ to ${H'}^{\CC}$. The resulting new pair is
called the Jordan--H\"older reduction of $(E,\varphi)$. It is
constructed using a recursive procedure in which certain choices are
made, and the main result of this subsection (see Proposition
\ref{prop:JH-unique}) is the proof that the resulting reduction is
canonical up to isomorphism.


Let $G'\subset G$ be an inclusion of complex connected Lie
subgroup with Lie algebras $\glie'\subset\glie$. Assume that the
normalizer $N_{G}(\glie')$ of $\glie'$ in $G$ is equal to $G'$.
Suppose that $E$ is a holomorphic principal $G$-bundle.

\begin{lemma}
\label{lemma:reduction-normalizer} The holomorphic reductions of
the structure group of $E$ to $G'$ are in bijection with the
holomorphic subbundles $F\subset E(\glie)$ of Lie subalgebras
satisfying this property:
\begin{quotation}
for any $x\in X$ and trivialization $E_x\simeq G$, the fiber
$F_x$, which we identify to a subspace of $\glie$ via the induced
trivialization $E(\glie)_x\simeq\glie$, is conjugate to $\glie'$.
\end{quotation}
\end{lemma}
\begin{proof}
Let $d=\dim\glie'$ and let $\Gr_d(\glie)$ denote the Grassmannian
of complex $d$-subspaces inside $\glie$. Let
$\OOO_{\glie'}=\{\Ad(h)(\glie')\mid h\in G\}\subset \Gr_d(\glie)$.
By assumption there is a biholomorphism $\OOO_{\glie'}\simeq
G/G'$. Furthermore, the set of vector bundles $F\subset E(\glie)$
satisfying the condition of the lemma is in bijection with the
holomorphic sections of $E(\OOO_{\glie'})$, so the result follows.
\end{proof}

We now apply this principle to a particular case. Let $P\subset
\HC$ be a parabolic subgroup, let $L\subset P$ be a Levi subgroup
and let $U\subset P$ be the unipotent radical. Denote $\ulie=\Lie
U$, $\plie=\Lie P$ and $\llie=\Lie L$. The adjoint action of $P$
on $\plie$ preserves $\ulie$ and using the standard projection
$P\to P/U\simeq L$ (see Section
\ref{ss:standard-parabolic-subgroups} and recall that $P$ is
isomorphic to $P_A$ for some choice of $A$) we make $P$ act
linearly on $\llie$ via the adjoint action. Hence $P$ acts
linearly on the exact sequence $0\to\ulie\to\plie\to\llie\to 0$.
We claim that $N_P(\llie)=L$. To check this we identify $P$ (up to
conjugation) with some $P_A$, then use (\ref{eq:def-llie-A}) and
(\ref{eq:def-ulie-A}) together with the surjectivity of the
exponential map $\ulie_A\to U_A$ to deduce that no nontrivial
element of $U$ normalizes $\llie$, and finally use the
decomposition $P=LU$.

\begin{lemma}
\label{lemma:parabolic-levi} Suppose that $E_{\sigma}$ is a
holomorphic principal $P$-bundle. The reductions of the structure
group of $E_{\sigma}$ from $P$ to $L\subset P$ are in bijection
with the splittings of the exact sequence of holomorphic vector
bundles
\begin{equation}
\label{eq:ex-seq-Lie-bdl} 0\to E_{\sigma}(\ulie)\to
E_{\sigma}(\plie)\to E_{\sigma}(\llie)\to 0
\end{equation}
given by holomorphic maps $E_{\sigma}(\llie)\to E_{\sigma}(\plie)$
which are fiberwise morphisms of Lie algebras.
\end{lemma}
\begin{proof}
Since $N_P(\llie)=L$, we may use Lemma
\ref{lemma:reduction-normalizer} with $G=P$ and $G'=L$. The
subalgebras $\glie'\subset\plie$ which are conjugate to $\plie$
are the same as the images of sections $\llie\to\plie$ of the
exact sequence $0\to\ulie\to\plie\to\llie\to 0$ which are
morphisms of Lie algebras. Hence the vector subbundles $F\subset
E(\plie)$ satisfying the requirements of Lemma
\ref{lemma:reduction-normalizer} can be identified with the images
of maps $E(\llie)\to E(\plie)$ which give a section of the
sequence (\ref{eq:ex-seq-Lie-bdl}) and which are fiberwise a
morphism of Lie algebras.
\end{proof}


Suppose that $(E,\varphi)$ is a $\alpha$-polystable pair which is
not $\alpha$-stable. By Proposition
\ref{proposition:stable-simple} there exists a semisimple non
central infinitesimal automorphism $s\in\aut^{ss}(E,\varphi)$. The
splitting $\hclie=\zlie\oplus\hclie_s$ (recall that
$\hclie_s=[\hclie,\hclie]$ is the semisimple part) is invariant
under the adjoint action of $\HC$ (which is connected by
assumption) hence we have
$\H^0(E(\hclie))=\H^0(E(\zlie))\oplus\H^0(E(\hclie_s))$ so
projecting to the second summand we can assume that
$s\in\H^0(E(\hclie_s))$.

As shown in the proof of Proposition
\ref{proposition:stable-simple}, the image of $s$ is contained in
an adjoint orbit in $\hclie$ which contains an element
$u=u_r+\imag u_i$ such that $u_r,u_i$ are commuting elements of
$\hlie$. Let $a\in\hlie_s=[\hlie,\hlie]$ be an infinitesimal
generator of the torus generated by $u_r$ and $u_i$ and let
$\HC_1$ be the complexification of $H_1:=Z_H(a)=\{h\in H\mid
\Ad(h)(a)=a\}$. Let $\psi_s:E\to\hclie$ be the antiequivariant map
corresponding to the section $s$. Then
$$E_1=\{e\in E\mid \psi_s(e)=u\}\subset E$$
is a $\HC_1$-principal bundle, which defines a reduction of the
structure group of $E$. We say that the pair $(E_1,\HC_1)$ is the
{\bf reduction of $(E,\HC)$ induced by $s$ and $u$}.

Define $B_1=\{v\in B\mid \rho(a)(v)=0\}$. The restriction of
$\rho$ to $H_1$ preserves $B_1$, so we have a subbundle
$E_1(B_1)\subset E_1(B)\simeq E(B)$. Let $\phi:E^L\to B$ be the
antiequivariant map inducing the section $\varphi\in
\H^0(E(B)\otimes L)$ (see Section
\ref{ss:principal-bundles-parabolic}). By the definition of the
infinitesimal automorphisms, for any $(e,l)\in E_1^L$ we have
$\rho(\psi_s(e))\phi(e,l)=0$. Now $\rho(\psi_s(e))=\rho(u_r+\imag
u_i)=\rho(u_r)+\imag\rho(u_i)$. Since $\rho$ restricted to $H$ is
Hermitian, $\rho(u_r)$ and $\rho(u_i)$ have purely imaginary
eigenvalues, and since $[\rho(u_r),\rho(u_i)]=0$ it follows that
$$\rho(\psi_s(e))\phi(e,l)=0\quad\Longleftrightarrow\quad
\rho(u_r)\phi(e,l)=\rho(u_i)\phi(e,l)=0
\quad\Longleftrightarrow\quad \rho(a)\phi(e,l)=0$$ for any
$(e,l)\in E^L$. This implies that $\phi(E_1^L)\subset B_1$, and
consequently $\varphi$ lies in the subbundle $E_1(B_1)\otimes
L\subset E(B)\otimes L$. To stress this fact we rename $\varphi$
with the symbol $\varphi_1$. To sum up: assuming that
$(E,\varphi)$ is $\alpha$-polystable but not $\alpha$-stable we
have obtained a subgroup $H_1=Z_H(a)\subset H$, a $H_1$-invariant
subspace $B_1\subset B$, and a new pair $(E_1,\varphi_1)$, where
$E_1$ is a $\HC_1$ principal bundle and $\varphi_1\in
\H^0(E_1(B_1)\otimes L)$. We denote the Lie algebras of $H_1$ and
its complexification by $\hlie_1$ and $\hclie_1$.

\begin{proposition}
The pair $(E_1,\phi_1)$ is $\alpha$-polystable.
\end{proposition}
\begin{proof}
Since $H_1$ is the centralizer of $a$ and $\alpha$ belongs to the
center of $\hclie$, we have $\alpha\in\hclie_1$. Hence the
statement of the proposition makes sense. We first prove that
$(E_1,B_1)$ is $\alpha$-semistable. Let $P_1\subset \HC_1$ be a
standard parabolic subgroup. By (\ref{item:stand-parab-1}) in
Lemma \ref{lemma:stand-parab} there is some $s\in\imag\hlie_1$
(satisfying $s=s_{\chi}$ for an appropriate antidominant character
$\chi$ of $P_1$) such that $P_1=\{g\in\HC_1\mid
e^{ts}ge^{-ts}\text{ is bounded as $t\to\infty$ }\}$. Since
$\imag\hlie_1\subset\imag\hlie$ it makes sense to define $P=\{g\in
\HC\mid e^{ts}ge^{-ts}\text{ is bounded as $t\to\infty$ }\}$,
which is a parabolic subgroup of $\HC$, and clearly $P_1\subset
P$. Hence, any reduction $\sigma_1$ of the structure group of
$E_1$ to $P_1$, say $(E_1)_{\sigma_1}\subset E_1$, gives
automatically a reduction $\sigma$ of the structure group of $E$
to $P$, specified by
$E_{\sigma}=(E_1)_{\sigma_1}\times_{P_1}P\subset
(E_1)_{\sigma_1}\times_{P_1}\HC=E$. Furthermore, any antidominant
character $\chi\in\imag\hlie$ of $P_1$ is an antidominant
character of $P$, and there is an equality
$\deg(E_1)(\sigma_1,\chi)=\deg(E)(\sigma,\chi)$. Finally, if the
section $\varphi_1$ belongs to
$\H^0(E_1(B_1)_{\sigma_1,\chi}^-\otimes L)$, then it also belongs
to $\H^0(E(B)_{\sigma,\chi}^-\otimes L)$. All this implies that
$(E_1,\varphi_1)$ is $\alpha$-semistable.

To prove that $(E_1,\varphi_1)$ is $\alpha$-polystable it remains
to show that if the reduction $\sigma_1$ and $\chi$ have been
chosen so that $\deg(E_1)(\sigma_1,\chi)-\la\alpha,\chi\ra=0$,
then there is a holomorphic reduction $\sigma_{L_1}$ of the
structure group of $(E_1)_{\sigma_1}$ to the Levi
$L_1=\{g\in\HC_1\mid \Ad(g)(s)=s\}$ such that
\begin{equation}
\label{eq:condition-on-phi-1}
\varphi_1\in\H^0(E(B_1)_{\sigma_{L_1},\chi}^0\otimes L).
\end{equation}
Define $L=\{g\in\HC\mid \Ad(g)(s)=s\}$, which is a Levi subgroup
of $P$, let $U_1\subset P_1$ and $U\subset P$ be the unipotent
radicals, and denote the corresponding Lie algebras by
$\ulie_1=\Lie U_1$, $\plie_1=\Lie P_1$, $\llie_1=\Lie L_1$,
$\ulie=\Lie U$, $\plie=\Lie P$, $\llie=\Lie L$. By Lemma
\ref{lemma:parabolic-levi} it suffices to check that there exists
a bundle morphism $w_1:(E_1)_{\sigma_1}(\llie_1)\to
(E_1)_{\sigma_1}(\plie_1)$ given fiberwise by morphisms of Lie
algebras, defining a splitting of the exact sequence
\begin{equation}
\label{eq:extension-class-poly-0} 0\to
(E_1)_{\sigma_1}(\ulie_1)\to (E_1)_{\sigma_1}(\plie_1)\to
(E_1)_{\sigma_1}(\llie_1)\to 0.
\end{equation}

Let $T\subset H$ be the closure of $\{e^{ta}\mid t\in\RR\}$, which
is a torus. Denote by $T^{\vee}=\Hom(T,S^1)$ the group of
characters of $T$. We have decompositions
$$\ulie=\bigoplus_{\eta\in T^{\vee}}\ulie_{\eta},\qquad
\plie=\bigoplus_{\eta\in T^{\vee}}\plie_{\eta},\qquad
\llie=\bigoplus_{\eta\in T^{\vee}}\llie_{\eta},$$ and since the
elements of $\HC_1$ fix $a$, the action of $\HC_1$ on $\ulie$,
$\plie$ and $\llie$ respects the splittings above. It follows that
we have a commutative diagram with exact rows $$\xymatrix{ 0\ar[r]
& E_\sigma(\ulie) \ar[r]\ar[d]^{\simeq} & E_\sigma(\plie)
\ar[r]\ar[d]^{\simeq} &
E_\sigma(\llie)\ar[d]^{\simeq} \ar[r] & 0 \\
0\ar[r] & (E_1){\sigma_1}(\ulie) \ar[r]\ar@{=}[d] &
(E_1)_{\sigma_1}(\plie) \ar[r]\ar@{=}[d] & (E_1)_{\sigma_1}(\llie)
\ar[r]\ar@{=}[d] &
0 \\
0\ar[r] & \bigoplus_{\eta\in
T^{\vee}}(E_1){\sigma_1}(\ulie_{\eta}) \ar[r] & \bigoplus_{\eta\in
T^{\vee}}(E_1)_{\sigma_1}(\plie_{\eta}) \ar[r] &
\bigoplus_{\eta\in T^{\vee}}(E_1)_{\sigma_1}(\llie_{\eta}) \ar[r]
& 0 }$$
Taking in the bottom row the summands corresponding to the trivial
character $\eta=0$ (the constant representation $T\to\{1\}\in
S^1$) we get the exact sequence (\ref{eq:extension-class-poly-0}).
By hypothesis the pair $(E,\varphi)$ is $\alpha$-polystable, so
there is a section $v:E_{\sigma}(\llie)\to E_{\sigma}(\plie)$ of
the top row, given fiberwise by morphisms of Lie algebras. Using
the isomorphisms and equalities in the diagram, this gives rise to
a section
$$w:\bigoplus_{\eta\in T^{\vee}}(E_1)_{\sigma_1}(\llie_{\eta})\to
\bigoplus_{\eta\in T^{\vee}}(E_1)_{\sigma_1}(\plie_{\eta})$$ of
the bottom row. Then $w=(w_{\eta\mu})_{\eta,\mu\in T^{\vee}}$,
where $w_{\eta\mu}:(E_1)_{\sigma_1}(\llie_{\eta})\to
(E_1)_{\sigma_1}(\plie_{\mu})$, and one checks that $w_1:=w_{00}$
is fiberwise a morphism of Lie algebras and that it gives the
desired splitting of the sequence
(\ref{eq:extension-class-poly-0}). To check
(\ref{eq:condition-on-phi-1}) we proceed as follows. First note
that $s_{\chi}$ belongs both to the center of $\llie_1$ and
$\llie$, hence it defines holomorphic sections
$s_{\sigma_1,\chi}\in\H^0((E_1)_{\sigma_1}(\llie_1))$ and
$s_{\sigma,\chi}\in\H^0(E_\sigma(\llie))$. Condition
(\ref{eq:condition-on-phi-1}) is equivalent to
\begin{equation}
\label{eq:condition-on-phi-2}
\rho(w_1(s_{\sigma_1,\chi}))(\varphi)=0
\end{equation}
(note that $(E_1)_{\sigma_1}(\plie_1)$ is a subbundle of
$(E_1)_{\sigma_1}(\hclie_1)\simeq E_1(\hclie_1)$, hence it acts
fiberwise on $E(B)\otimes L$). To prove this equality, we use
again the hypothesis that $(E,\varphi)$ is $\alpha$-polystable,
which implies that $\varphi\in\H^0(E(B)_{\sigma_L,\chi}^-\otimes
L)$, where $\sigma_L$ is the reduction specified by $w$. This is
equivalent to $\rho(w(s_{\sigma,\chi}))(\varphi)=0$, and this
implies (\ref{eq:condition-on-phi-2}) because
$s_{\chi}\in\llie_0\subset\bigoplus_{\eta\in
T^{\vee}}\llie_{\eta}$.
\end{proof}


Let $(E,\varphi)$ be a $\alpha$-polystable pair. Iterating the
procedure described in the previous subsection as many times as
possible we obtain a sequence of groups $H=H_0\supset H_1\supset
H_2\supset\dots$ and elements
$a_j\in(\hlie_{j-1})_s=[\hlie_{j-1},\hlie_{j-1}]$ such that
$H_j=Z_{H_{j-1}}(a_j)$, vector subspaces $B=B_0\supset B_1\supset
B_2\supset\dots$, and $\alpha$-polystable pairs
$(E,\varphi)=(E_0,\varphi_0),\ (E_1,\varphi_1),\dots$, where $E_j$
is a $\HC_j$-principal bundle over $X$ and contained in $E_{j-1}$,
and $\varphi_j\in \H^0(E_j(B_j)\otimes L)$. Since $\dim H_j<\dim
H_{j-1}$, this process has to eventually stop at some pair, say
$(E_r,\varphi_r)$, which will necessarily be $\alpha$-stable. We
say that $(E_r,\varphi_r,H_r,B_r)$ is the {\bf Jordan--H\"older}
reduction of $(E,\varphi,H,B)$. To justify this terminology we
need to prove that the construction is independent of the choices
made in the process. Note that the elements in the sequence
$\{a_0,a_1,\dots,a_l\}$ all belong to the initial Lie algebra
$\hlie$ and they commute pairwise. Hence they generate a torus
$T\subset H$, the closure of the set $\{\exp\sum t_ja_j\mid
t_0,\dots,t_l\in\RR\}$, and $H_l$ is the centralizer in $H$ of
$T_{(E,\varphi)}$. With this in mind, the following proposition
implies the uniqueness of the Jordan--H\"older reduction.

Let $H_s\subset H$ be the connected Lie subgroup whose Lie algebra
is $\hlie_s=[\hlie,\hlie]$.

\begin{proposition}
\label{prop:JH-unique} Let $(E,\varphi)$ be a $\alpha$-polystable
pair. Suppose that $T',T''\subset H_s$ are tori, and define $H'$
(resp. $H''$) to be the centralizer in $H$ of $T'$ (resp. $T''$).
Let $B'$ (resp. $B''$) be the fixed point set of the action of
$T'$ (resp. $T''$) on $B$, and assume that there are reductions
$E'\subset E$ (resp. $E''\subset E$) of the structure group of $E$
to ${H'}^{\CC}$ (resp. ${H''}^{\CC}$). Let $\phi:E^L\to B$ the
equivariant map corresponding to $\varphi$. Assume that
$\phi(E'^L)\subset B'\otimes L$ and $\phi(E''^L)\subset B''\otimes
L$. Denote by $\varphi'\in \H^0(E'(B')\otimes L)$ and
$\varphi''\in \H^0(E''(B'')\otimes L)$ the induced sections.
Finally, suppose that both $(E',\varphi')$ and $(E'',\varphi'')$
are $\alpha$-stable. Then there is some $g\in\HC$ such that
$H'^{\CC}= g^{-1}(H''^{\CC})g$, $E'=E''g$,
$T'^{\CC}=g^{-1}(T''^{\CC})g$ and $B'=\rho(g^{-1})B''$.
\end{proposition}

Before proving Proposition \ref{prop:JH-unique} we state and prove
two auxiliary lemmas.

\begin{lemma}
\label{lemma:commuting-central} Let $u',u''\in\hlie$ and let
$s',s''\in \H^0(E(\hclie))$ be sections such that $s'(x)$ (resp.
$s''(x)$) is conjugate to $\imag u'$ (resp. $\imag u''$) for any
$x\in X$. Let $(E',H'^{\CC})$ (resp. $(E'',H''^{\CC})$) be the
reductions of $(E,\HC)$ induced by $s'$ and $\imag u'$ (resp.
$s''$ and $\imag u''$).
\begin{enumerate}

\item Assume that $[s',s'']=0$. Let $\hlie''^{\CC}$ be the Lie
algebra of $H''^{\CC}$. Then we can naturally identify $s'$ with a
section of $E''(\hlie''^{\CC})$.

\item Let $\zlie''$ be the center of $\hlie''^{\CC}$. If
$s'\in\H^0(E''(\zlie''))$ then there is some $h\in\HC$ such that
$E''\subset E'h$ as subsets of $E$.

\end{enumerate}
\end{lemma}
\begin{proof}
Let $\psi',\psi'':E\to\hclie$ be the antiequivariant maps
corresponding to $s',s''$. The condition $[s',s'']=0$ implies that
for any $e\in E$ the elements $\psi'(e),\psi''(e)\in\hclie$
commute. Since $E''=(\psi'')^{-1}(\imag u'')$, this implies that,
for any $e\in E''$, $\psi'(e)$ commutes with $\imag u''$, so
$\psi'(e)$ belongs to $\hlie''^{\CC}$. This proves (1). We now
prove (2). First observe that, being a centralizer of a semisimple
element in $\hclie$, $H''^{\CC}$ is connected (see e.g. Theorem
13.2 in \cite{borel:1956}). Hence, the adjoint action of
$H''^{\CC}$ on $\hlie''^{\CC}$ fixes any element in $\zlie''$.
Take some element $e\in E''$. By hypothesis, there is some
$h\in\HC$ such that $\psi'(e)=\Ad(h^{-1})(\imag u')$, so $e\in
E'h$. The condition $s'\in\H^0(E''(\zlie''))$ implies that
$\psi'(e)\in\zlie''$ so, by the previous observation, for any
$g\in H''^{\CC}$ we have $\psi'(eg)=\Ad(g^{-1})\Ad(h^{-1})(\imag
u')=\Ad(h^{-1})(\imag u')$, hence $eg\in E'h$. It follows that
$E''\subset E'h$.
\end{proof}

For any $u\in\hlie$ we denote by $T_u\subset H$ the torus
generated by $u$, i.e., the closure of $\{\exp {tu}\mid
t\in\RR\}$, and $T_u^{\CC}$ denotes the complexification of $T_u$.

\begin{lemma}
\label{lemma:de-tal-toro-tal-astilla} Let
$u',u''\in\hlie_s=[\hlie,\hlie]$ and let $H'^{\CC}$ (resp.
$H''^{\CC}$) be the complexification of the centralizer $Z_H(u')$
(resp. $Z_H(u'')$). If there is some $g\in\HC$ such that
$H'^{\CC}=g^{-1} (H''^{\CC}) g$ then
$T_{u'}^{\CC}=g^{-1}T_{u''}^{\CC}g$.
\end{lemma}
\begin{proof}
The center of $\hlie'^{\CC}$ is $\zlie\oplus\Lie T_{u'}^{\CC}$,
and the sum is direct because $u'$ is assumed to belong to
$\hlie_s$. Similarly, the center of $\hlie''^{\CC}$ is
$\zlie\oplus\Lie T_{u''}^{\CC}$. Since $\HC$ is connected, its
adjoint action on $\zlie$ is trivial, and hence taking the center
of the Lie algebra in each side of the equality
$T_{u'}^{\CC}=g^{-1}T_{u''}^{\CC}g$ we deduce that $\Lie
T_{u'}^{\CC}=g^{-1}(\Lie T_{u''}^{\CC}) g$. This implies the
equality between the complexified tori.
\end{proof}


We now prove Proposition \ref{prop:JH-unique}.

\begin{proof}
Let $u',u''\in\hlie_s$ satisfy $T'=T_{u'}$ and $T''=T_{u''}$. The
existence of reductions of $E$ to the centralizers of $u'$ and
$u''$ gives rise to sections $s',s''\in
\aut^{ss}(E,\varphi)\subset \H^0(E(\hclie))$ such that $s'(x)$
(resp. $s''(x)$) is conjugate to $\imag s'$ (resp. $\imag s''$)
for any $x\in X$.

If $[s',s'']=0$ then by (1) Lemma \ref{lemma:commuting-central} we
can view $s'\in\aut^{ss}(E'',\varphi'')$ and
$s''\in\aut^{ss}(E',\varphi')$. Since by assumption
$(E'',\varphi'')$ and $(E',\varphi')$ are $\alpha$-stable, by
Proposition \ref{proposition:stable-simple} we deduce that $s'$ is
central in the centralizer of $s''$ and vice-versa. By (2) in
Lemma \ref{lemma:commuting-central} there exist $g,h\in\HC$ such
that $E'\subset E''g$ and $E''\subset E'h$. This implies that
$E'\subset E''g\subset E'hg$, but $E'\subset E'hg$ clearly implies
that $E'=E'hg$, which combined with the previous chain of
inclusions gives $E'=E''g$. It then follows that
$H'^{\CC}=g^{-1}(H''^{\CC})g$. By Lemma
\ref{lemma:de-tal-toro-tal-astilla} we have
$T_{s'}^{\CC}=g^{-1}T_{s''}^{\CC}g$. Finally, since the fixed
point set of $T_{s'}^{\CC}$ acting on $B$ coincides with the fixed
point set of $T_{s'}$ (and similarly for $T_{s''}^{\CC}$) we have
$B'=\rho(g^{-1})B''$.

Suppose now that $[s',s'']\neq 0$. There are holomorphic
splittings
\begin{equation}
\label{eq:decom-E-F} E(\hclie)=E_1\oplus\dots\oplus
E_p=F_1\oplus\dots\oplus F_q
\end{equation}
such that $\ad(s')|_{E_j}=\lambda_j\Id_{E_j}$ and
$\ad(s'')|_{F_k}=\mu_k\Id_{F_k}$, where the real numbers
$\lambda_1<\dots<\lambda_p$ (resp. $\mu_1<\dots<\mu_q$) are the
eigenvalues of $\ad(\imag s')$ (resp. $\ad(\imag s'')$).
Define for any $j$ the subbundles $F_{\leq j}=\bigoplus_{k\leq
j}F_k\subset E(\hclie)$ and $E_{\leq j}=\bigoplus_{k\leq
j}E_k\subset E(\hclie)$. Denote by $\pi_k:E(\hclie)\to E_k$ the
projection using the decomposition (\ref{eq:decom-E-F}). Let
$\EEE_{\leq k}$ (resp. $\EEE_k$, $\FFF_{\leq j}$, $\FFF_j$) be the
sheaf of local holomorphic sections of $E_{\leq k}$ (resp. $E_k$,
$F_{\leq j}$, $F_j$). Define for any $j$ the sheaf
$$\FFF_{\leq j}^{\sharp}=\bigoplus_{k=1}^p
\pi_k(\EEE_{\leq k}\cap \FFF_{\leq j}).$$ This is a subsheaf of
the sheaf associated to $E(\hclie)$, and we denote by $F_{\leq
j}^{\sharp}\subset E(\hclie)$ the subbundle obtained by taking the
saturation of $\FFF_{\leq j}^{\sharp}$.

By (1) in Lemma \ref{lemma:parabolic-representation} $s''$ induces
a holomorphic reduction $\sigma''\in\Gamma(E(\HC/P))$ of the
structure group of $E$ to $P=P_{\imag u''}$.

\begin{lemma}
The filtration $F_{\leq 1}^{\sharp}\subset\dots\subset F_{\leq
q}^{\sharp}=E(\hclie)$ also induces a reduction $\sigma^{\sharp}$
of the structure group of $E$ to $P$.
\end{lemma}
\begin{proof}
For any $t\in\RR$ there is a natural fiberwise action of $e^{ts'}$
on $E(\HC/P)$, which allows to define
$e^{ts'}\sigma''\in\Gamma(E(\HC/P))$. For the reader's
convenience, we recall how this is defined. For any $x\in X$ we
can identify $\sigma''(x)$ with an antiequivariant map
$\xi_{\sigma''}:E_x\to\HC/P$ (here $\HC$ acts on the left of
$\HC/P$). Similarly, $s'(x)$ corresponds to a map
$\psi:E_x\to\hclie$ which is antiequivariant and hence satisfies,
for any $f\in E_x$ and $g\in\HC$,
\begin{equation}
\label{eq:exponential-psi-eg} e^{t\psi(fg)}=g^{-1}e^{t\psi(f)}g.
\end{equation}
Then $e^{t s'}\sigma''(x)$ corresponds to the antiequivariant map
$\xi_{e^{t s'}\sigma''}:E_x\to\HC/P$ defined as
$$\xi_{e^{t s'}\sigma''}(f)=e^{t\psi(f)}\xi_{\sigma''}(f)
=\xi_{\sigma''}(f e^{-t\psi(f)}).$$ That $\xi_{e^{ts'}\sigma''}$
is antiequivariant follows from (\ref{eq:exponential-psi-eg}). For
each $x$ the action of $e^{t s'(x)}$ defines on the fiber
$E_x(\HC/P)$ a decomposition in Zariski locally closed
subvarieties $\{\CCC_{x,i}\}$, the Schubert cells. Each
$\CCC_{x,i}$ corresponds to a connected component $C_{x,i}\subset
E_x(\HC/P)$ of the fixed point set of the action of $\{e^{t
s'(x)}\mid t\in\RR\}$ on $E_x(\HC/P)$, and $\CCC_{x,i}$ is the set
of $z\in E_x(\HC/P)$ such that $e^{ts'(x)}z$ converges to
$C_{x,i}$ as $t\to\infty$. Since $s'$ is algebraic and, for any
$x$, $s'(x)$ is conjugate to the same element $\imag u'$, each
$\CCC_i=\bigcup_{x\in X}\CCC_{x,i}$ is a Zariski locally closed
subvariety of $E(\HC/P)$. Since $\sigma''$ is an algebraic section
of $E(\HC/P)$, there is a Zariski open subset $U\subset X$ such
that $\sigma''|_U$ is contained in a unique cell $\CCC_j\subset
E(\HC/P)$. Then for any $x\in U$ the limit
$\sigma^{\sharp}_x:=\lim_{t\to\infty}e^{t s'}\sigma''(x)\in
C_{x,j}\subset \CCC_j$ is well defined, and the filtration
$\{\FFF_{\leq j,x}^{\sharp}\}$ corresponds to $\sigma^{\sharp}_x$.
As $x$ moves along $U$ the elements $\sigma^{\sharp}_x$ describe
an algebraic section $\sigma^{\sharp}_U\in \Gamma(U;E(\HC/P))$.
Finally, $F_{\leq j}^{\sharp}$ results from extending the
reduction $\sigma^{\sharp}_U$ to an algebraic section
$\sigma^{\sharp}\in \Gamma(E(\HC/P))$, which exists and is unique
thanks to the properness of the flag variety $\HC/P$.
\end{proof}

Let $\chi$ be the antidominant character of $P$ corresponding to
$u''$, so that $s_{\chi}=\imag u''$.

\begin{lemma}
We have $\varphi\in \H^0(E(B)_{\sigma^{\sharp},\chi}^-\otimes L)$.
\end{lemma}
\begin{proof}
Let $U\subset X$ denote, as in the preceding lemma, a nonempty
Zariski open subset such that for any $x\in U$ we have
$\sigma^{\sharp}(x)=\lim_{t\to\infty}e^{ts'}\sigma''(x)$. By
continuity, it suffices to prove that for any $x\in U$
\begin{equation}
\label{eq:phi-es-a-dins} \varphi(x)\in
E(B)_{\sigma^{\sharp},\chi}^-\otimes L.
\end{equation}
The vector $\varphi(x)$ corresponds to an antiequivariant map
$\phi:E^L_x\to B$, whereas $\sigma^{\sharp}$ corresponds to an
antiequivariant map $\xi_{\sigma^{\sharp}}:E_x\to\HC/P$. Define
$P^{\sharp}_x=\xi_{\sigma^{\sharp}}^{-1}(P)\subset E_x$. Then
$P^{\sharp}_x$ is an orbit of the action of $P$ on $E_x$ on the
right (which can also be obtained by identifying $E(\HC/P)$ with
the quotient $E/P$). And (\ref{eq:phi-es-a-dins}) is equivalent to
requiring that $\phi(x)$ restricted to $(P^{\sharp}_x)^L$ is
contained in $B^-_{\chi}$. Define for any real $t$ the map
$\xi_{\sigma^t}:E_x\to\HC/P$ as
$\xi_{\sigma^t}(f)=\xi_{\sigma''}(fe^{-t\psi(f)})$, where
$\psi:E_x\to\hclie$ is the antiequivariant map corresponding to
$s'$. Let also $P^t_x$ be $\xi_{\sigma^t}^{-1}(P)$. By the
previous lemma, we have $\xi_{\sigma^\sharp}=\lim_{t\to\infty}
\xi_{\sigma^t}$, so we have $P^{\sharp}_x=\lim_{t\to\infty}P^t_x$
as orbits of $E_x/P$. By continuity, it suffices to check that for
any $t$ the restriction of $\phi(x)$ to $(P^t_x)^L$ is contained
in $B^-_{\chi}$.

Since $s',s''\in\aut(E,\varphi)$, we have
\begin{equation}
\label{eq:rho-fixes-varphi} \rho(e^{ts'})(\varphi)=\varphi
\end{equation}
and we also have $\varphi\in\H^0(E(B)_{\sigma'',\chi}^-\otimes
L)$. Defining $P''_x=\xi_{\sigma''}^{-1}(P)$ this implies that
\begin{equation}
\label{eq-varphi-negative} \phi(g,l)\in B^-_{\chi}\qquad \text {
for any $g\in P''_x$ and $l\in L_x$.}
\end{equation} Assume that $f\in P^t_x$ and
$l\in L_x$. Then
$\xi_{\sigma^t}(f)=\xi_{\sigma''}(fe^{-t\psi(f)})\in P$, so
$fe^{-t\psi(f)}\in P''_x$. Hence
$$\phi(f,l)=\phi(fe^{-t\psi(f)},l)\in B_{\chi}^-,$$
where the equality follows from (\ref{eq:rho-fixes-varphi}) and
the inclusion follows from (\ref{eq-varphi-negative}). This proves
that $\phi(x)$ maps $(P^t_x)^L$ inside $B^-_{\chi}$, so we are
done.
\end{proof}

Hence we can apply the $\alpha$-polystability condition, which in
view of Lemma \ref{lemma:parabolic-representation} and Remark
\ref{rem:parabolic-representation} reads
\begin{equation}
\label{eq:sum-deg-F-sharp} \deg(E)(\sigma^{\sharp},\chi)=\mu_q\deg
F_{\leq q}^{\sharp} +\sum_{j=1}^{q-1}(\mu_j-\mu_{j+1})\deg F_{\leq
j}^{\sharp}\geq 0
\end{equation}
(the $\la\alpha,\chi\ra$ term vanishes because we assume that
$s''$ is orthogonal to the center of $\hlie$). On the other hand,
since $s''\in\aut^{ss}(E,\varphi)$, the same arguments as in the
proof of Proposition \ref{proposition:stable-simple} imply that
\begin{equation}
\label{eq:sum-deg-F} \deg(E)(\sigma'',\chi)=\mu_q\deg F_{\leq q}
+\sum_{j=1}^{q-1}(\mu_j-\mu_{j+1})\deg F_{\leq j}=0.
\end{equation}
An easy computation shows that $\deg\FFF_{\leq j}^{\sharp}=\deg
F_{\leq j}$, whereas in general $\deg\FFF_{\leq j}^{\sharp}\leq
\deg F_{\leq j}^{\sharp}$ with equality if and only if $\FFF_{\leq
j}^{\sharp}=(\FFF_{\leq j}^{\sharp})^{\vee\vee}$, so that in
general $$\deg F_{\leq j}\leq\deg F^{\sharp}_{\leq j}.$$ Since
$\deg F_{\leq q}=\deg\FFF_{\leq q}^{\sharp}=\deg F_{\leq
q}^{\sharp}$ (because $\FFF_{\leq q}$ is equal to the sheaf
associated to $E(\hclie)$) and $\mu_j-\mu_{j+1}<0$ for any $1\leq
j\leq q-1$, we have
$$\deg(E)(\sigma'',\chi)\geq \deg(E)(\sigma^{\sharp},\chi),$$
which combined (\ref{eq:sum-deg-F-sharp}) and (\ref{eq:sum-deg-F})
yields $\deg(E)(\sigma'',\chi)=\deg(E)(\sigma^{\sharp},\chi)=0$.
By the previous comments, this equality implies $\FFF_{\leq
j}^{\sharp}=(\FFF_{\leq j}^{\sharp})^{\vee\vee}$ for any $j$, so
that $\FFF_{\leq j}^{\sharp}$ is the sheaf of local holomorphic
sections of a subbundle $F_{\leq j}^{\sharp}\subset E(\hclie)$.
This has the following consequence: if we define
$\FFF_l^{\sharp}=\bigoplus_k \pi_k(\FFF_l\cap\EEE_{\leq k})$, then
$\FFF_l^{\sharp}$ is also the sheaf of sections of a subbundle
$F_l^{\sharp}\subset E(\hclie)$ and we have $F_{\leq
j}^{\sharp}=\bigoplus_{l\leq j}F_l^{\sharp}$. In particular, we
obtain a decomposition $E(\hclie)=\bigoplus_{l\leq q}
F_l^{\sharp}$. Let $s^{\sharp}=\sum_j\mu_j\Id_{F_j^{\sharp}}\in
\H^0(E(\hclie))$. Then we have $[s',s^{\sharp}]=0$ and furthermore
$s^{\sharp}\in\aut^{ss}(E,\varphi)$. These two properties imply
that $s^{\sharp}\in\aut^{ss}(E',\varphi')$, so by Proposition
\ref{proposition:stable-simple} $s^{\sharp}$ is central in the
centralizer of $s'$. Similarly $s'$ is central in the centralizer
of $s^{\sharp}$, so we can proceed as in the first case and deduce
the statement of the theorem with $s''$ replaced by $s^{\sharp}$.
Reversing the roles of $s'$ and $s''$ we conclude the proof of
Proposition \ref{prop:JH-unique}.
\end{proof}

\subsection{Hitchin-Kobayashi correspondence}
\label{ss:Hitchin-Kobayashi-correspondence}

Choose a Hermitian metric $h_L$, on the complex line bundle $L$,
and denote by $F_L\in\Omega^2(X;\imag\RR)$ the curvature of the
corresponding Chern connection. Suppose that $E_h\subset E$
defines a reduction of the structure group of $E$ from $\HC$ to
$H$. Then the vector bundle $E(B)=E\times_{\HC}B$ can be
canonically identified with $E_h\times_H B$, and hence inherits a
Hermitian structure (obtained from the Hermitian structure on $B$,
which is preserved by $H$). So for any $\varphi\in
\H^0(E(B)\otimes L)$ it makes sense to define
$$\mu_h(\varphi):=\rho^*\left(-\frac{\imag}{2}
\varphi\otimes\varphi^{*_{h,h_L}}\right).$$ Here we identify
$\imag\varphi\otimes\varphi^{*_{h,h_L}}$ with a skew symmetric
section of $\End(E(B)\otimes L)^*=\End(E(B))^*$, hence a section
of $E_h(\ulie(B))^*$. The map $\rho^*:E_h(\ulie(B))^*\to
E_h(\hlie)^*$ is induced by the dual of the infinitesimal action
of $\hlie$ on $B$. Using the isomorphism $\hlie^*\simeq\hlie$
given by the non-degenerate pairing $\la\cdot,\cdot,\ra$ we view
$\mu_h(\varphi)$ as a section of $E_h(\hlie)$.

\begin{theorem}\label{hk-twisted-pairs}.
Let $(E,\varphi)$ be a $\alpha$-polystable pair. There exists a
reduction $h$ of the structure group of $E$ from $\HC$ to $H$,
given by a subbundle $E_h\subset E$, such that
\begin{equation}
\label{eq:Hitchin-Kobayashi}
\Lambda(F_h)+\mu_h(\varphi)=-\imag\alpha,
\end{equation}
where $F_h\in\Omega^2(X;E_h(\hlie))$ denotes the curvature of the
Chern connection on $E$ with respect to $h$ and
$\Lambda:\Omega^2(X)\to\Omega^0(X)$ is the adjoint of wedging with
the volume form on $X$. Furthermore, if $(E,\varphi)$ is
$\alpha$-stable then $h$ is unique. Conversely, if $(E,\varphi)$
is a pair which admits a solution to equation
(\ref{eq:Hitchin-Kobayashi}), then $(E,\varphi)$ is
$\alpha$-polystable.
\end{theorem}

\begin{proof}
Suppose first of all that $(E,\varphi)$ is $\alpha$-stable. Then
by Proposition \ref{proposition:stable-simple} we have
$\aut^{ss}(E,\varphi)=\H^0(E(\zlie))$, so $(E,\varphi)$ is simple
in the sense of Definition~3.8 in
\cite{bradlow-garcia-prada-mundet:2003}.

Hence we can apply
Theorem 4.1 of \cite{bradlow-garcia-prada-mundet:2003} to deduce
the existence and uniqueness of $h$. (Recall that the notion of
$\alpha$-stability given in the present paper coincides with the
one in \cite{bradlow-garcia-prada-mundet:2003} thanks to
(\ref{item:comparing-degree}) in Lemma
\ref{lemma:parabolic-representation}.) If $(E,\varphi)$ is
$\alpha$-polystable but not stable, then we consider the
Jordan--H\"older reduction $(E',\varphi',H',B')$ of
$(E,\varphi,H,B).$ Now the pair $(E',\varphi')$ is simple and we
can proceed as before to get a reduction $h'$ of the structure
group of $E'$ from ${H'}^{\CC}$ to $H'$ satisfying
(\ref{eq:Hitchin-Kobayashi}). But $h'$ also defines a reduction of
the structure group of $E$ from $\HC$ to $H$, by defining
$E_h:=E_{h'}\times_{H'}H\subset E_{h'}\times_{H'}\HC=E$. For this
choice of $h$, equation (\ref{eq:Hitchin-Kobayashi}) still holds.

The proof of the converse is standard. One first proves that if
$(E,\varphi)$ admits a solution to the equations then
$(E,\varphi)$ is $\alpha$-semistable (see for example
\cite{bradlow-garcia-prada-mundet:2003}). To prove
$\alpha$-polystability one can use the same strategy as in the
Hitchin--Kobayashi correspondence for vector bundles.
Namely, assume that $h\in E(\HC/H)$ defines a reduction of the
structure group to $H$, in such a way that equation
(\ref{eq:Hitchin-Kobayashi}) is satisfied. Assume also that
$P\subset\HC$ is a parabolic subgroup, that there is a holomorphic
reduction $\sigma$ of the structure group of $E$ to $P$, an
antidominant character $\chi$ of $P$ such that $\varphi$ is
contained in $E(B)_{\sigma,\chi}^-\otimes L$ and such that
\begin{equation}
\label{eq:grau-zero} \deg(E)(\sigma,\chi)-\la\alpha,\chi\ra=0.
\end{equation}
We want to prove that there is a further reduction $\sigma_L$ of
the structure group of $E$ from $P$ to $L$ and that $\varphi$ is
contained in $E(B)_{\sigma_L,\chi}^0\otimes L$.

Let $E_h\subset E$ be the principal $H$ bundle specified by $h$.
The reduction $\sigma$ corresponds to an antiequivariant map
$\xi:E\to\HC/P$, so that $\xi(f)$ is a parabolic subgroup of $\HC$
for each $f\in E$. Then, using the construction given in Lemma
\ref{lemma:tots-parabolics-P-s} we define an $H$-antiequivariant
map $\psi:E_h\to\imag\hlie$ by setting $\psi(f)=s_{\xi(f),\chi}$
for any $f\in E_h$. The map $\psi$ corresponds to a section of
$E_h(\imag\hlie)$, which we denote by
$$s_{h,\sigma,\chi}\in E_h(\imag\hlie).$$

For details on the following notions the reader can consult
\cite{mundet:2000}. Let $\EE$ be the $C^{\infty}$ $H$-principal
bundle underlying $E_h$, and let $\AAA$ be the set of connections
on $\EE$. Each element of $A\in \AAA$ defines a holomorphic
structure $\ov{\partial}_A$ on $\EE$. Let also $\SSS$ be the space
of smooth sections of $\EE\times_HB\otimes L$, and let $\GGG$ be
the gauge group of $E$. The space $\AAA\times\SSS$ has a natural
structure of infinite dimensional symplectic manifold, with
respect to which the action of $\GGG$ is Hamiltonian and
$(A,\phi)\mapsto
\mu(A,\phi):=\Lambda(F_h)+\mu_h(\varphi)+\imag\alpha$ can be
identified with a moment map for this action (see Section 4 in
\cite{mundet:2000}). Furthermore, $-\imag s_{h,\sigma,\chi}$ can
be identified with an element in the Lie algebra of the gauge
group $\GGG$.

We will now apply the notions of maximal weight $\lambda$ and the
function $\lambda_t$ (see Section 2.3 in \cite{mundet:2000}). Let
$A\in\AAA$ be the element giving rise to the
$\ov{\partial}$-operator which corresponds to the holomorphic
structure $E$. A simple computation tells that
(\ref{eq:grau-zero}) is equivalent to the maximal weight of
$-\imag s_{h,\sigma,\chi}$ on $(\ov{\partial}_A,\varphi)$ being
zero:
$$\lambda((\ov{\partial}_A,\varphi),-\imag s_{h,\sigma,\chi})=
\lim_{t\to\infty} \lambda_t((\ov{\partial}_A,\varphi),-\imag
s_{h,\sigma,\chi})=0.$$ Equation (\ref{eq:Hitchin-Kobayashi}) is
equivalent to the vanishing of the moment map of the action of
$\GGG$ at the pair $(\ov{\partial}_A,\phi)$. Hence we have
$\lambda_0((\ov{\partial}_A,\varphi),-\imag s_{h,\sigma,\chi})=0$,
and since $\lambda_t((\ov{\partial}_A,\varphi),-\imag s_{h,
\sigma,\chi})$ is nondecreasing as a function of $t$ it follows
that $\lambda_t((\ov{\partial}_A,\varphi),-\imag
s_{h,\sigma,\chi})=0$ for any $t$. This implies that
$e^{ts_{h,\sigma,\chi}}$ fixes the pair
$(\ov{\partial}_A,\varphi)$. That $\ov{\partial}_A$ is fixed
implies that $s_{h,\sigma,\chi}$ induces a holomorphic reduction
$\sigma_L$ of the structure group of $E$ to $L$, and that
$\varphi$ is fixed implies that $\varphi$ is contained in
$E(B)_{\sigma,\chi}^-\otimes L$.
\end{proof}

\subsection{Automorphism groups of polystable pairs}
\label{ss:simple-infinitesimally-simple}

In this section we prove that the automorphism group of an
$\alpha$-polystable pair is reductive. Let $(E,\varphi)$ be an
$L$-twisted pair. Let $\Aut(E,\varphi)$ denote the holomorphic
automorphisms of $(E,\varphi)$, i.e., the holomorphic gauge
transformations $g:E\to E$ such that $\phi\circ g^L=\phi$, where
$\phi:E^L\to B$ is the antiequivariant map corresponding to $\varphi$
and $g^L:E\times_XL\to E\times_XL$ is the transformation acting as $g$
in the $E$ factor and the identity in the $L$ factor.

The group $\Aut(E,\varphi)$ carries a natural structure of Lie
group with Lie algebra equal to $\aut(E,\phi)$.

\begin{lemma}
\label{lemma:reductive-stabilizer} Let $(E,\varphi)$ be an
$\alpha$-polystable pair. Then $\Aut(E,\varphi)$ is a reductive
Lie group.
\end{lemma}
\begin{proof}
If $(E,\varphi)$ is $\alpha$-polystable, then by Theorem
\ref{hk-twisted-pairs} there exists a reduction $h\in
\Gamma(E(\HC/H))$ of the structure group satisfying equation
(\ref{eq:Hitchin-Kobayashi}). By the arguments in the proof of
Theorem \ref{hk-twisted-pairs} this can be interpreted as the
vanishing of the moment map of the action of $\GGG$ (the gauge
group of $E_h$) on $\AAA\times\SSS$ at the point $(A,\varphi)$,
where $A$ is the Chern connection of $E$ and $h$. It follows
(see for example Proposition 1.6 in \cite{sjamaar:1995}) that
$\Aut(E,\phi)$ is the complexification of $\Aut(E,\phi)\cap\GGG$.
Any $g\in\Aut(E,\phi)\cap\GGG$ preserves simultaneously the
complex structure of $E$ and the reduction $h$, hence it also
preserves the Chern connection $A$. But the group of gauge
transformations in $\GGG$ preserving a given connection can be
identified with a closed subgroup of the automorphisms of the
fiber of $E_h$ at any given point, and consequently is a compact
Lie group. Hence $\Aut(E,\phi)\cap\GGG$ is a compact Lie group, so
by the previous argument $\Aut(E,\phi)$ is reductive.
\end{proof}

\section{$G$-Higgs  bundles  and  non-abelian Hodge  Theory}
\label{reps}

\subsection{$L$-twisted $G$-Higgs pairs, $G$-Higgs bundles and stability}
\label{sec:g-higgs-defs}

Let $G$ be a real reductive Lie group, let $H\subset G$ be a
maximal compact subgroup and let $\glie=\hlie\oplus\mlie$ be a
Cartan decomposition, so that the Lie algebra structure on $\glie$
satisfies
$$[\hlie,\hlie]\subset\hlie,\qquad
[\hlie,\mlie]\subset\mlie,\qquad [\mlie,\mlie]\subset\hlie.$$ The
group $H$ acts linearly on $\mlie$ through the adjoint
representation, and this action extends to a linear holomorphic
action of $\HC$ on $\mclie=\mlie\otimes\CC$. This is the \textbf{isotropy
representation:}
\begin{equation}
    \label{eq:isotropy-rep}
  \iota\colon \HC \to \GL(\mclie).
\end{equation}
Furthermore, the Killing form on $\glie$ induces
on $\mclie$ a Hermitian structure which is preserved by the action
of $H$.

Let $X$ be a closed Riemann surface  and let $L$ be a holomorphic line
bundle on $X$. Let $E(\mclie)= E \times_{\HC}\mclie$ be the
$\mclie$-bundle associated to $E$ via the isotropy representation.
Let $K$ be the canonical bundle of $X$.

\begin{definition}
  \label{def:g-higgs}
  An \textbf{$L$-twisted $G$-Higgs pair} on $X$ is a pair
  $(E,\varphi)$, where $E$ is a holomorphic $\HC$-principal bundle
  over $X$ and $\varphi$ is a holomorphic section of $E(\mclie)\otimes
  L$.  A \textbf{$G$-Higgs bundle} on $X$ is a $K$-twisted $G$-Higgs
  pair. Two $L$-twisted $G$-Higgs pairs $(E,\varphi)$ and
  $(E',\varphi')$ are \textbf{isomorphic} if there is an isomorphism
  $f\colon E \xra{\simeq} E'$ such that $\varphi' = (\iota(f)\otimes \Id_L)(\varphi)$.
\end{definition}

\begin{remark}
  \label{rem:G-Higgs-special-cases}
  When $G$ is compact $\liem=0$ and hence a $G$-Higgs pair is simply a
  holomorphic principal $G^\CC$-bundle.  When $G$ is complex, if
  $U\subset G$ is a maximal compact subgroup, the Cartan decomposition
  of $\lieg$ is $\lieg=\lieu +\imag \lieu$, where $\lieu$ is the Lie
  algebra of $U$. Then an $L$-twisted $G$-Higgs pair $(E,\varphi)$
  consists of a holomorphic $G$-bundle $E$ and $\varphi\in
  H^0(X,E(\lieg)\otimes L)$, where $E(\lieg)$ is the $\lieg$-bundle
  associated to $E$ via the adjoint representation.  These are the
  objects introduced originally by Hitchin \cite{hitchin:1987a} when
  $G=\SL(2,\CC)$ and $L=K$.
\end{remark}

An $L$-twisted $G$-Higgs pair is thus a particular case of the general
concept of an $L$-twisted pair introduced in Section~\ref{appendix}.
Hence $\alpha$-stability, semistability and polystability are defined
for any $\alpha\in\imag\hlie\cap\zlie$, where $\zlie$ is the center of
$\hlie^\CC$.

\subsection{Moduli
  spaces of $G$-Higgs bundles}
\label{sec:moduli-spaces}

In order to relate \emph{$G$-Higgs bundles} to
representations of the fundamental group of $X$ (or certain central
extension of the fundamental group) in $G$, one requires $\alpha$ to
lie also in the center of $\lieg$. Since we will be mostly concerned
with $G$-Higgs bundles for $G$ semisimple, this means simply
$\alpha=0$. This justifies the following terminology.

\begin{notation}
  A $G$-Higgs bundle $(E,\varphi)$ is said to be \textbf{stable} if it
  is $0$-stable. We define \textbf{semistability} and
  \textbf{polystability} of $G$-Higgs bundles similarly.
\end{notation}

Henceforth, we shall assume that $G$ is connected. Then the
topological classification of $H^\CC$-bundles $E$ on $X$ is given by a
characteristic class $c(E)\in \pi_1(H^\CC)=\pi_1(H)=\pi_1(G)$.  For a
fixed $d \in \pi_1(G)$, the \textbf{moduli space of polystable
  $G$-Higgs bundles} $\mathcal{M}_d(G)$ is by definition the set of
isomorphism classes of polystable $G$-Higgs bundles $(E,\varphi)$ such
that $c(E)=d$.  When $G$ is compact, the moduli space $\cM_d(G)$
coincides with $M_d(G^\CC)$, the moduli space of polystable
$G^\CC$-bundles with topological invariant $d$.

The moduli space $\cM_d(G)$ has the structure of a complex analytic
variety.  This can be seen by the standard slice method (see, e.g.,
Kobayashi \cite{kobayashi:1987}).  Geometric Invariant Theory
constructions are available in the literature for $G$ real compact
algebraic (Ramanathan \cite{ramanathan:1996}) and for $G$ complex
reductive algebraic (Simpson \cite{simpson:1994,simpson:1995}).  The
case of a real form of a complex reductive algebraic Lie group follows
from the general constructions of Schmitt
\cite{schmitt:2005,schmitt:2008}. We thus have the
following.
\begin{theorem}\label{alg-moduli}
The moduli space $\cM_d(G)$ is a complex analytic variety, which is
algebraic when $G$ is algebraic.
\end{theorem}

\begin{remark}
  Schmitt's construction (loc.\ cit.)  in fact applies in the more
  general setting of $L$-twisted $G$-Higgs pairs.
\end{remark}

\subsection{Deformation theory of $G$-Higgs bundles}
\label{sec:deformation-theory}

In this section we recall some standard facts about the
deformation theory of $G$-Higgs bundles.  A convenient reference
for this material is Biswas--Ramanan \cite{biswas-ramanan:1994}.

\begin{definition}\label{def:def-complex}
Let $(E,\varphi)$ be a $G$-Higgs bundle.  The \emph{deformation complex}
of $(E,\varphi)$ is the following complex of sheaves:
\begin{equation}\label{eq:def-complex}
  C^{\bullet}(E,\varphi)\colon E(\lieh^\CC) \xrightarrow{d\iota(\varphi)}
  E(\liem^\CC)\otimes K.
\end{equation}
\end{definition}
This definition makes sense because $\phi$ is a section of
$E(\lie{m}^{\CC})\otimes K$ and $[\lie{m}^{\CC},\lie{h}^{\CC}]
\subseteq \lie{m}^{\CC}$.

The following result generalizes the fact that the infinitesimal
deformation space of a holomorphic vector bundle $V$ is
isomorphic to $H^1(\End V)$.

\begin{proposition}
  \label{prop:deform}
  The space of infinitesimal deformations of a $G$-Higgs bundle
  $(E,\varphi)$ is naturally isomorphic to the hypercohomology group
  $\HH^1(C^{\bullet}(E,\varphi))$.
\end{proposition}

For any $G$-Higgs bundle there is a natural long exact sequence
\begin{equation}
  \label{eq:hyper-les}
  \begin{split}
  0 &\to \HH^0(C^{\bullet}(E,\varphi)) \to H^0 (E(\lieh^\CC))
  \xrightarrow{d\iota(\varphi)} H^0(E(\liem^\CC)\otimes K) \\
  &\to \HH^1(C^{\bullet}(E,\varphi))
  \to H^1 (E(\lieh^\CC)) \xrightarrow{d\iota(\varphi)}
  H^1(E(\liem^\CC)\otimes K) \to
  \HH^2(C^{\bullet}(E,\varphi)) \to 0.
  \end{split}
\end{equation}
As an immediate consequence we have the following result.
\begin{proposition}
\label{prop:inf-aut-sp-hypercoh}
  The infinitesimal automorphism space $\aut(E,\varphi)$ defined in
  Section~\ref{sec:inf-aut} is isomorphic to
  $\HH^0(C^{\bullet}(E,\varphi))$.
\end{proposition}

Let $d\iota\colon \lieh^{\CC} \to \End(\liem^{\CC})$ be the derivative
at the identity of the complexified isotropy representation $\iota =
\Ad_{|H^{\CC}}\colon H^{\CC} \to \Aut(\liem^{\CC})$ (cf.\
Section~\ref{sec:twisted-g-higgs-defs}). Let $\ker d\iota \subseteq
\lieh^{\CC}$ be its kernel and let $E(\ker d\iota) \subseteq
E(\lieh^{\CC})$ be the corresponding subbundle.  Then there is an
inclusion $H^0(E(\ker d\iota)) \into \HH^0(C^{\bullet}(E,\varphi))$.
\begin{definition}
  \label{def:i-simple}
  A $G$-Higgs bundle $(E,\varphi)$ is said to be
  \textbf{infinitesimally simple} if the infinitesimal automorphism
  space $\HH^0(C^{\bullet}(E,\varphi))$ is isomorphic to
$ H^0(E(\ker d\iota \cap \liez))$.
\end{definition}

Similarly, we have an inclusion $\ker\iota\cap Z(H^{\CC})\into
\Aut(E,\varphi)$, where $Z(H^\CC)$ is the center of $H^\CC$.

\begin{definition}
  \label{def:simple}
  A $G$-Higgs bundle $(E,\varphi)$ is said to be
  \textbf{simple} if $\Aut(E,\varphi)=\ker \iota \cap Z(H^\CC)$.
\end{definition}

As a consequence of Propositions \ref{prop:inf-aut-sp-hypercoh}
and \ref{proposition:stable-simple} we
 have the following.

\begin{proposition}
  \label{prop:stable-simple-G-higgs}
  Any stable $G$-Higgs bundle $(E,\varphi)$ with $\varphi\neq 0$
  is infinitesimally simple.
\end{proposition}

\begin{remark}
  \label{rem-gr-gc-ss}
  If $\ker d \iota = 0$, then $(E,\varphi)$ is infinitesimally simple
  if and only if the vanishing $\HH^0(C^{\bullet}(E,\varphi)) = 0$
  holds.  A particular case of this situation is when the group $G$
  is a complex semisimple group: indeed, in this case
  the isotropy representation is just the adjoint
  representation.
\end{remark}

Next we turn to the question of the vanishing of $\HH^2$ of the
deformation complex.  In order to deal with this question we shall
use Serre duality for hypercohomology (see e.g. Theorem 3.12 in
\cite{huybrechts:2006}), which says that there are natural
isomorphisms
\begin{equation}
  \label{eq:hyper-serre}
  \HH^i(C^{\bullet}(E,\varphi)) \cong
  \HH^{2-i}(C^{\bullet}(E,\varphi)^*\otimes K)^*,
\end{equation}
where the dual of the deformation complex (\ref{eq:def-complex}) is
\begin{displaymath}
  C^{\bullet}(E,\varphi)^*\colon E(\liem^\CC)\otimes K^{-1}
  \xrightarrow{-d\iota(\varphi)} E(\lieh^\CC).
\end{displaymath}
An important special case of this is when $G$ is a complex group.
\begin{proposition}
  \label{prop:complex-hyper-duality}
  Assume that $G$ is a complex group.  Then there is a natural
  isomorphism
  \begin{displaymath}
    \HH^2(C^{\bullet}(E,\varphi)) \cong \HH^0(C^{\bullet}(E,\varphi))^*.
  \end{displaymath}
\end{proposition}

\begin{proof}
  This is immediate from (\ref{eq:hyper-serre}) and the fact that the
  deformation complex is dual to itself, except for a sign in the map
  which does not influence the cohomology (cf.\
  Remark~\ref{rem:G-Higgs-special-cases}):
  \begin{equation}\label{eq:5}
    C^{\bullet}(E,\varphi)^*\otimes K\colon
    E(\lieg)\xrightarrow{-\ad(\varphi)}
    E(\lieg)\otimes K
  \end{equation}
  (note that when $G$ is complex $d\iota = \ad\colon\lieg\to\lieg$).
\end{proof}

\begin{remark}
  The isomorphism
  \begin{math}
    \HH^1(C^{\bullet}(E,\varphi)) \cong \HH^1(C^{\bullet}(E,\varphi))^*
  \end{math}
  is also important: it gives rise to the natural complex symplectic
  structure on the moduli space of $G$-Higgs bundles for complex
  groups $G$.
\end{remark}

We have the following key observation (cf.\ (\ref{eq:5});
again we are ignoring the irrelevant change of sign in the dual
complex).

\begin{proposition}
  \label{prop:def-g-gc}
  Let $G$ be a real group and let $G^{\CC}$ be its complexification.
  Let $(E,\varphi)$ be a $G$-Higgs bundle. Then there is an isomorphism
  of complexes:
  \begin{displaymath}
    C^{\bullet}_{G^{\CC}}(E,\varphi) \cong
    C^{\bullet}_{G}(E,\varphi) \oplus C^{\bullet}_{G}(E,\varphi)^*
    \otimes K,
  \end{displaymath}
  where $C^{\bullet}_{G^{\CC}}(E,\varphi)$ denotes the deformation
  complex of $(E,\varphi)$ viewed as a $G^{\CC}$-Higgs bundle, and
  $C^{\bullet}_{G}(E,\varphi)$ denotes the deformation complex of
  $(E,\varphi)$ viewed as a $G$-Higgs bundle.
\end{proposition}

\begin{corollary}
  \label{cor:hh-g-gc}
  With the same hypotheses as in the previous Proposition, there is an
  isomorphism
  \begin{displaymath}
    \HH^0(C^{\bullet}_{G^{\CC}}(E,\varphi)) \cong
    \HH^0(C^{\bullet}_{G}(E,\varphi)) \oplus
    \HH^2(C^{\bullet}_{G}(E,\varphi))^*.
  \end{displaymath}
\end{corollary}
\begin{proof}
  Immediate from the Proposition and Serre duality (\ref{eq:hyper-serre}).
\end{proof}

\begin{proposition}
  \label{prop:hh-vanishing}
  Let $G$ be a real semisimple group and let $G^{\CC}$ be its
  complexification.  Let $(E,\varphi)$ be a $G$-Higgs bundle which is
  stable viewed as a $G^{\CC}$-Higgs bundle.  Then the vanishing
  \begin{displaymath}
    \HH^0(C^{\bullet}_{G}(E,\varphi)) = 0 =
    \HH^2(C^{\bullet}_{G}(E,\varphi))
  \end{displaymath}
  holds.
\end{proposition}
\begin{proof}
  Since $G$ is semisimple, so is $G^{\CC}$.  Hence, in view of
  Remark~\ref{rem-gr-gc-ss}, the result follows at once from
  Corollary~\ref{cor:hh-g-gc} and
  Proposition~\ref{prop:stable-simple-G-higgs}.
\end{proof}

The following result on smoothness of the moduli space can be proved,
for example, from the standard slice method construction referred to
above.
\begin{proposition}
  \label{prop:smoothness}
  Let $(E,\varphi)$ be a stable $G$-Higgs bundle. If $(E,\varphi)$ is
  simple and
$$
\HH^2(C^{\bullet}_{G}(E,\varphi))=0,
$$ then $(E,\varphi)$
  is a smooth point in the moduli space. In particular, if
  $(E,\varphi)$ is a simple $G$-Higgs bundle which is stable as a
  $G^{\CC}$-Higgs bundle, then it is a smooth point in the moduli space.
\end{proposition}

Suppose now that we are in the situation of
Proposition~\ref{prop:smoothness} and that a local universal family
exists. Then the dimension of the
component of the moduli space containing $(E,\varphi)$ equals the
dimension of the infinitesimal deformation space
$\HH^1(C^{\bullet}_{G}(E,\varphi))$.   In view of
Proposition~\ref{prop:stable-simple-G-higgs},
Remark~\ref{rem-gr-gc-ss} and Proposition~\ref{prop:dim-moduli},
we also have
$\HH^0(C^{\bullet}_{G}(E,\varphi))=\HH^2(C^{\bullet}_{G}(E,\varphi))=0$.
So we have
$\HH^1(C^{\bullet}_{G}(E,\varphi))=-\chi(C^{\bullet}_{G}(E,\varphi))$.
A remarkable fact on this equality is that, whereas the left hand
side may depend on the choice of $(E,\phi)$, the right hand side
is independent of it, as we will see in the proposition below. We
shall refer to $-\chi(C^{\bullet}_{G}(E,\varphi))$ as the
\textbf{expected dimension} of the moduli space.

\begin{proposition}
  \label{prop:dim-moduli}
  Let $G$ be semisimple.  Then the expected dimension of the moduli space
  of $G$-Higgs bundles is
  \begin{math}
    (g-1)\dim G^{\CC}.
  \end{math}
\end{proposition}

\begin{proof}
Let $(E,\varphi)$ be any $G$-Higgs bundle.
  The long exact sequence (\ref{eq:hyper-les}) gives us
  \begin{displaymath}
    \chi(C^{\bullet}_{G}(E,\varphi)) - \chi(E(\lieh^{\CC})) +
    \chi(E(\liem^{\CC})\otimes K) = 0.
  \end{displaymath}
  Serre duality implies that $\chi(E(\liem^{\CC})\otimes K)=
  \chi(E(\liem^{\CC}))$ and
  from the Riemann--Roch formula we therefore obtain
 \begin{displaymath}
   -\chi(C^{\bullet}_{G}(E,\varphi)) = \deg(E(\liem^{\CC})) +
   (g-1)\rk(E(\liem^{\CC})) -
   \bigl(\deg(E(\lieh^{\CC}))+ (1-g)\rk(E(\lieh^{\CC})).
 \end{displaymath}
 Any invariant pairing on $\lieg^{\CC}$ (e.g. the Killing form)
 induces isomorphisms
 $E(\liem^{\CC})\simeq E(\liem^{\CC})^*$ and
 $E(\lieh^{\CC})\simeq E(\lieh^{\CC})^*$.
  Hence $\deg(E(\liem^{\CC})) = \deg(E(\lieh^{\CC})) = 0$,
  whence the result.
In particular, the value of $- \chi(C^{\bullet}_{G}(E,\varphi))$
is independent of the choice of $G$-Higgs bundle $(E,\varphi)$.
\end{proof}

\begin{remark}
  Note that the actual dimension of the moduli space (if non-empty)
  can be smaller than the expected dimension. This happens for example
  when $G=\SU(p,q)$ with $p\neq q$ and maximal Toledo invariant (this
  follows from the study of $\U(p,q)$-Higgs bundles in
  \cite{bradlow-garcia-prada-gothen:2003}) --- in this case there are
  in fact no stable $\SU(p,q)$-Higgs bundles.
\end{remark}

\subsection{$G$-Higgs bundles and Hitchin equations}

Let $G$ be a connected semisimple real Lie group. Let $(E,\varphi)$ be
a $G$-Higgs bundle over a compact Riemann surface $X$. By a slight
abuse of notation, we shall denote the $C^\infty$-objects underlying
$E$ and $\varphi$ by the same symbols. In particular, the Higgs field
can be viewed as a $(1,0)$-form: $\varphi \in
\Omega^{1,0}(E(\liem^{\CC}))$.  Given a reduction $h$ of structure
group to $H$ in the smooth $H^{\CC}$-bundle $E$, we denote by $F_h$
the curvature of the unique connection compatible with $h$ and the
holomorphic structure on $E$.  Let $\tau_h\colon
\Omega^{1,0}(E(\lieg^{\CC})) \to \Omega^{0,1}(E(\lieg^{\CC}))$ be defined
by the
compact conjugation of $\glie^\CC$ which is given fibrewise by the reduction
$h$, combined with complex conjugation on complex $1$-forms.

\begin{theorem} \label{higgs-hk}
  There exists a reduction $h$ of the structure group of $E$
  from $H^\CC$ to $H$ satisfying the Hitchin equation
  $$
  F_h -[\varphi,\tau_h(\varphi)]= 0
  $$
  if and only if $(E,\varphi)$ is polystable.
\end{theorem}

\begin{remark}
  the Hitchin equation is an equation of $2$-forms and makes sense
  without choosing metrics on $X$ and $K$, whereas the general
  Hermite--Einstein equation (\ref{eq:Hitchin-Kobayashi}) is an
  equation of scalars and requires a choice of metrics on $X$ and
  $L$. Nevertheless, for any choice of metric in the conformal class
  of the Riemann surface structure on $X$ (and hence on the
  holomorphic cotangent bundle $K$), the Hitchin equation is
  equivalent to (\ref{eq:Hitchin-Kobayashi}) for this choice of
  metric.
\end{remark}

Theorem \ref{higgs-hk} was proved by
Hitchin \cite{hitchin:1987a}  for $G=\SL(2,\CC)$ and Simpson
\cite{simpson:1988,simpson:1992}
for an arbitrary semisimple complex Lie group $G$.
The proof for an arbitrary reductive
real Lie group $G$ when $(E,\varphi)$ is stable
is given in  \cite{bradlow-garcia-prada-mundet:2003}, and the general
polystable case follows as a particular case
of the more general  Hitchin--Kobayashi correspondence given in
Theorem~\ref{hk-twisted-pairs}.

{}From the point of view of moduli spaces it is convenient
to fix  a $C^\infty$ principal   $H$-bundle
$\bE_H$ with fixed topological class $d\in \pi_1(H)$
and study  the moduli space of solutions to \textbf{Hitchin's equations}
for a pair $(A,\varphi)$ consisting of  an $H$-connection $A$ and
$\varphi\in \Omega^{1,0}(X,\bE_H(\mlie^\CC))$:
\begin{equation}\label{hitchin}
\begin{array}{l}
F_A -[\varphi,\tau(\varphi)]= 0\\
\dbar_A\varphi=0.
\end{array}
\end{equation}
Here $d_A$ is the covariant derivative associated to $A$ and $\dbar_A$
is the $(0,1)$ part of $d_A$, which defines a holomorphic structure on
$\bE_H$. Also $\tau$ is defined by the fixed reduction of structure
group $\bE_H \into \bE_H(H^{\CC})$. The gauge group $\HHH$ of $\bE_H$
acts on the space of solutions and the moduli space of solutions is
$$
\Mg_d(G):= \{ (A,\varphi)\;\;\mbox{satisfying}\;\;
(\ref{hitchin})\}/\HHH.
$$
Now,
Theorem \ref{higgs-hk} has as a consequence the following global
statement.

\begin{theorem} \label{hk}
There is a homeomorphism
$$
\cM_d(G)\cong \Mg_d(G)
$$
\end{theorem}

To explain this correspondence we interpret the moduli
space of $G$-Higgs bundles in terms of pairs $(\dbar_E, \varphi)$ consisting
of a $\dbar$-operator (holomorphic structure) on the $H^\CC$-bundle
$\bE_{H^\CC}$ obtained from
$\bE_H$ by the extension of structure group $H\subset H^\CC$, and
$\varphi\in \Omega^{1,0}(X,\bE_{H^\CC}(\mlie^\CC))$
satisfying $\dbar_E\varphi=0$.
Such pairs are in correspondence with  $G$-Higgs bundles $(E,\varphi)$,
where $E$ is the holomorphic $H^\CC$-bundle defined by the operator
$\dbar_E$ on $\bE_{H^\CC}$ and $\dbar_E\varphi=0$ is equivalent
to $\varphi\in H^0(X,E(\mlie^\CC)\otimes K)$.
The moduli space of polystable $G$-Higgs bundles $\cM_d(G)$ can now
be identified with the orbit space
$$
\{ (\dbar_E,\varphi)\;\;:\;\; \dbar_E\varphi=0,\;
\;\mbox{$(\dbar_E,\varphi)$ defines a polystable $G$-Higgs
bundle}\}/\HHH^\CC,
$$
where $\HHH^\CC$ is the gauge group of $\bE_{H^\CC}$, which is in fact
the complexification of $\HHH$.
Since  there is a one-to-one correspondence between
$H$-connections on $\bE_H$ and $\dbar$-operators on $\bE_{H^\CC}$,
the correspondence given in Theorem \ref{hk} can be interpreted
by saying that in the $\HHH^\CC$-orbit of a polystable $G$-Higgs
bundle $(\dbar_{E_0},\varphi_0)$ we can find another Higgs bundle
$(\dbar_E,\varphi)$
whose corresponding pair $(d_A,\varphi)$ satisfies
$F_A -[\varphi,\tau(\varphi)]= 0$, and this is unique up to $H$-gauge
transformations.

The infinitesimal deformation space of a solution $(A,\varphi)$ to
  Hitchin's equations can be described as the first cohomology group
  of a certain elliptic deformation complex. To do this, we follow
  Hitchin \cite[\S~5]{hitchin:1987a}. The linearized equations are:

\begin{displaymath}
    \begin{aligned}
      d_A(\dot{A})-[\dot{\varphi},\tau(\varphi)] - [\varphi,\tau(\dot{\varphi})]
      &= 0, \\
      \dbar_A\dot{\varphi}+[\dot{A}^{0,1},\varphi] &=0,
    \end{aligned}
  \end{displaymath}
  for $\dot{A} \in \Omega^1(X,\bE_H(\lieh))$ and
  $\dot{\varphi}\in \Omega^{1,0}(X,\bE_H(\liem^\CC))$.
The infinitesimal action of
$$
\psi\in  \Lie \HHH = \Omega^0(X,\bE_H(\lieh))
$$
 is
  \begin{displaymath}
    ({A},{\phi}) \mapsto (d_A\psi,[\phi,\psi]).
  \end{displaymath}
  Thus the deformation theory of Hitchin's equations is governed by
  the (elliptic) complex
  \begin{multline*}
    C^\bullet(A,\varphi)\colon \Omega^0(X,\bE_H(\lieh))
    \xra{d_0} \Omega^1(X,\bE_H(\lieh)) \oplus \Omega^{1,0}(X,\bE_H(\liem^\CC)) \\
    \xra{d_1} \Omega^2(X,\bE_H(\lieh)) \oplus \Omega^{1,1}(X,\bE_H(\liem^\CC)),
  \end{multline*}
  where the maps are
  \begin{displaymath}
    d_0(\psi) =
    (d_A\psi,[\varphi,\psi])
  \end{displaymath}
  and
  \begin{displaymath}
    d_1(\psi) = (d_A(\dot{A})-[\dot{\varphi},\tau(\varphi)] -
    [\phi,\tau(\dot{\varphi})], \dbar_A\dot{\varphi}+[\dot{A}^{0,1},\varphi]).
  \end{displaymath}
  The fact that $(A,\varphi)$ is a solution
  to the equations, together with the gauge invariance of the equations,
  guarantees that $d_1 \circ d_0 = 0$.
Denote by $H^i(C^\bullet(A,\varphi))$ the cohomology groups of the gauge
  theory deformation complex $C^\bullet(A,\varphi)$.

Let
$$
\Aut(A,\varphi):=\{h\in\HHH \;\;:\;\; h^*A=A,\;\;\mbox{and}\;\; \iota(h)(\varphi)=\varphi\}.
$$

Here $\iota: H\to \Aut(\liem)$ is the isotropy representation.
Clearly $Z(H)\cap \ker \iota\subset \Aut(A,\varphi)$.

\begin{definition}
Let $(A,\varphi)$ be a solution of (\ref{hitchin}). We say that
$(A,\varphi)$ is \textbf{irreducible} if and only if
$\Aut(A,\varphi)=Z(H)\cap \ker \iota$.
We say that $(A,\varphi)$ is  \textbf{infinitesimally irreducible}
if the Lie algebra of $\Aut(A,\varphi)$, which is identified with
$H^0(C^\bullet(A,\varphi))$ equals  $Z(\lieh)\cap \ker d\iota$.
\end{definition}

  \begin{proposition}
    Assume that $H^0(C^\bullet(A,\varphi)) = H^2(C^\bullet(A,\varphi)) = 0$ and that
    $(A,\varphi)$ is irreducible. Then $\Mg_d$ is smooth
    at $[A,\varphi]$ and the tangent space is
    \begin{displaymath}
      T_{[A,\varphi]}\Mg_d \cong H^1(C^\bullet(A,\varphi)).
    \end{displaymath}
  \end{proposition}

  For a proper understanding of many aspects of the geometry of the moduli space of
  Higgs bundles, one needs to consider the moduli space as the gauge
  theory moduli space $\Mg_d(G)$. On the other hand, the formulation of the
  deformation theory in terms of hypercohomology is very
  convenient. Fortunately, one has the following.
  \begin{proposition}
    At a smooth point of the moduli space, there is a natural
    isomorphism of infinitesimal deformation spaces
    \begin{displaymath}
      H^1(C^\bullet(A,\varphi)) \cong \HH^1(C^\bullet(E,\varphi)),
    \end{displaymath}
    where the holomorphic structure on the Higgs bundle $(E,\varphi)$ is
    given by $\dbar_A$.
  \end{proposition}

  As in Donaldson--Kronheimer \cite[\S~6.4]{donaldson-kronheimer:1990}
  this can be seen by using a Dolbeault resolution to calculate
  $\HH^1(C^\bullet(E,\varphi))$ and using harmonic representatives of
  cohomology classes, via Hodge theory.  For this reason we can freely
  apply the complex deformation theory described in
  Section~\ref{sec:deformation-theory} to the gauge theory situation.

  The following result is not essential for the present paper but we
  include it here for completeness. It can be deduced from the
  treatment of the Hitchin--Kobayashi correspondence given in
  Section~\ref{appendix}.
\begin{proposition}
Under the correspondence given by Theorem \ref{hk},
a stable
$G$-Higgs bundle corresponds to
an infinitesimally irreducible
solution to Hitchin equations,
while a  $G$-Higgs bundle which  is
stable and simple
is in correspondence with
 an irreducible solution.
\end{proposition}

\subsection{Surface group representations}

Let $X$ be a closed oriented surface of genus $g$ and let
\begin{displaymath}
  \pi_{1}(X) = \langle a_{1},b_{1}, \dotsc, a_{g},b_{g} \suchthat
  \prod_{i=1}^{g}[a_{i},b_{i}] = 1 \rangle
\end{displaymath}
be its fundamental group.  Let $G$ be a connected reductive real Lie group.
By a \textbf{representation} of $\pi_1(X)$ in
$G$ we understand a homomorphism $\rho\colon \pi_1(X) \to G$.
The set of all such homomorphisms,
$\Hom(\pi_1(X),G)$, can be naturally identified with the subset
of $G^{2g}$ consisting of $2g$-tuples
$(A_{1},B_{1}\dotsc,A_{g},B_{g})$ satisfying the algebraic equation
$\prod_{i=1}^{g}[A_{i},B_{i}] = 1$.  This shows that
$\Hom(\pi_1(X),G)$ is a real analytic  variety, which is algebraic
if $G$ is algebraic.

The group $G$ acts on $\Hom(\pi_1(X),G)$ by conjugation:
\[
(g \cdot \rho)(\gamma) = g \rho(\gamma) g^{-1}
\]
for $g \in G$, $\rho \in \Hom(\pi_1(X),G)$ and $\gamma\in
\pi_1(X)$. If we restrict the action to the subspace
$\Hom^+(\pi_1(X),g)$ consisting of \emph{reductive
representations}, the orbit space is Hausdorff (see Theorem 11.4
in \cite{Ri}).  By a \textbf{reductive representation} we mean one
that composed with the adjoint representation in the Lie algebra
of $G$ decomposes as a sum of irreducible representations. If $G$
is algebraic this is equivalent to the Zariski closure of the
image of $\pi_1(X)$ in $G$ being a reductive group. (When $G$ is
compact every representation is reductive.) Define the
\emph{moduli space of representations} of $\pi_1(X)$ in $G$ to be
the orbit space
\[
\mathcal{R}(G) = \Hom^{+}(\pi_1(X),G) / G. \]

One has the following (see e.g.\ Goldman \cite{goldman:1984}).

\begin{theorem}
The moduli space $\cR(G)$ has the structure of a real analytic variety, which
is algebraic if $G$ is algebraic and is a complex variety if $G$ is complex.
\end{theorem}

Given a representation $\rho\colon\pi_{1}(X) \to
G$, there is an associated flat $G$-bundle on
$X$, defined as
\begin{math}
  E_{\rho} = \widetilde{X}\times_{\rho}G
\end{math},
where $\widetilde{X} \to X$ is the universal cover and
$\pi_{1}(X)$ acts on $G$ via $\rho$. This gives in fact  an
identification between the set of equivalence classes of
representations  $\Hom(\pi_1(X),G) / G$ and the set of equivalence
classes of flat $G$-bundles, which in turn is parameterized by the
cohomology set $H^1(X,G)$. We can then assign  a topological
invariant  to a representation $\rho$ given by the  characteristic
class $c(\rho):=c(E_{\rho})\in \pi_1(G)$ corresponding to
$E_{\rho}$. To define this, let $\widetilde G$ be the universal
covering group of $G$. We have an exact  sequence
$$
1 \lto\pi_1(G)\lto \widetilde G \lto G \lto 1
$$
which gives rise to the (pointed sets) cohomology sequence
\begin{equation}\label{characteristic}
H^1(X, {\widetilde G}) \lto  H^1(X, {G})\stackrel{c}{\lto}   H^2(X, \pi_1(G)).
\end{equation}

Since $\pi_1(G)$ is abelian the orientation of $X$ defines an
isomorphism
$$
 H^2(X, \pi_1(G))\cong  \pi_1(G),
$$
and $c(E_\rho)$ is defined as the image of $E$ under the last map in
(\ref{characteristic}). Thus the class $c(E_\rho)$ measures  the
obstruction to lifting $E_\rho$ to a flat $\widetilde G$-bundle, and hence to
lifting $\rho$ to a representation of $\pi_1(X)$ in $\widetilde G$.
For a  fixed
 $d\in \pi_1(G)$, the
{\em moduli space of reductive representations} $\mathcal{R}_d(G)$
with topological invariant $d$
is defined as  the subvariety
\begin{equation}\label{eq:RdG}
  \mathcal{R}_d(G):=\{[\rho] \in \mathcal{R}(G) \suchthat
  c(\rho)=d\},
\end{equation}
where as usual $[\rho]$ denotes the $G$-orbit $G\cdot \rho$ of
$\rho\in\Hom^+(\pi_1(X),G)$.

One can study  deformations of a class of representations
$[\rho]\in \mathcal{R}_d(G)$ by means of group cohomology (see
\cite{goldman:1984}). The Lie algebra $\lieg$ is endowed  with the
structure of a $\pi_1(X)$-module by means of the composition
$$
\pi_1(X) \stackrel{\rho}{\lto}  G \stackrel{\Ad}{\lto}  \Aut(\lieg).
$$

\begin{definition}
Let  $\rho:\pi_1(X)\to G$ be a representation of $\pi_1(X)$ in
$G$. Let $Z_G(\rho)$ be the centralizer in $G$ of
$\rho(\pi_1(X))$. We say that  $\rho$ is {\bf irreducible} if and
only if it is reductive and $Z_G(\rho)=Z(G)$, where $Z(G)$ is the
center of $G$. We say that $\rho$ is an  {\bf infinitesimally
irreducible} representation if it is reductive and   $\Lie
Z_G(\rho)=\Lie Z(G)$.
\end{definition}

One has the following basic facts (\cite{goldman:1984}).

\begin{proposition}
  \begin{enumerate}
  \item The Zariski tangent space to $\mathcal{R}_d(G)$ at an
    equivalence class $[\rho]$ is isomorphic to the cohomology group
    $H^1(\pi_1(X),\lieg_{\Ad\circ \rho})$.

  \item $H^0(\pi_1(X),\lieg_{\Ad\circ \rho})\cong \Lie Z_G(\rho)$.

  \item $H^2(\pi_1(X),\lieg_{\Ad\circ \rho})\cong
    H^0(\pi_1(X),\lieg_{\Ad\circ \rho})^*$
  \end{enumerate}
\end{proposition}

From this one obtains the following (\cite{goldman:1984}).

\begin{proposition}
Let $G$ be a semisimple Lie group and let $\rho:\pi_1(X)\to G$ be
irreducible. Then the equivalence class $[\rho]$ is a smooth point
in $\mathcal{R}_d(G)$.
\end{proposition}

This is simply because $Z_G(\rho)=Z(G)$ is finite and hence
$$
H^0(\pi_1(X),\lieg_{\Ad\circ \rho})=H^2(\pi_1(X),\lieg_{\Ad\circ \rho})=0.
$$

An alternative way to study deformations of a representation is by
using  the corresponding flat connection. To explain this, let
$\bE$ be a $C^\infty$ principal $G$-bundle over $X$ with fixed
topological class $d\in\pi_1(G)=\pi_1(H)$. Let $D$ be a
$G$-connection on $\bE$ and let $F_D$ be its curvature. If $D$ is
flat, i.e.\ $F_D=0$, then the holonomy of $D$ around a closed loop
in $X$ only depends on the homotopy class of the loop and thus
defines a representation of $\pi_1(X)$ in $G$.  This gives an
identification\footnote{even when $G$ is complex algebraic, this
is
  merely a real  \emph{analytic} isomorphism, see Simpson
  \cite{simpson:1992,simpson:1994,simpson:1995}},
$$
\cR_d(G)\cong
\{\mbox{Reductive $G$-connections}\;\; D \suchthat
F_D=0\}/\GGG,
$$
where, by definition, a flat connection is reductive if the
corresponding representation of $\pi_1(X)$ in $G$ is reductive,
and $\GGG$ is the group of automorphisms of $\bE$ --- the
\textbf{gauge
  group}. We can now linearize the flatness condition near a flat
  connection $D$:
  \begin{displaymath}
    \ddt F(D + bt)_{t=0} = D(b)
  \end{displaymath}
  for $b \in \Omega^1(X, \bE(\lieg))$.

  Linearize the action of the gauge group $D\mapsto g\cdot D =
  gDg^{-1}$.
For $g(t) = \exp(\psi t)$ with $\psi \in \Omega^0(X,
  \bE(\lieg))$,
  \begin{displaymath}
    \ddt (g(t)\cdot D)_{t=0} = D(\psi).
  \end{displaymath}
  Thus the infinitesimal deformation space is $H^1$ of the complex
  \begin{displaymath}
    0 \to \Omega^0(X, \bE(\lieg)) \xra{D} \Omega^1(X, \bE(\lieg))
    \xra{D} \Omega^2(X,\bE(\lieg)) \to 0.
  \end{displaymath}
  Note that $F_D=D^2=0$ means that this is in fact a
  complex.

\subsection{Representations and $G$-Higgs bundles}

We assume now that $G$ is connected  and  semisimple.
With the notation of the  previous sections, we have the following
non-abelian Hodge Theorem for representations of the fundamental group
of a closed Riemann surface in a semisimple connected Lie group.

\begin{theorem}\label{na-Hodge}
Let $G$ be a connected semisimple real Lie group. There is a homeomorphism
$\mathcal{R}_d(G) \cong \mathcal{M}_d(G)$. Under this homeomorphism,
stable $G$-Higgs bundles correspond to  infinitesimally
irreducible representations, and stable and simple $G$-Higgs bundles
correspond to irreducible representations.
\end{theorem}

\begin{remark}
On the open subvarieties defined by the smooth
points of $\mathcal{R}_d$ and $\mathcal{M}_d$, this correspondence is in
fact an isomorphism of real analytic varieties.
\end{remark}

\begin{remark} There is a similar correspondence when $G$ is reductive
but not semisimple. In this case, it makes sense to consider
nonzero values of the stability parameter $\alpha$. The resulting
Higgs bundles can be geometrically interpreted in terms of
representations of the universal central extension of the
fundamental group of $X$, and the value of $\alpha$ prescribes the
image of a generator of the center in the representation.
\end{remark}

The proof of Theorem \ref{na-Hodge} is the combination of two
existence theorems for gauge-theoretic equations. To explain this, let
$\bE_G$ be, as above, a $C^\infty$ principal $G$-bundle over $X$ with
fixed topological class $d\in\pi_1(G)=\pi_1(H)$. Every $G$-connection
$D$ on $\bE_G$ decomposes uniquely as
$$
D=d_A + \psi,
$$
where $d_A$ is an $H$-connection on $\bE_H$ and
$\psi\in \Omega^1(X,\bE_H(\mlie))$.  Let
$F_A$ be the curvature of $d_A$.
We consider the following set of equations for the pair $(d_A,\psi)$:
\begin{equation}\label{harmonicity}
\begin{array}{l}
F_A +\frac{1}{2}[\psi,\psi]= 0\\
d_A\psi=0  \\
d_A^\ast\psi=0.
\end{array}
\end{equation}
These equations are invariant under the action of $\HHH$, the gauge group of
$\bE_H$. A theorem of Corlette \cite{corlette:1988}, and
Donaldson \cite{donaldson:1987} for $G=\SL(2,\CC)$,  says the following.
\begin{theorem}\label{corlette} There is a homeomorphism
$$
\{\mbox{Reductive $G$-connections}\;\; D \suchthat
F_D=0\}/\GGG\cong \{(d_A,\psi)\;\;\mbox{satisfying}\;\;
(\ref{harmonicity})\}/\HHH.
$$
\end{theorem}

The first two equations in (\ref{harmonicity}) are equivalent
to the flatness of $D=d_A+\psi$, and Theorem \ref{corlette}
simply says that in the $\GGG$-orbit of a reductive flat $G$-connection
$D_0$ we can find a flat $G$-connection $D=g(D_0)$ such that if we
write $D=d_A+\psi$, the
additional condition $d_A^\ast\psi=0$ is satisfied. This can be interpreted
more geometrically in terms of the reduction  $h=g(h_0)$ of $\bE_G$
to an $H$-bundle obtained by the action of $g\in\GGG$ on $h_0$.
The equation
$d_A^\ast\psi=0$ is equivalent to the harmonicity of the
$\pi_1(X)$-equivariant map $\widetilde X \to G/H$ corresponding to
the new reduction of structure group $h$.

To complete the argument, leading to Theorem \ref{na-Hodge}, we just need
Theorem \ref{higgs-hk} and the following simple result.
\begin{proposition}\label{prop:circle}
The correspondence $(d_A,\varphi)\mapsto (d_A,\psi:=\varphi-\tau(\varphi))$
defines  a homeomorphism
$$
\{(d_A,\varphi)\;\;\mbox{satisfying}\;\;
(\ref{hitchin})\}/\HHH\cong
\{(d_A,\psi)\;\;\mbox{satisfying}\;\;
(\ref{harmonicity})\}/\HHH.
$$
\end{proposition}

\subsection{Jordan--H\"older reduction of  $G$-Higgs bundles}
\label{sec:twisted-g-higgs-defs}

The purpose of this
subsection is to show that the Jordan--H\"older reduction (defined
in Section~\ref{ss:Jordan-Holder}) of a $G$-Higgs bundle is itself a
$G'$-Higgs bundle for a reductive subgroup $G'\subset G$.

\begin{proposition}
Let  $(E,\varphi)$ be  an $L$-twisted $G$-Higgs pair which is
$\alpha$-polystable but not $\alpha$-stable. Then  the
Jordan--H\"older reduction of $(E,\varphi)$ is an $L$-twisted
$G'$-Higgs pair for some reductive subgroup $G'\subset G$.
\end{proposition}

\begin{proof}
Recall from Section  \ref{ss:Jordan-Holder}  that in the
Jordan--H\"older reduction $(E',\varphi',H',(\mclie)')$ of
$(E,\varphi,H,\mclie)$ the subgroup $H'\subset H$ is defined as
the centralizer of a torus $T\subset H$ and that $(\mclie)'$ is
the fixed point set of $T$ acting on $\mclie$. So it suffices to
prove that the Lie algebra structure on $\hlie\oplus\mlie$ induces
a structure of Cartan pair on $(\hlie',(\mclie)'\cap\mlie)$. The
action of $T$ on $\hlie$ and $\mlie$ induces decompositions
$$\hlie=\bigoplus_{\eta\in T^{\vee}}\hlie_{\eta}
\qquad\text{and}\qquad \mlie=\bigoplus_{\eta\in
T^{\vee}}\mlie_{\eta},$$ where $T^{\vee}$ denotes the group of
characters of $T$ (for which we use additive notation). Then one
has, as usual,
$$[\hlie_{\eta},\hlie_{\mu}]\subset\hlie_{\eta+\mu},\qquad
[\hlie_{\eta},\mlie_{\mu}]\subset\mlie_{\eta+\mu},\qquad
[\mlie_{\eta},\mlie_{\mu}]\subset\hlie_{\eta+\mu}$$ for any pair
of characters $\eta,\mu\in T^{\vee}$. Taking $\eta=\mu=0$ and
observing that $\hlie'=\hlie_0$ and $(\mclie)'\cap\mlie=\mlie_0$,
it follows that
$$[\hlie',\hlie']\subset\hlie',\qquad
[\hlie',(\mclie)'\cap\mlie]\subset(\mclie)'\cap\mlie,\qquad
[(\mclie)'\cap\mlie,(\mclie)'\cap\mlie]\subset\hlie',$$ so that
$(\hlie',(\mclie)'\cap\mlie)$ is certainly a Cartan pair.

\end{proof}

\begin{remark}
  We can make a more precise statement: defining $G'$ as the
centralizer of $T$ inside $G$ we have proved that the
Jordan--H\"older reduction of $(E,\varphi)$ is an $L$-twisted
$G'$-Higgs pair.
\end{remark}

\section{Simplified stability of $G$-Higgs bundles}
\label{twisted-higgs}

In this section we give concrete examples of $G$-Higgs bundles for
various interesting cases of real reductive groups $G$ and we show how
the general stability conditions can be simplified to more workable
conditions in these particular cases. Even though our main interest
lies in $G$-Higgs bundles, we state and prove our results in the more
general setting of $L$-twisted $G$-Higgs pairs, since this requires no
extra work.

\subsection{$\Sp(2n,\CC)$-Higgs bundles}

\label{ss:L-twisted-sp(2n,C)}

Consider now the case $G=\Sp(2n,\CC)$. A maximal compact subgroup
of $G$ is $H=\Sp(2n)$ and hence  $H^\CC$ coincides with
$\Sp(2n,\CC)$. Now, if $\WW=\CC^{2n}$ is the fundamental
representation of $\Sp(2n,\CC)$ and $\omega$ denotes the standard
symplectic form on $\WW$, the isotropy representation space is
$$
\lie{m}^\CC=\splie(\WW)=\splie(\WW,\omega):= \{\xi\in\End(\WW)\mid
\omega(\xi\cdot,\cdot)+\omega(\cdot, \xi\cdot)=0\}\subset\End\WW,
$$
so it coincides with the adjoint representation of $\Sp(2n,\CC)$
on its Lie algebra. An $L$-twisted $\Sp(2n,\CC)$-Higgs pair is
thus a pair consisting of a rank $2n$ holomorphic symplectic
vector bundle $(W,\Omega)$ over $X$ (so $\Omega$ is a holomorphic
section of $\Lambda^2W^*$ whose restriction to each fiber of $W$
is non degenerate) and a section
$$
\Phi \in H^0(L\otimes \splie(W)),
$$
where $\splie(W)$ is the vector bundle whose fiber over $x$ is
given by $\splie(W_x,\Omega_x)$.

\begin{remark}
Since the center of $\splie(2n,\CC)$ is trivial, $\alpha=0$ is the only
possible value for which stability of an $L$-twisted
$\Sp(2n,\CC)$-Higgs pair  is defined.
\end{remark}

Define for any filtration by holomorphic subbundles
$$\WWW=(0=W_0\subsetneq W_1\subsetneq W_2\subsetneq\dots\subsetneq
W_k=W)$$ satisfying $W_{k-i}=W_i^{\perp_{\Omega}}$ for any $i$
(here $\perp_{\Omega}$ denotes the perpendicular with respect to
$\Omega$) the set
$$\Lambda(\WWW)=\{(\lambda_1,\lambda_2,\dots,\lambda_k)\in\RR^k
\mid \lambda_i\leq\lambda_{i+1}\text{ and }
\lambda_{k-i+1}+\lambda_i=0 \text{ for any $i$ }\}.$$ For any
$\lambda\in\Lambda(\WWW)$ define the following subbundle of
$L\otimes \End W$:
$$N(\WWW,\lambda)=
\bigcap_{\lambda_i\geq\lambda_j}
\{\Phi\in L\otimes\splie(W)\suchthat \Phi(W_i)\subset L\otimes W_j\}.$$ 
Define also
$$d(\WWW,\lambda)=\sum_{j=1}^{k-1}(\lambda_j-\lambda_{j+1})\deg W_j$$
(note that since $W$ carries a symplectic structure we have
$W\simeq W^*$ and hence $\deg W=\deg W_k=0$).

Following again Sections~\ref{ss:app-pol,sem,stability} and
Section~\ref{ss:app-filtrations}, the pair
$((W,\Omega),\Phi)$ is said to be

\begin{itemize}

\item {\bf semistable} if for any filtration $\WWW$ as above and
any $\lambda\in\Lambda(\WWW)$ such that $\Phi\in
H^0(N(\WWW,\lambda))$, the following inequality holds:
$d(\WWW,\lambda)\geq 0.$

\item \textbf{stable} if it is semistable and furthermore, for any
choice of filtration $\WWW$ and $\lambda\in\Lambda(\WWW)$ which is
not identically zero (so for which there is a $j < k$ such that
$\lambda_j < \lambda_{j+1}$), and such that $\Phi \in
H^0(N(\WWW,\lambda))$, we have $d(\WWW,\lambda) > 0$.

\item \textbf{polystable} if it is
 semistable and for any filtration $\WWW$ as above and
 $\lambda\in\Lambda(\WWW)$ satisfying $\lambda_i<\lambda_{i+1}$
 for each $i$, $\psi\in H^0(N(\WWW,\lambda))$ and
 $d(\WWW,\lambda)=0$, there is an isomorphism
 $$W\simeq W_1\oplus W_2/W_1\oplus\dots\oplus W_k/W_{k-1}$$
 such that the pairing via $\Omega$ any element of the summand
 $W_i/W_{i-1}$ with
 an element of the summand $W_j/W_{j-1}$ is zero unless $i+j=k+1$;
 furthermore, via the isomorphism above,
 $$\Phi\in H^0(\bigoplus_i
 L\otimes\Hom(W_i/W_{i-1},W_i/W_{i-1})).$$
\end{itemize}

We now prove a result which shows that the definition of
(semi-,poly-)stability which we have given coincides with the usual one
in the literature. Recall that if $(W,\Omega)$ is a symplectic
vector bundle, a subbundle $W'\subset W$ is said to be isotropic
if $W' \subset W'^{\perp_\Omega}$ and coisotropic if $W' \supset
W'^{\perp_\Omega}$.

\begin{theorem}
\label{thm:sp(2n,C)-stability} An $L$-twisted $\Sp(2n,\CC)$-Higgs
pair $((W,\Omega),\Phi)$ is semistable if and only if for any
isotropic subbundle $W'\subset W$ such that $\Phi(W')\subset
L\otimes W'$ we have $\deg W'\leq 0$. Furthermore,
$((W,\Omega),\Phi)$ is stable if for any nonzero and strict
isotropic subbundle $0\neq W'\subset W$ such that $\Phi(W')\subset
L\otimes W'$ we have $\deg W'<0$. Finally, $((W,\Omega),\Phi)$ is
polystable if it is semistable and for any nonzero and strict
isotropic (resp., coisotropic) subbundle $W'\subset W$ such that $\Phi(W')\subset
L\otimes W'$ and $\deg W'=0$ there is a coisotropic (resp.,
isotropic) subbundle
$W''\subset W$ such that $\Phi(W'')\subset L\otimes W''$ and
$W=W'\oplus W''$.
\end{theorem}

\begin{proof}
  The proof follows the same ideas as the (more complicated) proof of
  Theorem~\ref{thm:simplified-conditions-sp-2n-R} below, so we just give a
  sketch. Take an $L$-twisted $\Sp(2n,\CC)$-Higgs pair
  $((W,\Omega),\Phi)$, and assume that for any isotropic subbundle
  $W'\subset W$ such that $\Phi(W')\subset L\otimes W'$ we have $\deg
  W'\leq 0$. We want to prove that $((W,\Omega),\Phi)$ is
  semistable. Choose any filtration $\WWW=(0\subsetneq W_1\subsetneq
  W_2\subsetneq\dots\subsetneq W_k=W)$ satisfying
  $W_{k-i}=W_i^{\perp_{\Omega}}$ for any $i$. We have to understand
  the geometry of the convex set $$\Lambda(\WWW,\Phi)=
  \{\lambda\in\Lambda(\WWW)\mid \Phi\in
  N(\WWW,\lambda)\}\subset\RR^k.$$ Define for that $\JJJ=\{j\mid
  \Phi(W_j)\subset L\otimes W_j\}= \{j_1,\dots,j_r\}$. One checks
  easily that if $\lambda=(\lambda_1,\dots,\lambda_k)\in\Lambda(\WWW)$
  then
\begin{equation}
\label{eq:Lambda-JJJ} \lambda\in\Lambda(\WWW,\Phi)
\Longleftrightarrow \lambda_a=\lambda_b\text{ for any $j_i\leq
a\leq b\leq j_{i+1}$}.
\end{equation}
We claim that the set of indices $\JJJ$ is symmetric:
\begin{equation}
\label{eq:JJJ-symmetric} j\in\JJJ \Longleftrightarrow k-j\in\JJJ.
\end{equation}
To check this it suffices to prove that $\Phi(W_j)\subset L\otimes
W_j$ implies that $\Phi(W_j^{\perp_{\Omega}})\subset L\otimes
W_j^{\perp_{\Omega}}$. Suppose that this is not true, so that for
some $j$ we have $\Phi W_j\subset L\otimes W_j$ and there exists
some $w\in W_j^{\perp_{\Omega}}$ such that $\Phi w \notin L\otimes
W_j^{\perp_{\Omega}}$. Then there exists $v\in W_j$ such that
$\Omega(v,\Phi w)\neq 0$. However, since $\Phi\in
H^0(L\otimes\splie(W))$, we must have $\Omega(v,\Phi
w)=-\Omega(\Phi v,w)$, and the latter vanishes because  by
assumption $\Phi v$ belongs to $W_j$. So we have reached a
contradiction.

Let $\JJJ'=\{j\in\JJJ\mid 2j\leq k\}$ and define for any
$j\in\JJJ'$ the vector
$$L_j=-\sum_{c\leq j}e_c+\sum_{d\geq k-j+1}e_d,$$
where $e_1,\dots,e_k$ is the canonical basis of $\RR^k$. It
follows from (\ref{eq:Lambda-JJJ}) and (\ref{eq:JJJ-symmetric})
that $\Lambda(\WWW,\Phi)$ is the positive span of the vectors $\{
L_j\mid j\in\JJJ'\}$. Consequently, we have
$$d(\WWW,\lambda)\geq 0\text{ for any $\lambda\in\Lambda(\WWW,\Phi)$ }
\Longleftrightarrow d(\WWW,L_j)\geq 0\text{ for any $j$ }.$$ One
computes $d(\WWW,L_j)=-\deg W_{k-j}-\deg W_j$. On the other hand,
since we have an exact sequence $0\to W_{k-j}\to W^*\to W_j^*\to
0$ (the injective arrow is given by the pairing with $\Omega$) we
have $0=\deg W^*=\deg W_{k-j}+\deg W_j^*$, so $\deg W_{k-j}=\deg
W_j$ and consequently $d(\WWW,L_j)\geq 0$ is equivalent to $\deg
W_j\leq 0$, which holds by assumption. Hence $((W,\Omega),\Phi)$
is semistable.

The converse statement, namely, that if $((W,\Omega),\Phi)$ is
semistable then for any isotropic subbundle $W'\subset W$ such
that $\Phi(W')\subset L\otimes W'$ we have $\deg W'\leq 0$ is
immediate by applying the stability condition of the filtration
$0\subset W'\subset W'^{\perp_{\Omega}}\subset W$.

Finally, the proof of the second statement on stability is very
similar to case of semistability, so we omit it. The statement on
polystability is also straightforward.
\end{proof}

\subsection{$\SL(n,\CC)$-Higgs bundles}

\label{ss:L-twisted-sl(2n,C)}

If $G=\SL(n,\CC)$ then the maximal compact subgroup of $G$ is
$H=\SU(n)$ and hence  $H^\CC$ coincides with $\SL(n,\CC)$. Now, if
$\WW=\CC^{n}$ is the fundamental representation of $\SL(n,\CC)$,
the isotropy representation space is given by the traceless
endomorphisms of $\WW$
$$
\lie{m}^\CC=\sllie(\WW)=\{\xi\in\End(\WW)\mid \Tr\xi
=0\}\subset\End\WW,
$$
so it coincides again with the adjoint representation of
$\SL(n,\CC)$ on its Lie algebra. An $L$-twisted $\SL(n,\CC)$-Higgs
pair is thus a pair consisting of a rank $n$ holomorphic vector
bundle $W$ over $X$ endowed with a trivialization $\det
W\simeq\OOO$ and a holomorphic section
$$
\Phi \in H^0(L\otimes \End_0W),
$$
where $\End_0W$ denotes the bundle of traceless endomorphisms of
$W$.

\begin{remark}
Since the center of $\sllie(n,\CC)$ is trivial, $\alpha=0$ is the only
possible value for which stability of an $L$-twisted
$\SL(n,\CC)$-Higgs pair  is defined.
\end{remark}

Define for any filtration by holomorphic subbundles
$$\WWW=(0=W_0\subsetneq W_1\subsetneq W_2\subsetneq\dots\subsetneq
W_k=W)$$ the convex set
$$\Lambda(\WWW)=\{(\lambda_1,\lambda_2,\dots,\lambda_k)\in\RR^k
\mid \lambda_i\leq\lambda_{i+1} \text{ for any $i$ and } \sum_i
\rk W_i(\lambda_i-\lambda_{i+1})=0\}.$$ For any
$\lambda\in\Lambda(\WWW)$ define the following subbundle of
$L\otimes \End W$:
$$N(\WWW,\lambda)=
\bigcap_{\lambda_i\geq\lambda_j}
\{\Phi\in L\otimes\End_0(W)\suchthat \Phi(W_i)\subset L\otimes W_j\}.$$ 
Define also
$$d(\WWW,\lambda)=\sum_{j=1}^{k-1}(\lambda_j-\lambda_{j+1})\deg W_j$$
(since $\det W$ is trivial we have $\deg W=\deg W_k=0$).

Following again Sections~\ref{ss:app-pol,sem,stability} and
\ref{ss:app-filtrations}, $(W,\Phi)$ is said to be:

\begin{itemize}

\item {\bf semistable} if for any filtration $\WWW$ and
$\lambda\in\Lambda(\WWW)$ such that $\Phi\in
H^0(N(\WWW,\lambda))$, we have $d(\WWW,\lambda)\geq 0.$

\item \textbf{stable} if it is semistable and furthermore, for any
choice of filtration $\WWW$ and $\lambda\in\Lambda(\WWW)$ which is
not identically zero (so for which there is a $j < k$ such that
$\lambda_j < \lambda_{j+1}$), and such that $\Phi \in
H^0(N(\WWW,\lambda))$, we have $d(\WWW,\lambda) > 0$.

 \item \textbf{polystable} if it is
 semistable and for any filtration $\WWW$ as above and
 $\lambda\in\Lambda(\WWW)$ satisfying $\lambda_i<\lambda_{i+1}$
 for each $i$, $\psi\in H^0(N(\WWW,\lambda))$ and
 $d(\WWW,\lambda)=0$, there is an isomorphism
 $$W\simeq W_1\oplus W_2/W_1\oplus\dots\oplus W_k/W_{k-1}$$
with respect to which
$$\Phi\in H^0(\bigoplus_i L\otimes\Hom(W_i/W_{i-1},W_i/W_{i-1})).$$
\end{itemize}

Again we have a result as Theorem \ref{thm:sp(2n,C)-stability}
implying that the present notions of (semi)stability coincide with
the usual ones.

\begin{theorem}
\label{thm:sl(n,C)-stability} An $L$-twisted $\SL(n,\CC)$-Higgs
pair $(W,\Phi)$ is semistable if and only if for any subbundle
$W'\subset W$ such that $\Phi(W')\subset L\otimes W'$ we have
$\deg W'\leq 0$. Furthermore, $(W,\Phi)$ is stable if for any
nonzero and strict subbundle $W'\subset W$ such that
$\Phi(W')\subset L\otimes W'$ we have $\deg W'<0$. Finally,
$(W,\Phi)$ is polystable if it is semistable and for each
subbundle $W'\subset W$ such that $\Phi(W')\subset L\otimes W'$
and $\deg W'=0$ there is another subbundle $W''\subset W$
satisfying $\Phi(W'')\subset L\otimes W''$ and $W=W'\oplus W''$.
\end{theorem}

The proof of Theorem \ref{thm:sl(n,C)-stability} is very similar
to that of Theorem \ref{thm:sp(2n,C)-stability}, so we omit it.

\subsection{$\Sp(2n,\RR)$-Higgs bundles}
\label{sec:L-twist}

\label{ss:L-twisted-sp(2n,R)}

Let $G=\Sp(2n,\RR)$. The maximal compact subgroup of $G$ is
$H=\U(n)$ and hence  $H^\CC=\GL(n,\CC)$. Now, if  $\VV=\CC^n$ is
the fundamental representation of $\GL(n,\CC)$, then the isotropy
representation space is:
$$
\lie{m}^\CC=S^2\VV\oplus  S^2\VV^*.
$$
An $L$-twisted $\Sp(2n,\RR)$-Higgs pair is thus a pair consisting
of a rank $n$ holomorphic vector bundle $V$ over $X$ and a section
$$
\varphi = (\beta,\gamma) \in H^0(L\otimes S^2V\oplus L \otimes
S^2V^*).
$$
In the particular case when $L = K$, we obtain the notion of an
$\Sp(2n,\RR)$-Higgs bundle.

\begin{notation}
  \label{not:symmetrized-tensor}
  If $W$ is a vector bundle and $W',W''\subset W$ are subbundles, then
  $W'\otimes_S W''$ denotes the image of $W'\otimes W''\subset W\otimes W$
  under the symmetrization map $W\otimes W\to S^2W$ (this
  should be defined in sheaf theoretical terms to be sure that
  $W'\otimes_SW''$ is indeed a subbundle, since the intersection of
  $W'\otimes W''$ and the kernel of the symmetrization map might
  change dimension from one fiber to the other). Also, we denote by
  $W'^{\perp}\subset W^*$ the kernel of the restriction map $W^*\to
  W'^*$.
\end{notation}

\subsubsection{(Semi,poly)stability for $\Sp(2n,\RR)$ bundles}
\label{sss:semi-poly-stability-sp-2n-R}
Let $\alpha$ be a real number. In order to state the
$\alpha$-(semi,poly)stability condition for an $L$-twisted
$\Sp(2n,\RR)$-Higgs pair, we need to introduce some notation.

Let $k$ be an integer satisfying $k\geq 1$. We call a 
{\bf filtration
of $V$ of length $k-1$} any strictly increasing filtration by holomorphic
subbundles
$$  \VVV=(0=V_0\subsetneq V_1\subsetneq V_2\subsetneq\dots\subsetneq
V_{k-1}\subsetneq V_{k}=V).$$
Let $\lambda=(\lambda_1<\lambda_2<\dots<\lambda_{k})$ be a strictly
increasing sequence of $k$
real numbers.
Define
$$N_{\beta}(\VVV,\lambda)=
\sum_{\lambda_i+\lambda_j\leq 0 \atop 1\leq i,j\leq k} L\otimes (V_i\otimes_S V_j)
\qquad
\text{and}
\qquad
N_{\gamma}(\VVV,\lambda)=
\sum_{\lambda_i+\lambda_j\geq 0 \atop 1\leq i,j\leq k} L\otimes (V_{i-1}^{\perp}\otimes_S
V_{j-1}^{\perp}),$$
and let \footnote{This is the same as the bundle $L\otimes
E(B)_{\sigma,\chi}^-$ of
Section~\ref{ss:app-filtrations}; we use the
notation $N(\VVV,\lambda)$ for convenience.}
\begin{equation}
  \label{eq:NV-lambda-def}
N(\VVV,\lambda)=N_{\beta}(\VVV,\lambda)\oplus N_{\gamma}(\VVV,\lambda).
\end{equation}
Define also\footnote{This expression is equal to
  $\deg(E)(\sigma,\chi)- \la\alpha,\chi\ra$ of
  Section~\ref{ss:app-filtrations}.}
\begin{equation}
\label{eq:dV-lambda-alpha}
d(\VVV,\lambda,\alpha)=
\lambda_k(\deg V_k-\alpha n_k)
+\sum_{j=1}^{k-1}(\lambda_j-\lambda_{j+1})(\deg V_j-\alpha n_j),
\end{equation}
where $n_j=\rk V_j$.

We say that the pair $(\VVV,\lambda)$ is {\bf trivial} if
the length of $\VVV$ is $0$ and $\lambda_1=0$.

According to Section~\ref{ss:app-filtrations} (see also
\cite{bradlow-garcia-prada-mundet:2003}) the $\alpha$-(semi)stability
condition for an $L$-twisted $\Sp(2n,\RR)$-Higgs pair can now be
stated as follows.

\begin{proposition}
  \label{prop:sp2n-alpha-stability}
The pair $(V,\varphi)$ is $\alpha$-(semi)stable if for any integer $k\geq 1$,
any filtration $\VVV$ of length $k-1$ of $V$ and any strictly\footnote{The consideration of \emph{strictly} increasing
  sequences is justified by
  Lemma~\ref{lemma:N-d-reduced-filtration-2} below.}
 increasing sequence $\lambda$
of $k$ real numbers such that $(\VVV,\lambda)$ is not trivial and
$\varphi\in H^0(N(\VVV,\lambda))$, the inequality
\begin{equation}
  d(\VVV,\lambda,\alpha)(\geq)> 0 \label{eq:des}
\end{equation}
holds (i.e., the inequality is required to be strict when defining stability, but it
can be equality when defining semistability).

The pair $(V,\varphi)$ is $\alpha$-polystable if it is $\alpha$-semistable
and for any integer $k\geq 1$,
any filtration $\VVV$ of length $k-1$ of $V$ and any
strictly   
increasing sequence $\lambda$
of $k$ real numbers such that $(\VVV,\lambda)$ is not trivial,
$\varphi\in H^0(N(\VVV,\lambda))$ and $d(\VVV,\lambda,\alpha)=0$,
there is an isomorphism of holomorphic bundles
$$\sigma:V\to V_1\oplus V_2/V_1\oplus\dots V_k/V_{k-1}$$
such that
$V_j=\sigma^{-1}(V_1\oplus V_2/V_1\oplus\dots V_j/V_{j-1})$ for each $1\leq j\leq k$, and such that
$$\beta\in
H^0\left(\sum_{\lambda_i+\lambda_j=0} L\otimes \sigma^{-1}(V_i/V_{i-1})\otimes_S \sigma^{-1}(V_j/V_{j-1}) \right)$$
and
$$\gamma\in
H^0\left(\sum_{\lambda_i+\lambda_j=0} L\otimes \sigma^*((V_i/V_{i-1})^*)\otimes_S
\sigma^*((V_j/V_{j-1})^*)\right).$$
In the particular case when $L = K$, we obtain the (semi,poly)stability
conditions for $\Sp(2n,\RR)$-Higgs bundles by setting $\alpha = 0$ above.
\end{proposition}

The case in which $(\VVV,\lambda)$ is trivial may of course be included
in the definition of $\alpha$-semistability, since in this situation
$d(\VVV,\lambda,\alpha)=0$ automatically.

\subsubsection{When the inclusions in $\VVV$ are not strict}
\label{sss:not-strict-inclusions}

Let $k\geq 1$ be an integer, and let $\VVV$ denote a sequence of not necessarily
strict inclusions of holomorphic
subbundles
$$\VVV=(0=V_0\subset V_1\subset V_2\subset\dots\subset V_{k-1}\subset V_{k}=V).$$
Let
$\lambda=(\lambda_1<\lambda_2<\dots<\lambda_{k})$ be a strictly
increasing sequence of $k$. The fact that the inclusions $V_i\subset V_{i+1}$
may be equalities does not prevent the definitions (\ref{eq:NV-lambda-def})
and (\ref{eq:dV-lambda-alpha}) from making sense. On the other hand, we can
relate this more general situation to the one considered previously via the
following construction.

\newcommand{\prev}{\operatorname{prev}}

Define
\begin{equation}
\label{eq:definicio-I}
I=\{i_1,\dots,i_r\}=\{i\in\{1,\dots,k\}\mid V_{i-1}\neq V_i\}
\end{equation}
and let
$$(\VVV,\lambda)':=((0:=V_{i_0}\subsetneq V_{i_1}\subsetneq\dots\subsetneq V_{i_{r-1}}\subsetneq V_{i_r}=V),
(\lambda_{i_1}<\dots\lambda_{i_r})).$$
Define also, for any $1\leq k\leq r$:
$$\prev(i_k):=\left\{\begin{array}{ll} 0 & \qquad\text{if $k=1$} \\ i_{k-1} & \qquad\text{otherwise} \end{array}\right.$$

\begin{lemma}
\label{lemma:N-d-reduced-filtration}
We have
$N(\VVV,\lambda)=N((\VVV,\lambda)')$ and
$d(\VVV,\lambda,\alpha)=d((\VVV,\lambda)',\alpha)$.
\end{lemma}
\begin{proof}
The equality
$d(\VVV,\lambda,\alpha)=d((\VVV,\lambda)',\alpha)$ is obvious.
We prove separately $N_{\beta}(\VVV,\lambda)=N_{\beta}((\VVV,\lambda)')$
and $N_{\gamma}(\VVV,\lambda)=N_{\gamma}((\VVV,\lambda)')$.
Using the set $I$ defined in (\ref{eq:definicio-I}) we may write
$$N_{\beta}((\VVV,\lambda)')=
\sum_{\lambda_i+\lambda_j\leq 0 \atop i,j\in I} L\otimes V_i\otimes_S V_j
\qquad
\text{and}
\qquad
N_{\gamma}((\VVV,\lambda)')=
\sum_{\lambda_i+\lambda_j\geq 0 \atop i,j\in I} L\otimes V_{\prev(i)}^{\perp}\otimes_S
V_{\prev(j)}^{\perp},$$
The inclusion $N_{\beta}((\VVV,\lambda)')\subset N_{\beta}(\VVV,\lambda)$
is obvious. Since $V_{\prev(i)}=V_{i-1}$ for any $i$, we have
$$N_{\gamma}((\VVV,\lambda)')=
\sum_{\lambda_i+\lambda_j\geq 0 \atop i,j\in I} L\otimes V_{i-1}^{\perp}\otimes_S
V_{j-1}^{\perp},$$
from which we deduce $N_{\gamma}((\VVV,\lambda)')\subset N_{\gamma}(\VVV,\lambda)$.

To prove the inclusion $N_{\beta}(\VVV,\lambda)\subset N_{\beta}((\VVV,\lambda)')$
note that if $1\leq i,j\leq k$ satisfy $\lambda_i+\lambda_j\leq 0$ and we denote by $1\leq i'\leq k$
(resp. $1\leq j'\leq k$) the smallest integer such that $V_{i'}=V_i$ (resp. $V_{j'}=V_j$) then we
have $i',j'\in I$ and $\lambda_{i'}\leq\lambda_i$, $\lambda_{j'}\leq\lambda_j$, so
$\lambda_{i'}+\lambda_{j'}\leq 0$, and clearly $V_i\otimes_S V_j=V_{i'}\otimes_S V_{j'}$.

Similarly, $N_{\gamma}(\VVV,\lambda)\subset N_{\gamma}((\VVV,\lambda)')$ is proved by observing that
if $1\leq i,j\leq k$ satisfy $\lambda_i+\lambda_j\geq 0$ and we denote by $1\leq i'\leq k$
(resp. $1\leq j'\leq k$) the biggest integer such that $V_{i'-1}=V_{i-1}$ (resp. $V_{j'-1}=V_{j-1}$) then we
have $i',j'\in I$ and $\lambda_{i'}\geq\lambda_i$, $\lambda_{j'}\geq\lambda_j$, so
$\lambda_{i'}+\lambda_{j'}\geq 0$, and clearly $V_{i-1}^{\perp}\otimes_S V_{j-1}^{\perp}=V_{i'-1}^{\perp}\otimes_S V_{j'-1}^{\perp}$.
\end{proof}

\subsubsection{When the numbers in $\lambda$ are not strictly increasing}
\label{sss:not-strictly-increasing}

Now assume that $\VVV$ is a filtration of $V$ of length $k-1$, with strict inclusions
as in Subsection \ref{sss:semi-poly-stability-sp-2n-R}:
$$  \VVV=(0=V_0\subsetneq V_1\subsetneq V_2\subsetneq\dots\subsetneq
V_{k-1}\subsetneq V_{k}=V),$$
and that $\lambda$ is a sequence of nondecreasing real numbers
$\lambda=(\lambda_1\leq \lambda_2\leq \dots\leq \lambda_{k})$. Again,
the definitions of $N(\VVV,\lambda)$ and $\deg(\VVV,\lambda,\alpha)$ make sense
and they can be related to the case considered in Subsection \ref{sss:semi-poly-stability-sp-2n-R}.
Indeed, if we define
\begin{equation}
\label{eq:definicio-J}
J=\{j_1,\dots,j_r\}=\{i\in\{1,\dots,k-1\}\mid \lambda_{i}\neq \lambda_{i+1}\}\cup\{k\}
\end{equation}
and
$$(\VVV,\lambda)'':=((0:=V_{j_0}\subsetneq V_{j_1}\subsetneq\dots\subsetneq V_{j_{r-1}}\subsetneq V_{j_r}=V),
(\lambda_{j_1}<\dots\lambda_{j_r})),$$
then we have:

\begin{lemma}
\label{lemma:N-d-reduced-filtration-2}
$N(\VVV,\lambda)=N((\VVV,\lambda)'')$ and
$d(\VVV,\lambda,\alpha)=d((\VVV,\lambda)'',\alpha)$.
Furthermore,
$$(\VVV,\lambda)''\text{ is trivial }\qquad\Longleftrightarrow\qquad
\lambda_1=\dots=\lambda_k=0.$$
\end{lemma}
\begin{proof}
The equality
$d(\VVV,\lambda,\alpha)=d((\VVV,\lambda)'',\alpha)$ is obvious.
Define, for any $2\leq k\leq r$,
$\prev(j_k)=j_{k-1}$, and let $\prev(j_1)=0$.
We prove separately $N_{\beta}(\VVV,\lambda)=N_{\beta}((\VVV,\lambda)'')$
and $N_{\gamma}(\VVV,\lambda)=N_{\gamma}((\VVV,\lambda)'')$.
We have
$$N_{\beta}((\VVV,\lambda)'')=
\sum_{\lambda_i+\lambda_j\leq 0 \atop i,j\in J} L\otimes V_i\otimes_S V_j
\qquad
\text{and}
\qquad
N_{\gamma}((\VVV,\lambda)'')=
\sum_{\lambda_i+\lambda_j\geq 0 \atop i,j\in J} L\otimes V_{\prev(i)}^{\perp}\otimes_S
V_{\prev(j)}^{\perp}.$$
The inclusion $N_{\beta}((\VVV,\lambda)'')\subset N_{\beta}(\VVV,\lambda)$
is obvious.
To prove the inclusion $N_{\beta}(\VVV,\lambda)\subset N_{\beta}((\VVV,\lambda)'')$
note that if $1\leq i,j\leq k$ satisfy $\lambda_i+\lambda_j\leq 0$ and we denote by $1\leq i'\leq k$
(resp. $1\leq j'\leq k$) the biggest integer such that $\lambda_{i'}=\lambda_i$ (resp. $\lambda_{j'}=\lambda_j$) then we
have $\lambda_{i'}+\lambda_{j'}\leq 0$ and, since $i\leq i'$ and $j\leq j'$,
$V_i\otimes_S V_j\subset V_{i'}\otimes_S V_{j'}$.

We prove $N_{\gamma}((\VVV,\lambda)'')\subset N_{\gamma}(\VVV,\lambda)$. Let $i,j\in J$ satisfy
$\lambda_i+\lambda_j\geq 0$. Define $i'=\prev(i)+1$ and $j'=\prev(j)+1$. Then
$\lambda_{i'}+\lambda_{j'}=\lambda_i+\lambda_j\geq 0$ and
$V_{\prev(i)}^{\perp}\otimes_S
V_{\prev(j)}^{\perp}=V_{i'-1}^{\perp}\otimes_S
V_{j'-1}^{\perp}$, so the inclusion follows.

We now prove $N_{\gamma}(\VVV,\lambda)\subset N_{\gamma}((\VVV,\lambda)'')$.
Suppose that $1\leq i,j\leq k$ satisfy $\lambda_i+\lambda_j\geq 0$. Define
$i'=\min\{k\in J\mid k\geq i\}$ and
$j'=\min\{k\in J\mid k\geq j\}$. Then we have
$\lambda_{i'}+\lambda_{j'}=\lambda_i+\lambda_j\geq 0$
and, since $\prev(i')\leq i-1$ and $\prev(j')\leq j-1$, we have
an inclusion $V_{i'-1}^{\perp}\otimes_S
V_{j'-1}^{\perp}\subset V_{\prev(i)}^{\perp}\otimes_S
V_{\prev(j)}^{\perp}$.

Finally, the characterization of the case in which $(\VVV,\lambda)''$ is trivial is obvious.
\end{proof}

Of course, we could consider a more general situation by combining this subsection with
the previous one, allowing pairs $(\VVV,\lambda)$ where neither the inclusions in $\VVV$
nor the inequalities in $\lambda$ are assumed to be strict. But since we will not need
this degree of generality, we will not study this case here.

\subsubsection{Simplified (semi,poly)stability for $\Sp(2n,\RR)$ bundles}

In the next theorem we prove that the notions of (semi,poly)stability can
be simplified, in the sense that it is enough to consider filtrations of
length $\leq 2$ (cf.\ Theorem~2.8.4.13 of Schmitt~\cite{schmitt:2008}).

%

\begin{theorem}
\label{thm:simplified-conditions-sp-2n-R}
Let $(V,\varphi)$ be an $L$-twisted $\Sp(2n,\RR)$-Higgs pair.
In the following statements we consider subbundles of $V$ denoted
as $V_j$, and we denote by $n_j$ the rank of $V_j$.

\begin{enumerate}
\item The pair $(V,\varphi)$ is $\alpha$-semistable
  if and only if for any filtration of holomorphic subbundles
  $0\subset V_1\subset V_2\subset V$ such that
  \begin{equation}
\varphi = (\beta,\gamma) \in H^0(L\otimes((S^2V_2 +  V_1\otimes_S
V)\oplus (S^2V_1^{\perp} +  V_2^{\perp}\otimes_S V^*)))
\label{eq:cond-filt}
  \end{equation}
we have
$\deg V-\deg V_2-\deg V_1\geq \alpha(n-n_2-n_1)$.

\item The pair $(V,\varphi)$ is $\alpha$-stable if and only if the
following condition is satisfied.
For any filtration of holomorphic subbundles $0\subset
  V_1\subset V_2\subset V$, except the filtration $0=V_1\subset V_2=V$ such that
$$\varphi\in H^0(L\otimes((S^2V_2 + V_1\otimes_S V)\oplus
(S^2V_1^{\perp} + V_2^{\perp}\otimes_S V^*)))$$ we have
$\deg V-\deg V_2-\deg V_1> \alpha(n-n_2-n_1).$

\item The pair $(V,\varphi)$ is $\alpha$-polystable if and only if it is
semistable and: for any filtration of holomorphic subbundles $0\subset
  V_1\subset V_2\subset V$ distinct from the filtration $0=V_1\subset V_2=V$ such that
$$\varphi\in H^0(L\otimes((S^2V_2 + V_1\otimes_S V)\oplus
(S^2V_1^{\perp} + V_2^{\perp}\otimes_S V^*)))$$ and
$\deg V-\deg V_2-\deg V_1=\alpha(n-n_2-n_1)$,
there exists an isomorphism of holomorphic vector bundles
$$\sigma:V\to V_1\oplus V_2/V_1\oplus V/V_2$$
satisfying these properties:
\begin{enumerate}
\item $V_1=\sigma^{-1}(V_1)$, $V_2=\sigma^{-1}(V_1\oplus V_2/V_1)$,

\item $ \beta\in H^0(K\otimes 
(S^2(\sigma^{-1}(V_2/V_1)) \oplus
\sigma^{-1}(V_1)\otimes_S\sigma^{-1}(V/V_2)))$,

\item $ \gamma\in H^0(K\otimes
(S^2(\sigma^*(V_2/V_1)^*)\oplus
\sigma^*(V_1^*)\otimes_S\sigma^*(V/V_2)^*)))$.
\end{enumerate}

\end{enumerate}
\end{theorem}

In the statement of Theorem \ref{thm:simplified-conditions-sp-2n-R} the
inclusions in the filtration $0\subset V_1\subset V_2\subset V$ are not
necessarily strict, in contrast to the filtrations considered in the 
definition given in Proposition \ref{prop:sp2n-alpha-stability}. However,
one can still interpret the conditions appearing in Theorem \ref{thm:simplified-conditions-sp-2n-R}
as particular cases of those appearing in Proposition \ref{prop:sp2n-alpha-stability},
so what Theorem \ref{thm:simplified-conditions-sp-2n-R} says is roughly
that many of the conditions involved in the characterization of
(semi,poly)stability given by Proposition \ref{prop:sp2n-alpha-stability} are
redundant. This can be seen using Lemma \ref{lemma:N-d-reduced-filtration}:
if we define
$$\VVV=(V_1\subset V_2\subset V_3=V)\qquad\text{and}\qquad \lambda=(\lambda_1=-1,\lambda_2=0,\lambda_3=1)$$
then we have
$$N(\VVV,\lambda)=L\otimes((S^2V_2 + V_1\otimes_S V)\oplus
(S^2V_1^{\perp} + V_2^{\perp}\otimes_S V^*))$$
and
$$d(\VVV,\lambda,\alpha)=\deg V-\deg V_2-\deg V_1-\alpha(n-n_2-n_1).$$
Finally, we have
$$(\VVV,\lambda)'\text{ is trivial }\qquad \Longleftrightarrow \qquad V_1=0,\,V_2=V,$$
which explains the exception considered in the characterizations of stability and polystability.

The advantage of considering non necessarily strict inclusions is that they
allow a more compact characterization of (semi,poly)stability. If we insisted
on using strict inclusions, then we would be forced to consider separately
filtrations of length $0$, $1$ and $2$, and our simplified conditions would
become more complicated.

The preceding arguments prove that the (semi,poly)stability conditions in Proposition
\ref{prop:sp2n-alpha-stability} imply the (semi,poly)stability "simplified" conditions  
stated in Theorem \ref{thm:simplified-conditions-sp-2n-R}. It thus suffices to check the
other implication, namely that the simplified conditions appearing in Theorem \ref{thm:simplified-conditions-sp-2n-R}
imply the conditions stated in Proposition \ref{prop:sp2n-alpha-stability}.

\begin{remark}
  \label{rem:phi-zero}
When $\varphi = 0$, the
$\alpha$-semistability of $(V,\varphi)$ is equivalent to $\alpha= \mu (V)$, where
$\mu (V) = \deg V / \rk V$ is the slope of $V$, and $V$ being
(semi)stable. Let us see how to deduce this from Theorem \ref{thm:simplified-conditions-sp-2n-R}.
If $0 = V_1 = V_2$, then the condition
  \eqref{eq:cond-filt} is equivalent to $\beta = 0$ and the inequality
\begin{equation}
\label{eq:des-simples}
\deg V-\deg V_2-\deg V_1\geq \alpha(n-n_2-n_1)
\end{equation}
  reads $\deg V \geq \alpha n$.  If $V_1 = V_2 = V$,
  then \eqref{eq:cond-filt} is equivalent to $\gamma = 0$ and the
  inequality \eqref{eq:des-simples} says that $\deg V \leq \alpha n$.
  Consequently, if $\varphi = (\beta,\gamma) = 0$, then $\alpha$-semistability
  implies $\alpha  = \deg V / \rk V = \mu(V)$.
In this case, taking $V_1 = 0$ and
  $V_2 \subset V$ any subbundle, the condition (\ref{eq:des-simples})
  is equivalent to $\mu (V_2) \leq \mu (V)$, so $V$ is semistable.  Conversely,
  if $V$ is semistable and $\alpha= \mu(V)$, then the condition (\ref{eq:des-simples}) is satisfied
  for any filtration $0 \subset V_1 \subset V_2 \subset V$. Indeed, (\ref{eq:des-simples})
  is equivalent to
  $$\deg V-\alpha n\geq (\deg V_1-\alpha n_1)+(\deg V_2-\alpha n_2).$$
  But $\deg V-\alpha n=0$ because $\mu(V)=\alpha$, and since $V$ is semistable we
  have $\mu(V_j)=\deg V_j/n_j\leq \mu(V)=\alpha$, which is equivalent to
  $\deg V_j-\alpha n_j\leq 0$.

In contrast,
$(V,\varphi)$ with $\varphi=0$ will never be stable as an
$\Sp(2n,\RR)$-Higgs bundle, even when $V$ is a stable vector bundle.

\end{remark}

The proof of Theorem~\ref{thm:simplified-conditions-sp-2n-R} will
be given in Section~\ref{sec:proof-thm:stabilitynotion}.

\subsection{Some results on convex sets}
\label{ss:convex-sets}
Let $W$ be an $n$ dimensional vector space over $\RR$. We denote
the convex hull of any subset $S\subset W$ by 
$$\CH(S)\subset W.$$
Let $h_1,h_2,\dots,h_l$ be elements of the dual space $W^*$. We
assume that $l\geq n$ and that the first $n$ elements
$h_1,\dots,h_n$ are a basis of $W^*$. Define for any $h\in W^*$
the set
$$\{h\leq a\}=\{v\in W\mid h(v)\leq a\}\subset W,$$
and define $\{h=a\}\subset W$ similarly.

Consider the convex subset of $W$
$$C=\bigcap_i\{h_i\leq 0\}$$
(here and below if no range is specified for the index then it is
supposed to be the whole set $\{1,\dots,l\}$).

\begin{remark}\label{rem:subspaces}
  The fact that $\{h_1,\dots,h_l\}$ span $W^*$ is equivalent to the
  condition that $C$ does not contain any positive dimensional vector
  subspace of $W$.  Indeed, if $h \in W^*$ and $Z \subset W$ is a
  subspace contained in $\{h \leq 0\}$, then $Z$ is contained in $\{h
  = 0\}$. Consequently any vector subspace of $W$ contained in $C$ has
  to lie in $\bigcap_i\{h_i = 0\} = 0$.
\end{remark}

Given a subset $A\subset W$ we denote by $W_A\subset W$ the smallest affine 
subspace containing $A$, and we denote by
$$\partial_{\RR}A\subset W_A\subset W$$
the boundary (in the sense of set topology) of $A$ {\it viewed as a subspace of $W_A$}.
More concretely, $\partial_{\RR}A$ is the set of points which are both the limit of a
sequence of points in $A$ and of a sequence of points in $W_A\setminus A$.
The proof of the following statement is straightforward and left to the reader:
\begin{equation}
\label{eq:frontera-compacte}
A\subset W\text{ compact and convex}\qquad\Longrightarrow\qquad A=\CH(\partial_{\RR}A).
\end{equation}
The following is also easy to prove:
\begin{equation}
\label{eq:frontera-interseccio}
W'\subset W\text{ affine subspace}\qquad\Longrightarrow\qquad \partial_{\RR}(A\cap W')\subset\partial_{\RR} A
\quad \text{for any $A\subset W$}.
\end{equation}

\begin{lemma}
$C=\CH(\partial_{\RR} C)$.
\end{lemma}
\begin{proof}
For any $\alpha\leq 0$ define
$C_\alpha=C\cap\{h_1+\dots+h_n=\alpha\}$. Since $C$ is closed,
for any $x\in C$
we have $h_i(x)\leq 0$, and furthermore $h_1,\dots,h_n$ is a basis
of $W^*$, we deduce that $C_\alpha$ is compact. Hence, by (\ref{eq:frontera-compacte}),
$C_\alpha=\CH(\partial_{\RR} C_\alpha)$. Now take any $x\in C$ and set
$\alpha=h_1(x)+\dots+h_n(x)$. Then $x\in C_\alpha=\CH(\partial
C_\alpha)\subset \CH(\partial_{\RR} C)$, where the second inclusion
follows from (\ref{eq:frontera-interseccio}). This proves that
$C\subset\CH(\partial_{\RR} C)$. The other inclusion follows from the
fact that $C$ is convex and closed.
\end{proof}

It follows from the definition of $C$ that we have 
$$\partial_{\RR} C=\bigcup_i C_i,$$ 
where $C_i=\{h_i=0\}\cap
C.$ For any $i$ the collection of elements
$h_1,\dots,h_l$ induce elements $h'_1,\dots,h'_l$ on the dual of
$\{h_i=0\}$ which obviously span. Hence we may apply again the
lemma to $C_i$ and deduce that $C_i=\CH(\partial_{\RR} C_i)$. Proceeding
recursively, we deduce that $C$ is the convex hull of the union of
the sets
$$C_I=\bigcap_{i\in I} \{h_i=0\}\cap C$$
where $I$ runs over the collection of subsets of $\{1,\dots,l\}$
satisfying:
\begin{equation}
\label{propietatI} \text{$|I|=n-1$ and the vectors $\{ h_i\mid
i\in I\}$ are linearly independent.}
\end{equation}

Each such subset $C_I$ is a halfline.

\begin{lemma}
Fix a basis $e_1,\dots,e_n$ of $W$, and denote by
$e_1^*,\dots,e_n^*$ the dual basis. Assume that any $h_i$ can be
written either as $e_a^*-e_b^*$ or $\pm(e_a^*+e_b^*)$ for some
indices $a,b$ depending on $i$. Then for any $I$ satisfying
(\ref{propietatI}) there are disjoint subsets
$P,N\subset\{1,\dots,n\}$ so that defining the element
$c_I=\sum_{i\in P} e_i-\sum_{j\in N} e_j$ we have $C_I=\RR_{\geq
0}c_I$. \label{lemma:eixos}
\end{lemma}
\begin{proof}
Pick some $I$ satisfying (\ref{propietatI}), so that
$C_I=\bigcap_{i\in I}\{h_i=0\}$ is one dimensional, and let
$c_I\in W$ be an element such that $C_I=\RR_{\geq 0}c_I$. Write
$c_I=\sum\lambda_j e_j$ and take some nonzero
$\lambda\in\{\lambda_1,\dots,\lambda_n\}$. Define
$P_{\lambda}=\{j\mid \lambda_j=\lambda\}$ and $N_{\lambda}=\{j\mid
\lambda_j=-\lambda\}$. We want to prove that for any $j\notin
P_{\lambda}\cup N_{\lambda}$, $\lambda_j=0$. Suppose the contrary.
Then
$$c_I'=\sum_{j\in P_{\lambda}\cup N_{\lambda}}2\lambda_je_j
+\sum_{j\notin P_{\lambda}\cup N_{\lambda}}\lambda_je_j$$ does not
belong to $\RR c_I$. However, for any pair of indices $a,b$ we
clearly have
$$(e_a^*-e_b^*)c_I=0\Longrightarrow (e_a^*-e_b^*)c'_I=0
\qquad\text{ and }\qquad (e_a^*+e_b^*)c_I=0\Longrightarrow
(e_a^*+e_b^*)c'_I=0.$$ This implies by our assumption that
$c_I'\in \bigcap_{i\in I}\{h_i=0\}=C_I$, in contradiction with the
fact that $C_I$ is one dimensional.
\end{proof}

\subsection{Proof of Theorem~\ref{thm:simplified-conditions-sp-2n-R}}
\label{sec:proof-thm:stabilitynotion}
As already mentioned, when $\varphi = 0$ the pair $(V,0)$ is
$\alpha$-semistable if and only if $\alpha = \mu (V)$ and $V$ is
semistable. Thus, by Remark~\ref{rem:phi-zero}, it suffices to
consider the case $\varphi \neq 0$.  Moreover, by
Lemmas~\ref{lemma:N-d-reduced-filtration} and
\ref{lemma:N-d-reduced-filtration-2} we may consider non-strictly
increasing filtrations and sequences $\lambda$ in the
(semi-,poly-)stability condition.

Let $\VVV$ be any filtration of
$V$, and define
$$\Lambda(\VVV,\varphi)=\{\lambda\in\RR^k\mid \lambda_1\leq\dots\leq\lambda_k,
\ \varphi\in N(\VVV,\lambda)\}.$$ The pair $(V,\varphi)$ is
$\alpha$-semistable if for any $\lambda\in \Lambda(\VVV,\varphi)$
we have
$$d(\VVV,\lambda,\alpha)\geq 0.$$
But $d(\VVV,\lambda,\alpha)$ is clearly a linear function on
$\lambda$, so to check semistability it suffices to verify that
$d(\VVV,\lambda,\alpha)\geq 0$ for any $\lambda$ belonging to a
set $\Lambda'\subset\RR^k$ whose convex hull is
$\Lambda(\VVV,\varphi)$. Define for any $i,j$ the subbundles
$$D_{i,j}=V_i\otimes_S V_j+V_{i-1}\otimes_S V+V\otimes_S V_{j-1}\subset S^2V$$
and
$$D^*_{i,j}=V_{i-1}^{\perp}\otimes_S V_{j-1}^{\perp}+
V_{i}^{\perp}\otimes_S V^*+V^*\otimes_S V_{j}^{\perp}\subset
S^2V^*.$$
A tuple $\lambda_1\leq\dots\leq\lambda_k$ belongs to
$\Lambda(\VVV,\varphi)$ if and only if these two conditions holds:
\begin{itemize}
\item for any $i,j$ such that $\beta$ is contained in
$H^0(L\otimes D_{i,j})$ but is not contained in the sum
$H^0(L\otimes D_{i-1,j})+H^0(L\otimes D_{i,j-1})$, we have
 $\lambda_i+\lambda_j\leq 0.$
\item for any $i,j$ such that $\gamma$ is contained in
$H^0(L\otimes D_{i,j}^*)$ but is not contained in the sum
$H^0(L\otimes D_{i+1,j}^*)+H^0(L\otimes D_{i,j+1}^*)$, we have
 $\lambda_i+\lambda_j\geq 0.$
\end{itemize}

Hence $\Lambda(\VVV,\varphi)\subset\RR^k$ is the intersection of
halfspaces of the form $\{\lambda_i-\lambda_{i+1}\leq 0\}$ and,
$\{\lambda_a+\lambda_b\leq 0\}$ (for at least one pair $(a,b)$, if
$\beta\neq 0$) or $\{\lambda_c+\lambda_d\geq 0\}$ (for at least one
pair $(c,d)$, if $\gamma \neq 0$). Since the only nonzero vector
subspace included in the set
$\Lambda=\{\lambda_1\leq\dots\leq\lambda_k\}$ is the line generated by
$(1,\dots,1)$ and the set $\Lambda(\VVV,\varphi)$ is contained in
$\Lambda$ and furthermore satisfies at least one equation of the form
$\lambda_a+\lambda_b\geq 0$ or $\lambda_c+\lambda_d\leq 0$, it follows
that $\Lambda(\VVV,\varphi)$ does not contain any nonzero vector
subspace.

So by the arguments Section~\ref{ss:convex-sets}, the space
$\Lambda(\VVV,\varphi)$ is the convex hull of a collection of half
lines of the form $\RR_{\geq 0} \lambda_I$, and by Lemma
\ref{lemma:eixos} we can assume that the coordinates of
$\lambda_I$ are $0$ and $\pm 1$. But if
$\lambda_I\in\Lambda(\VVV,\varphi)$ we necessarily must have
$c_I=(-1,\dots,-1,0,\dots,0,1,\dots,1)$, say $a$ copies of $-1$,
$b$ of $0$ and $k-(a+b)$ of $1$. Consider first the case when $0
<a < a+b < k$. Define now the filtration
$$\VVV'=(0\subsetneq V_{a}\subsetneq V_{a+b}\subsetneq V).$$
One can easily check that
$$d(\VVV,\lambda_I,\alpha)=d(\VVV',(-1,0,1),\alpha)
=\deg V-\deg V_a-\deg V_{a+b}-\alpha(n-n_a-n_{a+b}),$$ and that
$N(\VVV,\lambda) = L\otimes((S^2V_{a+b}+ V_a\otimes_S V)\oplus
(S^2V_a^{\perp}+ V_{a+b}^{\perp}\otimes_S V^*))$.

It therefore follows that to check semistability, it suffices to
consider filtrations of length at most two, with $\lambda=(-1,0,1)$
(or a suitable degeneration). This concludes the proof of (1) of
Theorem~\ref{thm:simplified-conditions-sp-2n-R}. The proofs of (2) and
(3) follow in essentially the same way.
\qed

\subsection{Jordan--H\"older
reduction for polystable  $\Sp(2n,\RR)$-Higgs bundles}
\label{sec:polystability}

It follows from Section~\ref{ss:Jordan-Holder} that any
$\alpha$-polystable $G$-Higgs pair admits a Jordan--H\"older
reduction.  In order to state this result in the case of
$G=\Sp(2n,\RR)$, we need to describe some special $\Sp(2n,\RR)$-Higgs
bundles arising from $G$-Higgs bundles associated to certain real
subgroups $G \subseteq \Sp(2n,\RR)$.

\subsubsection*{The subgroup $G = \U(n)$}

Observe that a $\U(n)$-Higgs bundle is nothing but a holomorphic
vector bundle $V$ of rank $n$. The standard inclusion
$\upsilon^{\U(n)}\colon\U(n) \into \Sp(2n,\RR)$ gives the
correspondence
\begin{equation}
  \label{eq:un-sp2n}
  V \mapsto \upsilon^{\U(n)}_*V = (V,0)
\end{equation}
associating the $\Sp(2n,\RR)$-Higgs bundle $\upsilon^{\U(n)}_*V =
(V,0)$ to the holomorphic vector bundle $V$.

\subsubsection*{The subgroup $G = \U(p,q)$}

In the following we assume that $p,q \geq 1$. As is easily seen, a
$\U(p,q)$-Higgs bundle (cf.\
\cite{bradlow-garcia-prada-gothen:2003}) is given by the data
$(\tilde{V},\tilde{W},\tilde\varphi=\tilde\beta+\tilde\gamma)$,
where $\tilde{V}$ and $\tilde{W}$ are holomorphic vector bundles
of rank $p$ and $q$, respectively, $\tilde\beta\in H^0(K\otimes
\Hom(\tilde{W},\tilde{V}))$ and $\tilde\gamma\in H^0(K\otimes
\Hom(\tilde{V},\tilde{W}))$. Let $n=p+q$. The imaginary part of
the standard indefinite Hermitian metric of signature $(p,q)$ on
$\CC^n$ is a symplectic form, and thus there is an inclusion
$\upsilon^{\U(p,q)}\colon\U(p,q) \into \Sp(2n,\RR)$. At the level
of $G$-Higgs bundles, this gives rise to the correspondence
\begin{equation}
\label{eq:upq-sp2n}
(\tilde{V},\tilde{W},\tilde\varphi=\tilde\beta+\tilde\gamma)
\mapsto \upsilon^{\U(p,q)}_*(\tilde{V},\tilde{W},\tilde\varphi)
= (V, \varphi=\beta+\gamma),
\end{equation}
where
\begin{displaymath}
    V = \tilde{V} \oplus \tilde{W}^*,\quad
    \beta =
    \begin{pmatrix}
      0 & \tilde{\beta} \\
      \tilde\beta & 0
    \end{pmatrix}\quad\text{and}\quad
    \gamma =
    \begin{pmatrix}
      0 & \tilde\gamma \\
      \tilde\gamma & 0
    \end{pmatrix}.
\end{displaymath}

In the following we shall occasionally slightly abuse language, saying
simply that $\upsilon^{\U(n)}_*V$ is a $\U(n)$-Higgs bundle and that
$\upsilon^{\U(p,q)}_*(\tilde{V},\tilde{W},\tilde\varphi)$ is a
$\U(p,q)$-Higgs bundle.

Another piece of convenient notation is the following. Let
$(V_i,\varphi_i)$ be $\Sp(2n_i,\RR)$-Higgs bundles and let $n = \sum
n_i$. We can define an $\Sp(2n,\RR)$-Higgs bundle $(V,\varphi)$ by
setting
\begin{displaymath}
  V = \bigoplus V_i \quad\text{and}\quad
  \varphi = \sum \varphi_i
\end{displaymath}
by using the canonical inclusions $H^0(K \otimes (S^2V_i \oplus
S^2V_i^*)) \subset H^0(K \otimes (S^2V \oplus S^2V^*))$. We shall
slightly abuse language and write $(V,\varphi) = \bigoplus
(V_i,\varphi_i)$, referring to this as \textbf{the direct sum} of the
$(V_i,\varphi_i)$.

With all this understood, we can state our
structure theorem on polystable $\Sp(2n,\RR)$-Higgs bundles from
Section~\ref{ss:Jordan-Holder} as follows.

\begin{proposition}
  \label{prop:polystable-spnr-higgs}
  Let $(V,\varphi)$ be a polystable $\Sp(2n,\RR)$-Higgs bundle. Then
  there is a decomposition
  $$(
  V,\varphi) = (V_1,\varphi_1) \oplus \dots \oplus (V_k,\varphi_k),
  $$ unique up to reordering, such that each
  $(V_i,\varphi_i)$ is a stable $G_i$-Higgs bundle, where $G_i$ is one
  of the following groups: $\Sp(2n_i,\RR)$ (with $\varphi\neq 0$), $\U(n_i)$ or
  $\U(p_i,q_i)$.
\end{proposition}

\subsection{$\GL(n,\RR)$-Higgs bundles}
\label{sec:gln-higgs}

Throughout this paper we have worked with $G$-Higgs pairs for
connected $G$. Thus the theory developed does not apply directly to
the disconnected group $G = \GL(n,\RR)$. However, the definitions
$G$-Higgs pair and stability makes perfect sense for $G=\GL(n,\RR)$
and, with this stability condition, we certainly expect the
Hitchin--Kobayashi correspondence to be valid. However, for
disconnected groups in general, there are some rather subtle questions
involved in developing the Hitchin--Kobayashi correspondence (in
particular in identifying the kind of parabolic subgroups which should
be considered) and we hope to come back to this question in the
future. For the present paper we shall limit ourselves to analysing
the natural stability condition for $L$-twisted $\GL(n,\RR)$-Higgs
pairs. We study this case because, when $L=K^2$, it is important in
the study of maximal degree $\Sp(2n,\RR)$-Higgs bundles (see
\cite{garcia-prada-gothen-mundet:2009b}).

A maximal compact subgroup of $\GL(n,\RR)$ is $H=\O(n)$ and hence
$H^\CC=\O(n,\CC)$. Now, if  $\WW$ is the standard $n$-dimensional
complex vector  space representation of $\O(n,\CC)$, then the
isotropy representation space is:
$$
\lie{m}^\CC=S^2\WW.
$$
An $L$-twisted $\GL(n,\RR)$-Higgs pair over $X$  is thus a pair
$((W,Q),\psi)$ consisting of a holomorphic $\O(n,\CC)$-bundle,
i.e.\ a rank $n$ holomorphic vector bundle $W$ over $X$ equipped
with a non-degenerate quadratic form $Q$, and a section
$$
\psi \in H^0(L\otimes S^2W).
$$
An equivalent definition is the following. Denote by $S^{2}_{Q}W$ the bundle
of endomorphisms $\xi$ of $W$ which are symmetric with respect to $Q$ i.e.\ 
such that $Q(\xi\,\cdot,\cdot)=Q(\cdot,\xi\,\cdot)$.
An $L$-twisted $\GL(n,\RR)$-Higgs pair over $X$  is thus a triple
$(W,Q,\psi)$ consisting of a holomorphic $\O(n,\CC)$-bundle,
i.e.\ a rank $n$ holomorphic vector bundle $W$ over $X$ equipped
with a non-degenerate quadratic form $Q$, and a section
$$
\psi \in H^0(L\otimes S^2_{Q}W).
$$
Note that when $\psi = 0$ a twisted $\GL(n,\RR)$-Higgs  pair is
simply  an orthogonal bundle.

\begin{remark}
Since the center of $\olie(n)$ is trivial, $\alpha=0$ is the only
possible value for which stability of an $L$-twisted
$\GL(n,\RR)$-Higgs pair  is defined.
\end{remark}

In order to state the stability condition for twisted
$\GL(n,\RR)$-Higgs pairs, we first introduce some notation. For any
filtration of vector bundles
\begin{displaymath}
  \WWW=(0=W_0\subsetneq W_1\subsetneq W_2\subsetneq\dots\subsetneq W_k=W)
\end{displaymath}
of satisfying $W_j = W_{k-j}^{\perp_Q}$ (here $W_{k-j}^{\perp_Q}$
denotes the orthogonal complement of $W_{k-j}$ with respect to
$Q$) define
$$\Lambda(\WWW)=\{(\lambda_1,\lambda_2,\dots,\lambda_k)\in\RR^k
\mid \lambda_i\leq \lambda_{i+1}\text{ and
}\lambda_i+\lambda_{k-i+1}=0\text{ for any $i$ }\}.$$ Define for
any $\lambda\in\Lambda(\WWW)$ the following bundle.
\begin{displaymath}
    N(\WWW,\lambda)=
    \sum_{\lambda_i+\lambda_j\leq 0} L\otimes W_i\otimes_S W_j.
\end{displaymath}
Also  we define
\begin{displaymath}
  d(\WWW,\lambda) = \sum_{j=1}^{k-1}(\lambda_j-\lambda_{j+1})\deg W_j
\end{displaymath}
(note that the quadratic form $Q$ induces an isomorphism $W\simeq
W^*$ so $\deg W=\deg W_k=0$).

According to Section~\ref{ss:app-filtrations} (see also
\cite{bradlow-garcia-prada-mundet:2003}) the stability conditions (for
$\alpha=0$) for an $L$-twisted $\GL(n,\RR)$-Higgs pair can now be
stated as follows.
\begin{proposition}
  \label{prop:gln-alpha-stability}
an $L$-twisted $\GL(n,\RR)$-Higgs pair
$(W,Q,\psi)$ is   semistable if for
  all filtrations $\WWW$ as above and all
  $\lambda\in\Lambda(\WWW)$ such that
  $\psi \in H^0(N(\WWW,\lambda))$, we have
  $d(\WWW,\lambda) \geq 0 $.

  The pair $(W,\psi)$ is stable if it is semistable and for any choice of
  filtration $\WWW$ and nonzero $\lambda\in\Lambda(\WWW)$ such that
  $\psi \in H^0(N(\WWW,\lambda))$, we have $d(\WWW,\lambda) > 0 $.

 The pair $(W,\psi)$ is polystable if it is
 semistable and for any filtration $\WWW$ as above and
 $\lambda\in\Lambda(\WWW)$ satisfying $\lambda_i<\lambda_{i+1}$
 for each $i$, $\psi\in H^0(N(\WWW,\lambda))$ and
 $d(\WWW,\lambda)=0$, there is an isomorphism
 $$W\simeq W_1\oplus W_2/W_1\oplus\dots\oplus W_k/W_{k-1}$$
 such that pairing via $Q$ any element of the summand $W_i/W_{i-1}$ with
 an element of the summand $W_j/W_{j-1}$ is zero unless $i+j=k+1$;
 furthermore, via this isomorphism,
 $$\psi\in H^0(\bigoplus_{\lambda_i+\lambda_j=0} L\otimes
 (W_i/W_{i-1})\otimes_S (W_j/W_{j-1})).$$
\end{proposition}

There is a simplification of the stability condition for $L$-twisted
$\GL(n,\RR)$-Higgs pairs analogous to
Theorem~\ref{thm:simplified-conditions-sp-2n-R}.

\begin{theorem}\label{thm:orthogonal-stability}
  Let $(W,Q,\psi)$ be a $L$-twisted  $\GL(n,\RR)$-Higgs pair. Then 
$(W,Q,\psi)$ is 
{semistable} if and only if $  \deg(W')\leq 0$ for any isotropic and 
$\psi$-invariant subbundle $W'\subset W$. Furthermore,
$(W,Q,\psi)$ is  
{stable} if and only if it is semistable and $\deg(W')<0$ for any isotropic and 
$\psi$-invariant strict subbundle $0\neq W' \subset W$.
Finally, $(W,Q,\psi)$ is  
 {polystable} if and only if it is semistable and, for any isotropic 
 (respectively.\ coisotropic) and $\psi$-invariant strict
  subbundle $0\neq W'\subset W$ such that $\deg(W')=0$, there is another
  coisotropic (resp.\ isotropic) and $\psi$-invariant subbundle 
$0\neq W''\subset W$ such that $W\simeq W'\oplus W''$.
\end{theorem}

\begin{proof}
  The proof is analogous to the proofs of
  Theorem~\ref{thm:simplified-conditions-sp-2n-R}. Take an $L$-twisted
  $\GL(n,\RR)$-Higgs pair $((W,Q),\psi)$, and assume that for any
  isotropic subbundle $W' \subset W$ such that $\psi \in H^0(S^2
  W'^{\perp_Q} \oplus W' \otimes_S W \otimes L)$ the inequality $ \deg
  W' \leq 0$ holds. We also assume that $\psi$ is nonzero, for
  otherwise the result follows from the usual characterization of
  (semi)stability for $\SO(n,\CC)$-principal bundles due to Ramanathan
  (see \cite{ramanathan:1975}). We want to prove that $((W,Q),\psi)$
  is semistable. Choose any filtration $\WWW=(0\subsetneq
  W_1\subsetneq W_2\subsetneq\dots\subsetneq W_k=W)$ satisfying
  $W_{k-i}=W_i^{\perp_{\Omega}}$ for any $i$. Consider the convex set
$$\Lambda(\WWW,\psi)= \{\lambda\in\Lambda(\WWW)\mid \psi\in
N(\WWW,\lambda)\}\subset\RR^k.$$ Define for any $i,j$ the
subbundle
$$D_{i,j}=W_i\otimes_S W_j+W_{i-1}\otimes_S W+W\otimes_S W_{j-1}\subset S^2W.$$
A tuple $\lambda=(\lambda_1,\dots,\lambda_k)\in\Lambda(\WWW)$
belongs to $\Lambda(\WWW,\psi)$ if and only if:
\begin{quote}
for any $i,j$ such that $\psi$ is contained in $H^0(L\otimes
D_{i,j})$ but is not contained in the sum $H^0(L\otimes
D_{i-1,j})+H^0(L\otimes D_{i,j-1})$, we have
 $\lambda_i+\lambda_j\leq 0.$
\end{quote}
Hence $\Lambda(\WWW,\psi)$ is the intersection of $\Lambda(\WWW)$
with the set of points in $\RR^k$ satisfying a collection of
inequalities of the form $\lambda_a+\lambda_b\leq 0$ and
$\lambda_c+\lambda_d\geq 0$ (the latter follow from the
restrictions $\lambda_i+\lambda_{k-i+1}=0$). Since $\Lambda(\WWW)$
does not contain any line, a fortiori $\Lambda(\WWW,\psi)$ neither
does, so (using Lemma \ref{lemma:eixos}) $\Lambda(\WWW,\psi)$ is
the convex hull of a set of half lines $\{\RR_{\geq 0}L_i\mid
i\in\III\}$, where $L_i=(-1,\dots,-1,0,\dots,0,1,\dots,1)$
contains $i$ copies of $-1$ and $i$ copies of $1$. Consequently,
we have
$$d(\WWW,\lambda)\geq 0\text{ for any $\lambda\in\Lambda(\WWW,\psi)$ }
\Longleftrightarrow d(\WWW,L_i)\geq 0\text{ for any $i\in\III$
}.$$ It follows from the definition that $N(\WWW,L_i)=W_i\otimes_S
W+ S^2W_{k-i}$ and since $W_{k-i}=W_i^{\perp_Q}$ the condition
$L_i\in \Lambda(\WWW,\psi)$ can be translated into the condition
$$\psi \in H^0(S^2 W_i^{\perp_Q}  \oplus W_i \otimes_S W \otimes L).$$
One computes $d(\WWW,L_i)=-\deg W_{k-i}-\deg W_i$. On the other
hand, since we have an exact sequence $0\to W_{k-i}\to W^*\to
W_i*\to 0$ (the injective arrow is given by the pairing with the
quadratic form $Q$) we have $0=\deg W^*=\deg W_{k-i}+\deg W_i^*$,
so $\deg W_{k-i}=\deg W_i$ and consequently $d(\WWW,L_i)\geq 0$ is
equivalent to $\deg W_i\leq 0$, which holds by assumption. Hence
$((W,Q),\psi)$ is semistable.

The converse statement, namely, that if $((W,Q),\psi)$ is
semistable then for any isotropic subbundle $W'\subset W$ such
that $\Phi(W')\subset L\otimes W'$ we have $\deg W'\leq 0$ is
immediate by applying the stability condition of the filtration
$0\subset W'\subset W'^{\perp_{Q}}\subset W$.

Finally, the proof of the second statement on stability is very
similar to the  case of semistability, so we omit it. The
statement on polystability is also straightforward.
\end{proof}

\begin{remark}
  The condition $\psi \in H^0(S^2 W_1^{\perp_Q}
  \oplus W_1 \otimes_S W \otimes L)$ is equivalent to
  $\tilde{\psi}(W_1) \subseteq W_1 \otimes L$, where
  $\tilde{\psi}=\psi\circ Q \colon W \to W\otimes L$.
\end{remark}


\begin{thebibliography}{99}

\bibitem{banfield}
D. Banfield, \emph{Stable pairs and principal bundles},
Quart. J. Math. \textbf{51} (2000), 417--436.


\bibitem{biswas-ramanan:1994}
I.~Biswas and S.~Ramanan, \emph{An infinitesimal study of the moduli of
 {H}itchin pairs}, J. London Math. Soc. (2) \textbf{49} (1994), 219--231.


\bibitem{borel:1956} A. Borel, \emph{Groupes linéaires algébriques},
  Ann. of Math. (2) {\bf 64} (1956), 20--82.

\bibitem{borel:1991} \bysame, \emph{Linear algebraic groups}, Second
  edition, Graduate Texts in Mathematics, vol.~126, Springer--Verlag,
  New York, 1991.


\bibitem{bradlow-garcia-prada-gothen:2003}
S.~B. Bradlow, O. Garc{\'\i}a-Prada, and P~.B. Gothen, \emph{Surface group representations and $\mathrm{U}(p,q)$-{H}iggs
  bundles}, J. Differential Geom. \textbf{64} (2003), 111--170.

\bibitem{bradlow-garcia-prada-mundet:2003}
S.~B. Bradlow, O.~Garc{\'{\i}}a-Prada, and I.~Mundet~i Riera, \emph{Relative
  {H}itchin-{K}obayashi correspondences for principal pairs}, Quart. J. Math.
  \textbf{54} (2003), 171--208.


\bibitem{corlette:1988}
K.~Corlette, \emph{Flat ${G}$-bundles with canonical metrics}, J. Differential
  Geom. \textbf{28} (1988), 361--382.


\bibitem{donaldson:1987}
S.~K. Donaldson, \emph{Twisted harmonic maps and the self-duality equations},
  Proc. London Math. Soc. (3) \textbf{55} (1987), 127--131.


\bibitem{donaldson-kronheimer:1990}
S.~K. Donaldson and P.~B. Kronheimer, \emph{The geometry of four-manifolds},
  Oxford Mathematical Monographs, The Clarendon Press Oxford University Press,
  New York, 1990.

\bibitem{garcia-prada-gothen-mundet:2009b}
O.~Garc{\'\i}a-Prada, P.~B. Gothen, and I.~Mundet i~Riera,
  \emph{Higgs bundles and surface group representations in the
real symplectic group}, v4: August 2012,
  \texttt{arXiv:0809.0576 [math.AG]}.

\bibitem{goldman:1984}
W.~M. Goldman, \emph{The symplectic nature of fundamental groups of
surfaces},
Adv. Math. \textbf{54} (1984), No. 2,  200--225.

\bibitem{hitchin:1987a}
N.~J. Hitchin, \emph{The self-duality equations on a {R}iemann surface}, Proc.
  London Math. Soc. (3) \textbf{55} (1987), 59--126.

\bibitem{huybrechts:2006} D.~Huybrechts, \emph{Fourier--Mukai
transforms in algebraic geometry}, Oxford University Press, 2006.

\bibitem{kobayashi:1987}
S. Kobayashi, \emph{Differential Geometry of Complex Vector Bundles},
Princeton University Press, 1987.

\bibitem{mundet:2000}
I.~Mundet i Riera, \emph{A Hitchin--Kobayashi correspondence for K\"ahler
fibrations}, J. Reine Angew. Math. {\bf 528} (2000), 41--80.

\bibitem{narasimhan-seshadri:1965}
M.~S. Narasimhan and C.~S. Seshadri, \emph{Stable and unitary vector bundles on
  a compact {R}iemann surface}, Ann. of Math.(2) \textbf{82} (1965), 540--567.

\bibitem{ramanathan:1975}
A.~Ramanathan, \emph{Stable principal bundles on a compact {R}iemann surface},
  Math. Ann. \textbf{213} (1975), 129--152.

\bibitem{ramanathan:1996} \bysame, \emph{Moduli for principal
    bundles over algebraic curves: I and II}, Proc. Indian
  Acad. Sci. Math. Sci.  \textbf{106} (1996), 301--328 and 421--449.

 \bibitem{Ri}
 R.W. Richardson, Conjugacy classes of $n$-tuples in Lie algebras
 and algebraic groups, {\em Duke Math. J.} {\bf 57} (1988) 1--35.

\bibitem{schmitt:2005}
A.~H.~W. Schmitt, \emph{Moduli for decorated tuples for sheaves and
representation spaces for quivers},   Proc. Indian Acad. Sci.
Math. Sci.  \textbf{115} (2005), 15--49.

\bibitem{schmitt:2008} \bysame, \emph{Geometric invariant theory
    and decorated principal bundles}, Z\"urich Lectures in Advanced
  Mathematics, European Mathematical Society, 2008.

\bibitem{serre:1987} J.~P. Serre, {\em Complex Semisimple Lie
Algebras}, Springer--Verlag, 1987.

\bibitem{simpson:1988}
C.~T. Simpson, \emph{Constructing variations of {H}odge structure using
  {Y}ang-{M}ills theory and applications to uniformization}, J. Amer. Math.
  Soc. \textbf{1} (1988), 867--918.

\bibitem{simpson:1992}
\bysame, \emph{Higgs bundles and local systems}, Inst. Hautes {\'E}tudes Sci.
  Publ. Math. \textbf{75} (1992), 5--95.

\bibitem{simpson:1994}
\bysame, \emph{Moduli of representations of the fundamental group of a smooth
  projective variety {I}}, Publ. Math., Inst. Hautes \'Etud. Sci. \textbf{79}
  (1994), 47--129.

\bibitem{simpson:1995}
\bysame, \emph{Moduli of representations of the fundamental group of a smooth
  projective variety {II}}, Publ. Math., Inst. Hautes \'Etud. Sci. \textbf{80}
  (1995), 5--79.

\bibitem{sjamaar:1995} R. Sjamaar,
Holomorphic slices, symplectic reduction and multiplicities of
representations, {\em Ann. of Math.} (2) {\bf 141} (1995),
87--129.

\bibitem{uhlenbeck-yau:1988}
K. Uhlenbeck, S.T. Yau, \emph{On the existence of Hermitian-Yang
Mills connections in stable vector bundles.} Frontiers of the
mathematical sciences: 1985 (New York, 1985). Comm. Pure Appl.
Math. \textbf{39} (1986), no. S, suppl., S257--S293.

\end{thebibliography}
\end{document}